\documentclass[preprint,3p,10pt]{elsarticle}
\usepackage{graphicx}
\usepackage{amsmath, amssymb, amsthm}
\usepackage{subfig}
\usepackage{xcolor}
\usepackage{mathrsfs}
\usepackage{multirow}
\usepackage{float}
\usepackage{bm}
\usepackage{mathrsfs}
\usepackage{booktabs}
\usepackage{appendix}
\usepackage{grffile}
\usepackage{algorithm}
\usepackage{algorithmic}
\usepackage{cases}
\usepackage[colorlinks,linkcolor=blue,anchorcolor=blue,citecolor=red]{hyperref}

\numberwithin{equation}{section}

\newcommand{\dd}{\mathop{}\!\mathrm{d}}
\newcommand{\pp}{\partial}

\newcommand{\fpp}[2]{\frac{\pp #1}{\pp #2}}
\newcommand{\mF}{\mathcal{F}}

\newcommand{\bx}{{\boldsymbol{x}}}
\newcommand{\bu}{{\boldsymbol{u}}}
\newcommand{\bbR}{{\mathbb{R}}}
\newcommand{\Pn}{P_N}
\newcommand{\bA}{{\bf A}}
\newcommand{\bB}{{\bf B}}

\newcommand{\mO}{\mathcal{O}}
\newcommand{\mE}{\mathcal{E}}

\newcommand{\bbu}{\bar{\bm u}}

\newtheorem{remark}{Remark}

\usepackage{ulem}

\begin{document}

\title{
An asymptotic-preserving  method for the three-temperature radiative transfer model}


\author{Ruo Li\fnref{fn1}}
\fntext[fn1]{CAPT, LMAM \& School of Mathematical Sciences,
    Peking University, Beijing, China, email: {\tt
      rli@math.pku.edu.cn}.}

\author{Weiming Li\fnref{fn2}}
\fntext[fn2]{Institute of Applied Physics and Computational 
Mathematics, Beijing, China, email:{\tt li\_weiming@iapcm.ac.cn}.}

\author{Shengtong Liang\fnref{fn3}}
\fntext[fn3]{School of Mathematical Sciences, Peking University, Beijing, China, email: {\tt liangshengtong@pku.edu.cn}.}

\author{Yuehan Shao\fnref{fn4}}
\fntext[fn4]{School of Mathematical Sciences, Peking University, Beijing, China, email: {\tt Joshao\_s@pku.edu.cn}.}

\author{Min Tang\fnref{fn5}}
\fntext[fn5]{Institute of Natural Sciences, School of Mathematics, Shanghai Jiaotong University, Shanghai, China, email: {\tt tangmin@sjtu.edu.cn}.}

\author{Yanli Wang\fnref{fn6}}
\fntext[fn6]{Beijing Computational Science Research Center, Beijing, China, email: {\tt ylwang@csrc.ac.cn}.
}

\begin{abstract}
We present an asymptotic-preserving (AP) numerical method for solving the three-temperature radiative transfer model, which holds significant importance in inertial confinement fusion. A carefully designed splitting method is developed that can provide a general framework of extending AP schemes for the gray radiative transport equation to the more complex three-temperature radiative transfer model. The proposed scheme captures two important limiting models: the three-temperature radiation diffusion equation (3TRDE) when opacity approaches infinity and the two-temperature limit when the ion-electron coupling coefficient goes to infinity. We have rigorously demonstrated the AP property and energy conservation characteristics of the proposed scheme and its efficiency has been validated through a series of benchmark tests in the numerical part.

\end{abstract}


\maketitle

\section{Introduction}

The radiative transfer equation holds significant importance in inertial confinement fusion, and accurately modeling thermal radiative transfer presents a formidable challenge. In many radiative transfer scenarios, the timescale required for electron and ion temperatures to approach thermal equilibrium exceeds the timescale needed to achieve quasi-neutral states, as noted by Enaux et al. in \cite{enaux2020numerical}. Consequently, it becomes necessary to employ a three-temperature system, wherein radiation, ion, and electron energies are considered separately. For example, in the study conducted by Shiva et al. \cite{effect2016}, a three-temperature radiation hydrodynamic model is employed to show the influence of electron thermal radiation on laser ablative shock waves propagating through ambient atmospheric air.

Numerical simulations of the three-temperature radiation transport model (3TRTM) pose substantial challenges, as highlighted by Evans et al. in \cite{EVANS20071695}. These challenges arise primarily due to the following reasons: 1) The intrinsic nonlinearity of the problem: The radiation intensity exhibits a nonlinear dependence on electron temperature, and the governing equations for electrons and ions manifest nonlinearity in diffusion and source terms. 2) High dimensionality of radiation intensity: The radiation intensity is a function of seven independent variables. 3) Restrictive time step constraints: The speed of light imposes stringent limitations on the time step length during simulations, as discussed by Enaux et al. in \cite{enaux2020numerical}. 4) Multiscale nature of physical parameters: The physical parameters, including opacity and the ion-electron coupling coefficient, exhibit significant variability across several orders of magnitude within the computational domain of certain benchmark tests.

The 3TRTM has two important limits: the three-temperature radiation diffusion equation (3TRDE) emerges as opacity approaches infinity, and the two-temperature limit arises when the ion-electron coupling coefficient tends towards infinity. The 3TRDE depends solely on space and time, resulting in a reduction in the dimension of independent variables. Nevertheless, it inherits the same challenges owing to nonlinearity and the multiscale ion-electron coupling coefficient. In the two-temperature limit, where ion and electron temperatures coincide, the model shares a similar form as the gray radiation transport model (GRTM), which only considers radiation-electron interactions.

Both of these limiting models, while simpler than the 3TRTM, present formidable challenges and have gathered significant attention from researchers. The 3TRDE comprises three reaction-diffusion-type equations with nonlinear diffusion and stiff nonlinear reactions. Efforts are made to develop robust and efficient schemes capable of handling nonlinearity in arbitrary geometries. For example, a positivity-preserving and conservative scheme is proposed in \cite{3T2022, Posi2020} and the discretization on unstructured grids can be found in \cite{3T2015, Combin2016}. The pursuit of maximum principle-preserving schemes is discussed in \cite{finite2019}, while semi-analytic solutions for the 3TRDE are explored in \cite{Mcclarren2011}. To enhance computational efficiency, adaptation, and multigrid methods \cite{Jiang2007, Parallel2004}, and certain iterative nonlinear solvers \cite{Anderson2017, Operator2019, Enaux2022, tang2021} have been proposed.

GRTM, being a high-dimensional system, faces similar challenges as the 3TRTM. Numerical methods for solving GRTM have been extensively studied and usually fall into two categories: stochastic and deterministic methods. The widely adopted stochastic method is the implicit Monte Carlo method (IMC), originally introduced for the two-temperature model by Fleck and Cummings (1971) and subsequently improved in \cite{shiyi2020}. It has also been extended to the three-temperature model in \cite{EVANS20071695,shiyi2023}. However, the IMC method is limited by the statistical noise, and its ability to handle multiscale parameters for GRTE and 3TRTM has not been thoroughly explored. An AP scheme means that the discrete numerical scheme of the microscopic model will converge to a discretization of the macroscopic model as the scale parameter approaches the limit. To handle multiscale parameters and the high speed of light, numerous AP schemes for GRTE have been proposed in the literature. For instance, in \cite{semi2008Ryan}, McClarren et al. designed a semi-implicit time integration method that explicitly treats the streaming term while handling the material coupling term implicitly with a linearized emission source. The unified gas kinetic scheme, as developed in \cite{sun2015asymptotic1}, extends the unified gas kinetic scheme for the linear radiative transport equation to the nonlinear GRTE by employing a linearized iterative solver for the material thermal temperature. An AP scheme based on a prediction-correction approach that allows for larger time steps and can capture the front position of the degenerate case is proposed in \cite{tang2021jcp}. A high-order AP scheme for GRTE is detailed in \cite{xiongtao}. However, the literature on efficient AP solvers for 3TRTM is much less.  

In this work, we introduce an AP splitting method to address the complexities of the 3TRTM. A carefully designed splitting method is developed that can provide a general framework for extending AP schemes for the GRTE to the more complex 3TRTM. The main idea involves dividing the nonlinear interactions between radiation, electrons, and ions temperatures in a specific way, such that the resulting split systems can effectively capture the diffusion limit when opacity tends towards infinity and the two-temperature limits when the collision time coefficient approaches infinity. More precisely, the original 3TRTM is divided into two parts. The first part consists of a seven-dimensional equation governing radiation density flux, featuring a modified source term, along with two ordinary differential equations (ODEs) for electron and ion temperatures. The second part is composed of three equations governing radiation density, electron, and ion temperatures, with all three functions being macroscopic, depending solely on space and time. To relax the time step requirement, a fully implicit solver is employed for both parts, solving nonlinear systems through carefully designed alternating iterations. While the first system is characterized by high dimensionality and nonlinearity, it is almost the same as the GRTE and can be effectively addressed by extending the efficient AP solvers proposed in the existing literature for GRTE, as discussed in the preceding paragraph. In this paper, we extend the AP scheme introduced in \cite{Fu2022} for GRTE to 3TRTM and apply it to the alternating iteration method, eliminating the need to solve a large nonlinear system. On the other hand, the structure of the second system shows similarity to the 3TRDE, making it both easier to solve and the various robust and efficient existing schemes \cite{tang2021jcp} can be utilized to solve it. For simplicity, we employ a fully implicit solver. The nonlinear system obtained by implicit discretization is also solved using the alternating iteration method. We subsequently establish the AP properties in two distinct limits: the 3TRDE and GRTE, and validate the energy conservation of the numerical scheme. 

The rest of this paper is organized as follows. In Sec. \ref{sec:model_AP}, the details of 3TRTM and the two asymptotic limits are discussed. The AP temporal splitting method is proposed in Sec. \ref{sec:split}. The nonlinear iterative solvers for both systems are discussed in Sec. \ref{sec:time} with particular emphasis on the first high-dimensional system. The spatial discretization together with a discussion of the AP properties is shown in Sec. \ref{sec:IMEX}. To validate the performance of the scheme, numerical simulations are conducted in Sec. \ref{sec:num}, including tests of the AP properties and tests on the homogeneous model problem, the Marshak wave problem with and without conduction terms, and a two-dimensional Riemann problem. Finally, some conclusion remarks are made in Sec. \ref{sec:conclusion}. 

\section{Three-temperature radiative transfer model}
\label{sec:model_AP}

In this section, we will present the detailed formulation of 3TRTM and focus on the derivation of the two limiting models: 3TRDE and the two-temperature limit.

\subsection{The 3TRTM} 
\label{sec:model}

In 3TRTM, the temperatures for ions ($T_i$), electrons ($T_e$), and the radiation field ($\psi$) are coupled together to represent the evolution of the radiation field, electron energy, and ion energy density in a hot plasma. The detailed model after nondimensionalization is as follows \cite{EVANS20071695}.
\begin{subequations} 
\label{eq:model}
\begin{align} 
    & \frac{\epsilon^2}{c}\frac{\partial \psi}{\partial t}+ \epsilon \Omega\cdot \nabla \psi=\frac{1}{|\Omega|}\sigma acT_e^4 - \sigma \psi + \frac{Q_r}{|\Omega|}, \label{eq:rte1}\\
    & \frac{\epsilon^2 \partial E_e}{\partial t}= \epsilon^2 \nabla \cdot D_e\nabla T_e+c\kappa(T_i-T_e)+ \sigma \int_{\mathbb{S}^2}\psi \dd \Omega -\sigma ac T_e^4 + Q_e,    \label{eq:rte2}  \\
    & \frac{\epsilon^2 \partial E_i}{\partial t}= \epsilon^2 \nabla \cdot D_i\nabla T_i +c\kappa(T_e-T_i) + Q_i. \label{eq:rte3} 
\end{align}
\label{eq:rte}
\end{subequations}
Here $\psi(t,\bx,\Omega)$ represents the radiation intensity, where $\Omega$ denotes the angular variable located on the unit sphere $\mathbb{S}^2$, and $|\Omega|$ is the surface area of $\Omega$, equaling $4\pi$. The dimensionless parameter $\epsilon$ is analogous to the Knudsen number in the Boltzmann equation, representing the ratio of the mean free path to the characteristic length. $\sigma(\bx, T_e)$ denotes the opacity, a material-specific coefficient. $c$ stands for the speed of light. The radiation constant $a$ is defined as 
\begin{equation}
    a=\frac{8\pi k^4}{15h^3c^3},
\end{equation}
where $h$ represents Planck's constant, and $k$ denotes the Boltzmann constant.
$E_e(t, \bx)$ and $T_e(t, \bx)$ refer to the energy and temperature of the electrons, while $E_i(t, \bx)$ and $T_i(t, \bx)$ represent the energy and temperature of the ions. The relationship between material temperature $T_s(t,\bx)$ and material energy density $E_s(t,\bx)$ is 
\begin{equation}
\label{eq:cv}
    E_s=C_{v_s}T_s> 0,\qquad s=e,i,
\end{equation}
where $C_{v_s} > 0$ represents the heat capacity. $Q_r, Q_e$, and $Q_i$ are the source terms for radiation, electrons, and ions, respectively. The thermal diffusion coefficients $D_s(\bx, T_s)$ in the conduction terms $\nabla \cdot D_s \nabla T_s$ (where $s=e,i$) exhibit nonlinear dependence on temperature as 
\begin{equation}
    D_s=K_{s}T_s^{5/2},\qquad s=e,i,
\end{equation}
where $K_s$ represents the plasma-dependent coefficient \cite{enaux2020numerical}. The timescale of ion-electron coupling is characterized by the collision time coefficient $\kappa$, which can take various forms, as documented in the literature \cite{3T2021, EVANS20071695, 3T2015}. In this paper, $\kappa$ is treated as a constant but can vary significantly over several orders of magnitude. The radiation energy $E_r$ is defined as
\begin{equation}
    E_r=\frac{1}{c}\int_{\mathbb{S}^2}\psi\dd\Omega,
\end{equation}
and the radiation temperature $T_r$ follows the relation
\begin{equation}
\label{eq:T_r}
    E_r =aT_r^4,\qquad T_r\geqslant  0.
\end{equation}

In Eq. \eqref{eq:rte}, \eqref{eq:rte1} describes the evolution of radiation intensity, while \eqref{eq:rte2} and \eqref{eq:rte3} govern the evolution of electron and ion temperatures. These equations couple conduction and energy exchange among radiation, electrons, and ions. Integrating \eqref{eq:rte1} over $\Omega$ and combining it with \eqref{eq:rte2} and \eqref{eq:rte3} yields the conservation of total energy, expressed as
\begin{equation}
\label{eq:rteenergy}
    \fpp{E_e}{t}+ \fpp{ E_i}{t}+  \fpp{E_r}{t} + \nabla \cdot \bm{F}_r = \nabla\cdot D_e\nabla T_e+\nabla\cdot D_i\nabla T_i + Q_r/\epsilon^2 + Q_e/\epsilon^2 + Q_i/\epsilon^2.
\end{equation}
Here, $\bm{F}_r(t, \bx)$ represents the radiation flux
\begin{equation}
\label{eq:fluxF}
    \bm{F}_r(t, \bx) = \int_{\mathbb{S}^2} \Omega \psi \dd \Omega. 
\end{equation}

\subsection{Asymptotic behavior}
\label{sec:limit}
In the 3TRTM \eqref{eq:model}, two limiting models emerge as the dimensionless parameter $\epsilon$ approaches zero or when the collision time coefficient $\kappa$ tends to infinity. These limits provide valuable approximations to the original model and have been thoroughly investigated in \cite{enaux2020numerical}. For the sake of generality, we assume the source term $Q_s$ (where $s=r,e,i$) to be zero. In this section, we will present a detailed discussion of these two limits.

\paragraph{Diffusion limit}
As $\epsilon$ approaches zero, the radiation intensity $\psi$ becomes independent of the angular variable $\Omega$ and can be described by the radiation temperature $T_r$. Let $T_r^{(0)},T_e^{(0)}$, and $T_i^{(0)}$ represent the leading-order terms of $T_r,T_e$, and $T_i$, respectively. This leads to an equilibrium diffusion limit, where the corresponding local temperature $T^{(0)}$ satisfies a nonlinear diffusion equation, expressed as follows \cite{fleck1971}.
\begin{subequations} 
\begin{align}
    &T_r^{(0)}=T_e^{(0)}=T_i^{(0)}=T^{(0)}, \qquad E_r(T^{(0)})  = a \left(T^{(0)}\right)^4,\\
    &\begin{aligned}
    &\fpp{E_e(T^{(0)})}{t}+\fpp{E_i(T^{(0)})}{t}+\fpp{E_r(T^{(0)})}{t}\\
    & \hspace{2cm} + \nabla\cdot \bm{F}(t, \bx) = \nabla\cdot D_e\nabla T^{(0)}+\nabla\cdot D_i\nabla T^{(0)}.
    \end{aligned}
\end{align}
\label{eq:rtelimit}
\end{subequations}
In this approximation, the radiative flux $\boldsymbol{F}(t,\bx)$ is related to the temperature through Fick's law of diffusion, expressed as
\begin{equation}
\label{eq:flux_F}
    \boldsymbol{F}(t,\bx) = \int_{\mathbb{S}^2}\Omega\psi\dd\Omega= -\frac{ac}{3\sigma}\nabla (T^{(0)})^4.
\end{equation}
Eq. \eqref{eq:rtelimit} is referred to as the equilibrium diffusion limit \cite{enaux2020numerical}. This limit provides a simplified description of the behavior of the system when the ratio of the mean free path to the characteristic length $\epsilon$ becomes very small. 

\paragraph{Two temperature limit}
When the collision time coefficient $\kappa$ approaches infinity, the system converges to a two-temperature limit, resulting in a simplified two-temperature model. In this limit, the electron temperature $T_e$ and the ion temperature $T_i$ become equal, denoted as the material temperature $T$, and the radiation transfer equations \eqref{eq:rte}  can be reduced to a 2T (two-temperature) radiation transfer equation by summing up \eqref{eq:rte2} and \eqref{eq:rte3}. The resulting 2T radiation transfer equation is
\begin{subequations}
\label{eq:limit_kappa}
\begin{align}      
    & \frac{\epsilon^2}{c}\frac{\partial \psi}{\partial t}+ \epsilon \Omega\cdot \nabla \psi=\frac{1}{4\pi}\sigma acT^4 - \sigma \psi, \\
    & \epsilon^2 C_v\frac{ \partial T}{\partial t}= \epsilon^2 \nabla \cdot D_T\nabla T + \sigma \int_{\mathbb{S}^2}\psi \dd \Omega -\sigma ac T^4,     \\
    & C_v = C_{v_e} + C_{v_i}, \qquad D_T = D_{e} + D_{i},
\end{align}
\end{subequations}
where $C_v$ ($D_T$) represents the heat capacity (thermal diffusion coefficient), considering both the heat capacity (diffusion coefficient) of electrons $C_{v_e}$ ($D_e$) and that of ions $C_{v_i}$ ($D_i$).

\section{Asymptotic preserving time splitting scheme}
\label{sec:split}

In many physically meaningful benchmark tests, $\epsilon$ can vary from very small values to $\mathcal{O}(1)$, and $\kappa$ can span several orders of magnitude, ranging from $\mathcal{O}(1)$ to very large values. As a result, the numerical scheme should maintain uniform accuracy for all these extreme cases. Therefore, it is important to design an efficient scheme that is AP in both the diffusion limit as $\epsilon$ approaches zero and the 2T limit as $\kappa$ tends towards infinity. Moreover, the preservation of total energy \eqref{eq:rteenergy} is of most importance for some benchmark tests, especially when the regularity of the solution is low \cite{enaux2020numerical}. In this section, we propose an AP conservative time-splitting method. 

In this method, the 3TRTM is split into two systems. The first part is close to the GRTE and can be effectively addressed by extending the existing efficient AP solvers. The second part is a macroscopic model, where the radiation density, electron temperature, and ion temperature depend solely on space and time.

\subsection{The AP splitting scheme}
\label{sec:AP_split}

The splitting scheme is proposed in the framework of micro-macro decomposition \cite{split2021,3D2006}. The radiation intensity is decomposed into a macroscopic density part that is independent of the angular variable $\Omega$ and the remaining part:
\begin{equation}
\label{eq:MM}
    \psi(t,\bx,\Omega) = \rho(t,\bx) + g(t,\bx,\Omega), 
\end{equation}
with 
\begin{equation}
\label{eq:MM_rhog}
    \rho(t,\bx) =  \langle \psi(t,\bx,\Omega) \rangle,   \qquad g(t,\bx,\Omega) = \psi (t,\bx,\Omega)- \rho(t,\bx),  \qquad \langle \cdot \rangle = \frac{1}{|\Omega|} \int_{\mathbb{S}^2} \cdot  \dd \Omega.
\end{equation}
By integrating \eqref{eq:rte1} over the angular space, one can derive the equation for $\rho$ that depends on $g$. Subtracting the equation for $\rho$ from \eqref{eq:rte1} leads to the following micro-macro decomposition system: 
\begin{subequations}
\label{eq:MMRTE}
\begin{numcases} {}
    \frac{\epsilon^2}{c} \fpp{\rho}{t} + \epsilon \nabla \cdot \langle\Omega g \rangle  = -\sigma\left( \rho  -  \frac{ a c T_e^4}{|\Omega|}\right), \label{eq:MMRTE1}\\
    \frac{\epsilon^2}{c} \fpp{g}{t} +\epsilon \nabla \cdot [\Omega (\rho + g) - \langle \Omega g \rangle] = -\sigma g, \label{eq:MMRTE2}\\
    \frac{\epsilon^2 \partial E_e}{\partial t}= \epsilon^2 \nabla \cdot D_e\nabla T_e+c\kappa(T_i-T_e)+ \sigma ( |\Omega|\rho - ac T_e^4), \label{eq:MMRTE3}   \\
    \frac{\epsilon^2 \partial E_i}{\partial t}= \epsilon^2 \nabla \cdot D_i\nabla T_i +c\kappa(T_e-T_i). \label{eq:MMRTE4}
\end{numcases}
\end{subequations}
The computational time $[0,t_{\rm end}]$ is divided into time intervals $[t^n,t^{n+1}]$ for $n=0,1,2,\cdots$. Utilizing this micro-macro decomposition, the AP time-splitting scheme is formulated as follows.
\begin{enumerate}
\item With the approximated solution $(\rho^{n},g^{n},E_e^{n},E_i^{n})$ at time $t^n$, we solve the first system
\begin{subequations}
\label{eq:MM_split_2_1}
\begin{align} 
    \label{eq:MM_split_2_1_1}
    &\frac{\epsilon^2}{c}\frac{\partial \rho}{\partial t} 
    + \epsilon \nabla \cdot  \langle \Omega g \rangle = -\frac{\sigma}{2}\left(\rho - \frac{1}{|\Omega|} acT_e^4\right),\\            
    \label{eq:MM_split_2_1_0}
    & \frac{\epsilon^2}{c} \fpp{g}{t} +\epsilon \nabla \cdot [\Omega (\rho + g) - \langle \Omega g \rangle] = -\sigma g, \\
    \label{eq:MM_split_2_1_2}
    & \epsilon^2\fpp{E_e}{t}=\frac{1}{2}c\kappa (T_i-T_e)+ \frac{\sigma}{2}(|\Omega| \rho - acT_e^4),\\
    \label{eq:MM_split_2_1_3}
    &\epsilon^2 \fpp{E_i}{t}=\frac{1}{2}c\kappa (T_e-T_i),
\end{align}
\end{subequations} 
for one time step $\Delta t = t^{n+1}-t^n$ and get $(\rho^{n+1,\ast},g^{n+1,\ast},E_e^{n+1,\ast},E_i^{n+1,\ast})$. 
\item Using $(\rho^{n+1,\ast},g^{n+1,\ast},E_e^{n+1,\ast},E_i^{n+1,\ast})$ as the initial condition, we update the second system
\begin{subequations} 
\label{eq:MM_split_2_2}
\begin{align}
    \label{eq:MM_split21}
    &\frac{\epsilon^2}{c}\fpp{\rho}{t}= - \frac{\sigma}{2} \left(\rho - \frac{acT_e^4}{|\Omega|} \right), \\
    &\frac{\epsilon^2}{c}\fpp{g}{t}=0,\\
    \label{eq:MM_split22}
    &\epsilon^2 \fpp{E_e}{t}= \epsilon^2 \nabla \cdot  D_e \nabla T_e +\frac{1}{2}c\kappa(T_i-T_e)+ \frac{\sigma}{2}(|\Omega| \rho - acT_e^4),\\
    \label{eq:MM_split23}
    &\epsilon^2 \fpp{E_i}{t}=\epsilon^2 \nabla \cdot D_i\nabla T_i +\frac{1}{2}c\kappa(T_e-T_i).
\end{align}
\end{subequations}
for one time step $\Delta t = t^{n+1}-t^n$ and get $(\rho^{n+1},g^{n+1},E_e^{n+1},E_i^{n+1})$. 
\end{enumerate}

In the following, we will demonstrate that this splitting approach effectively captures both the diffusion limit and the two-temperature limit while conserving the total energy. 

\paragraph{Diffusion limit}  It is straightforward to demonstrate that as $\epsilon$ tends to zero in \eqref{eq:MM_split_2_1}, the leading-order terms $T_r^{(0)},T_e^{(0)}$, and $T_i^{(0)}$, satisfy
\begin{equation}
\label{eq:diff_limit_0}
    T_e^{(0)}=T_i^{(0)}= T^{(0)}, \quad E_r(T_r) =\frac{1}{c}\int_{\mathbb{S}^2}\psi\dd\Omega=\frac{|\Omega|}{c}\rho= a T_e^4+\mathcal{O}(\epsilon)=E_r(T^{(0)})+\mathcal{O}(\epsilon), \quad g = \mathcal{O}(\epsilon). 
\end{equation}
Thus, from the relation between $E_r$ and $T_r$ in \eqref{eq:T_r}, we have $T^{(0)}_r=T^{(0)}$. By adding up \eqref{eq:MM_split_2_1_1}, \eqref{eq:MM_split_2_1_2}, and  \eqref{eq:MM_split_2_1_3} together, the source terms on the right-hand sides of these three equations cancel each other, resulting in the following conservative equation:
\begin{equation}
\label{eq:diff_limit_1}
    \frac{E_r(T^{(0)})}{\partial t}+\frac{\partial E_e(T^{(0)})}{\partial t}+ \frac{\partial E_i(T^{(0)})}{\partial t}=-\frac{1}{\epsilon}\nabla \cdot  \langle \Omega g \rangle + \mathcal{O}(\epsilon).
\end{equation}
From \eqref{eq:MM_split_2_1_0}, one gets
\begin{equation}
\label{eq:diff_limit_2}
    g = -\frac{\epsilon}{\sigma} \nabla \cdot [\Omega (\rho+g) - \langle \Omega g\rangle]
    -\frac{\epsilon^2}{c\sigma}\frac{\partial g}{\partial t}.
\end{equation}
Multiplying both sides of \eqref{eq:diff_limit_2} by $\Omega$ and integrating with respect to angular variable, one obtains
\begin{equation}
    -\frac{1}{\epsilon}\nabla\cdot \langle\Omega g\rangle=\frac{1}{\sigma} \nabla \cdot [\langle\Omega\Omega \rangle \rho ]+\mathcal{O}(\epsilon).
\end{equation}
Here, we have used the fact that $g(t,\bx,\Omega)$ is of $\mathcal{O}(\epsilon)$. Therefore, \eqref{eq:diff_limit_1} can be expressed as follows.
\begin{equation}
    \label{eq:MM_limit_split_2_1}
    \begin{aligned}
        &\frac{E_r(T^{(0)})}{\partial t}+\frac{\partial E_e(T^{(0)})}{\partial t}+ \frac{\partial E_i(T^{(0)})}{\partial t}=\nabla \cdot \left(\frac{ac}{3\sigma} \nabla (T^{(0)})^4\right).
    \end{aligned}
\end{equation}
Similarly, as $\epsilon\to 0$, the leading-order terms of \eqref{eq:MM_split_2_2} yield 
\begin{equation}
    T_r^{(0)}=T_e^{(0)}=T_i^{(0)}= T^{(0)}, \qquad \frac{|\Omega|}{c}\rho=E_r(T^{(0)})+\mathcal{O}(\epsilon).
\end{equation}
By adding up \eqref{eq:MM_split21}, \eqref{eq:MM_split22}, and \eqref{eq:MM_split23}, the interaction terms cancel out, and the limiting model for \eqref{eq:MM_split_2_2} can be expressed as 
\begin{equation}
\label{eq:MM_limit_split_2_2}
    \begin{aligned}
    &\frac{E_r(T^{(0)})}{\partial t}+\frac{\partial E_e(T^{(0)})}{\partial t}+ \frac{\partial E_i(T^{(0)})}{\partial t}= \nabla \cdot  D_e \nabla T^{(0)} + \nabla \cdot  D_i  \nabla T^{(0)}.
    \end{aligned}
\end{equation}
It is evident that the combination of  \eqref{eq:MM_limit_split_2_1} and \eqref{eq:MM_limit_split_2_2} constitutes a time-splitting scheme for the diffusion limit equation \eqref{eq:rtelimit}.

\paragraph{Two temperature limit}
As $\kappa$ tends toward infinity, let $T$ be the material temperature in the limiting case. and it holds that 
\begin{equation}
\label{eq:two_limit_0}
    T_e = T + \mathcal{O}(\kappa^{-1}),\quad T_i = T + \mathcal{O}(\kappa^{-1}).     
\end{equation}
We can derive the limiting equations for \eqref{eq:MM_split_2_1} by summing up \eqref{eq:MM_split_2_1_2} and \eqref{eq:MM_split_2_1_3}. The resulting limiting equations are 
\begin{equation}
\label{eq:MM_limit_kappa_2_1}
    \begin{aligned}
    &\frac{\epsilon^2}{c}\frac{\partial \rho}{\partial t} 
    + \epsilon \nabla \cdot \langle \Omega  g \rangle = -\frac{\sigma}{2}\left(\rho - \frac{acT^4}{|\Omega|} \right),\\             
    &  \frac{\epsilon^2}{c} \fpp{g}{t} +\epsilon \nabla \cdot  [\Omega (\rho +g) - \langle \Omega g \rangle] = -\sigma g, \\
    & \epsilon^2\fpp{E_T}{t}  =  \frac{\sigma}{2}\big(|\Omega| \rho - acT^4\big),
    \end{aligned}
\end{equation}
where $ E_T = (C_{v_e} + C_{v_i}) T$. For the second system \eqref{eq:MM_split_2_2}, as $\kappa\to\infty$, one obtains \eqref{eq:two_limit_0}. By summing up \eqref{eq:MM_split22} and \eqref{eq:MM_split23}, the limiting system of \eqref{eq:MM_split_2_2} becomes
\begin{equation}
\label{eq:MM_limit_kappa_2_2}
    \begin{aligned}
    &\frac{\epsilon^2}{c}\fpp{\rho}{t}= - \frac{\sigma}{2} \left(\rho -\frac{acT^4}{|\Omega| } \right), \\
    &\epsilon^2 \fpp{E_T}{t} = \epsilon^2 \nabla \cdot  D_T\nabla T + \frac{\sigma}{2}(|\Omega| \rho - acT^4),
    \end{aligned}
\end{equation} 
where $E_T = \big(C_{v_e} + C_{v_i}\big) T$ and $D_T = D_{e}(T) + D_{i}(T)$. It is evident that the combination of equations \eqref{eq:MM_limit_kappa_2_1} and \eqref{eq:MM_limit_kappa_2_2} constitutes a time-splitting method for the two-temperature limit model \eqref{eq:limit_kappa}.

\subsection{Comparison with another splitting method}

Compared to the original system, in the first system \eqref{eq:MM_split_2_1}, spatial derivatives only exist in \eqref{eq:MM_split_2_1_1} and \eqref{eq:MM_split_2_1_0}, while the equations for the electron and ion temperatures reduce into ODEs as in \eqref{eq:MM_split_2_1_2} and \eqref{eq:MM_split_2_1_3}. On the other hand, the second part \eqref{eq:MM_split_2_2} is a macroscopic system. In some benchmark problems, the conduction terms $\nabla \cdot D_s \nabla T_s, s = e, i$ are set to be zero  \cite{shiyi2023} and the second system \eqref{eq:MM_split_2_2} reduces to an ODEs system that can be easily solved. The specific choices of 1) dividing the term $\sigma\big(\rho-\frac{acT_e^4}{|\Omega|}\big)$ in \eqref{eq:MMRTE1} and \eqref{eq:MMRTE3} into two $\frac{\sigma}{2}\big(\rho-\frac{acT_e^4}{|\Omega|}\big)$, and placing them in \eqref{eq:MM_split_2_1_1}, \eqref{eq:MM_split_2_1_2}, \eqref{eq:MM_split21}, and \eqref{eq:MM_split22}; 2) dividing the term $c\kappa(T_i-T_e)$ in \eqref{eq:MMRTE3} and \eqref{eq:MMRTE4} into two $\frac{c\kappa}{2}(T_i-T_e)$, and placing them in \eqref{eq:MM_split_2_1_2}, \eqref{eq:MM_split_2_1_3}, \eqref{eq:MM_split22}, and \eqref{eq:MM_split23}; are crucial for achieving the AP properties.

A splitting scheme for the macroscopic diffusion model is proposed in \cite{split2021, 3D2006}. Extending this approach, one can design a splitting scheme for the 3T model by a straightforward extension. We split the 3T model into a classical two-temperature radiative transfer system and a macroscopic system for the electron and ion temperatures, as shown below:
\begin{align}
\label{eq:split_1_1}
    & \text{Microscopic part:}   \left\{
    \begin{aligned}
    &  \frac{\epsilon^2}{c}\frac{\partial \psi}{\partial t}+\epsilon \Omega \cdot \nabla \psi  =   -\sigma (\psi - \frac{1}{|\Omega|} acT_e^4), \\
    &  \epsilon^2 \frac{\partial E_e}{\partial t}= \sigma( |\Omega| \rho  -  ac T_e^4), \\
    &  \epsilon^2 \frac{\partial E_i}{\partial t} = 0. 
    \end{aligned} 
    \right. \\
\label{eq:split_1_2}
    & \text{Macroscopic part:} \left\{
    \begin{aligned}
    & \fpp{\psi}{t} = 0, \\
    &\epsilon^2 \frac{\partial E_e}{\partial t}=\epsilon ^2\nabla \cdot D_e \nabla T_e +c\kappa(T_i-T_e),\\
    & \epsilon^2 \frac{\partial E_i}{\partial t}= \epsilon^2 \nabla \cdot D_i \nabla T_i
    +c\kappa(T_e-T_i).
    \end{aligned} 
    \right.
\end{align}
Here \eqref{eq:split_1_1} has the same form as the GRTE, while \eqref{eq:split_1_2} is a macroscopic system that depends solely on $T_e$ and $T_i$. However, as $\epsilon$ approaches zero, the diffusion limit of the splitting system \eqref{eq:split_1_1} and \eqref{eq:split_1_2} becomes
\begin{align}
    &  \left\{
    \begin{aligned}
    & \fpp{E_r(T^{(0)})}{t}+\fpp{E_e(T^{(0)})}{t} = \nabla \cdot \frac{ac}{3\sigma} \nabla \left(T^{(0)}\right)^4, \\
    &  T_e = T_r = T^{(0)}, \quad E_r(T^{(0)}) = a\left(T^{(0)}\right)^4,
    \end{aligned} 
    \right.\\
    &  \left\{
    \begin{aligned}
    & \fpp{E_e(T^{(0)})}{t}+\fpp{E_i(T^{(0)})}{t} = \nabla \cdot D_e \nabla T_e  + \nabla \cdot D_i \nabla T_i,\\
    & T_e = T_i = T^{(0)}.
    \end{aligned} 
    \right.
\end{align}
This method is no longer a time-splitting approach for the diffusion limit equation \eqref{eq:rtelimit}. 

\section{Temporal and angular discretization for the subsystems}
\label{sec:time}

In this section, we introduce the temporal and angular discretization based on the splitting method in Sec. \ref{sec:split}. We consider the nonlinearity in the coefficients $\sigma$ and $D_s$ and employ a fully implicit solver for both systems to ensure accuracy through the AP properties discussed in the previous section. However, the computational cost of nonlinear iterations for a high-dimensional system is extremely high, and efficient iterative solvers must be designed. We adopt two alternating iteration methods, so that, under a larger time step, these iterative methods can still converge quickly.

\subsection{Splitting of \texorpdfstring{$P_N$}{PN} system for angular discretization}
\label{sec:Pn}

We discretize the angular variable by $\Pn$ method, where the radiative intensity is approximated by a series of basis functions \cite{implicit2016Laboure, semi2008Ryan, On2008}. Without loss of generality, we consider the following 3TRTM in slab geometry:
\begin{subequations}
\label{eq:rte_1D}
\begin{align}
\label{eq:rte_1D1}
    & \frac{\epsilon^2}{c}\frac{\partial \psi}{\partial t}+\epsilon  \mu \frac{\partial \psi}{\partial x}= -\sigma\psi+\frac{1}{2}\sigma acT_e^4, \\
\label{eq:rte_1D2}
    & \frac{\epsilon^2 \partial E_e}{\partial t}=\epsilon^2  \fpp{}{x}\left(D_e \fpp{T_e}{x}\right)+c\kappa(T_i-T_e)+ \sigma \int_{-1}^{1} \psi \dd \mu -\sigma ac T_e^4,\\
\label{eq:rte_1D3}
    &  \frac{\epsilon^2 \partial E_i}{\partial t}= \epsilon^2  \fpp{}{x}\left(D_i \fpp{T_i}{x} \right)+c\kappa(T_e-T_i),
\end{align}
\end{subequations}
where $\mu \in [-1, 1]$. The details of how to project a spatially 3D RTE with angular variables defined on the surface of a 3D sphere to slab geometry can be found in \cite{ComputationEELewis}, and we omit the details here.

In the framework of the $P_N$ method \cite{On2008}, the radiative intensity $\psi$ is approximated as 
\begin{equation}
    \label{eq:Pn}
    \psi(t, x, \mu) \approx \sum_{l = 0}^M \frac{2l+1}{2} \psi_l(t, x) P_l(\mu),
\end{equation}
where $P_l(\mu)$ represents the Legendre polynomial of degree $l$ and $M\in\mathbb{N}$ is the expansion order. Utilizing the following orthogonality property of the Legendre polynomial: 
\begin{equation}
\label{eq:orth}
    \int_{-1}^1P_{m}(\mu)P_n(\mu) \dd\mu=\frac{2}{2n+1}\delta_{m,n},\quad m,n\in\mathbb{N},
\end{equation}
the expansion coefficients $\psi_l(t, x)$ are
\begin{equation}
    \label{eq:coe}
    \psi_l(t, x) = \int_{-1}^1 P_l(\mu) \psi(t, x, \mu) \dd \mu, \qquad l = 0, \cdots, M. 
\end{equation}
Therefore, it holds that
\begin{equation}
\psi_0=2\rho.
\end{equation}
By substituting \eqref{eq:Pn} into the governing equations of the 1D model \eqref{eq:rte_1D1} and utilizing the orthogonality of the basis functions \eqref{eq:orth}, as well as the following recursive property: 
\begin{equation}
\label{eq:rec}
    P_{k+1}(\mu) = \frac{2k+1}{k+1}\mu P_k(\mu) - \frac{k}{k+1}P_{k-1}(\mu),\quad k\in\mathbb{N}^+,
\end{equation}
the $P_N$ system can be written as follows. 
\begin{subequations}
\label{eq:pn_sys}
\begin{align}
\label{eq:pn_sys1}
    &\frac{2\epsilon^2}{c}\frac{\partial \rho}{\partial t} + \epsilon \frac{\partial \psi_1}{\partial x}+2\sigma\rho=\sigma acT_e^4,\\
\label{eq:pn_sys2}
    &\frac{\epsilon^2}{c}\frac{\partial \psi_1}{\partial t} + \epsilon \left(\frac{2}{3} \fpp{\rho}{x}+\frac{2}{3}\fpp{\psi_{2}}{x}\right) +\sigma\psi_1=0,\\
\label{eq:pn_sys3}
    &\frac{\epsilon^2}{c}\frac{\partial \psi_l}{\partial t} + \epsilon \left(\frac{l}{2l+1} \fpp{\psi_{l-1}}{x}+\frac{l+1}{2l+1}\fpp{\psi_{l+1}}{x}\right) +\sigma\psi_l=0, \quad 2\leqslant l \leqslant M-1,\\
\label{eq:pn_sys4}
    &\frac{\epsilon^2}{c}\frac{\partial \psi_M}{\partial t} +  \frac{\epsilon M}{2M+1} \fpp{\psi_{M-1}}{x} +\sigma\psi_M=0. 
\end{align}
\end{subequations}
Comparing \eqref{eq:pn_sys} with \eqref{eq:MMRTE1} and \eqref{eq:MMRTE2}, the equation for the microscopic part $g$ becomes a system for $(\psi_1, \psi_2,\cdots, \psi_M)^T=(\psi_1,\bar{\bu}^T)^T$. Let coefficient matrices be ${\bf A}=(A_{i,j})  \in \bbR^{(M-1)\times (M+1)}$  and ${\bf B}=(B_{i,j}) \in \bbR^{(M-1) \times (M+1)}$, where
\begin{equation}
\label{eq:coe_AB}
    A_{l-1,l} = \frac{l}{2l+1}, \quad l = 2, \cdots, M,\qquad B_{l-1, l+2} = \frac{l+1}{2l + 1}, \quad l = 2, \cdots, M-1,
\end{equation}
and all other entries are zero. Eq. \eqref{eq:pn_sys} can then be rewritten as
\begin{subequations}
\label{eq:ma_form}
\begin{align}
\label{eq:ma_form_1}
    & \frac{2\epsilon^2}{c}\frac{\partial \rho}{\partial t} + \epsilon \frac{\partial \psi_1}{\partial x} =  -\sigma(2\rho - acT_e^4),\\
    & \frac{\epsilon^2}{c}\frac{\partial \psi_1}{\partial t} + \frac{2\epsilon}{3}\fpp{\rho}{x}+ \frac{2\epsilon}{3}\fpp{\psi_2}{x}=-\sigma\psi_1,\\
\label{eq:ma_form_2}
    & \frac{\epsilon^2}{c} \fpp{\bbu}{t} + \epsilon {\bf A}\fpp{\bu}{x} +\epsilon {\bf B}\fpp{\bu}{x} = -\sigma \bbu,
\end{align}
\end{subequations}
where $\bu=(2\rho,\psi_1,\bbu^T)^T$. The governing equation for the electron becomes
\begin{equation}
\label{eq:rte_1D2_1}
    \epsilon^2 \fpp{E_e}{t} = \epsilon^2 \fpp{}{x} \left(D_e \fpp{T_e}{x} \right)+ c \kappa (T_i - T_e) +  \sigma (2\rho - acT_e^4),
\end{equation}
and the equation for the conservation of the total energy \eqref{eq:rteenergy} is written as 
\begin{equation}
\label{eq:Pn_energy}
   \fpp{E_e}{t} + \fpp{E_i}{t} + \frac{2}{c} \fpp{\rho}{t} = - \frac{1}{\epsilon}\fpp{\psi_1}{x} + \fpp{}{x} \left(D_e \fpp{T_e}{x}\right) + \fpp{}{x} \left(D_i \fpp{T_i}{x}\right).
\end{equation}
Eqs. \eqref{eq:ma_form}, \eqref{eq:rte_1D2_1} and \eqref{eq:rte_1D3} together form the $\Pn$ system for the 1D 3TRTM. 

Utilizing the splitting method introduced in the previous section, the two parts of the splitting approach for this $P_N$ system are as follows.
\begin{subequations}
\label{eq:PN_full}
\begin{align}
\label{eq:PN_full_1}
    \text{First part (microscopic part): }&\left\{
    \begin{aligned}
    &\frac{2\epsilon^2}{c}\fpp{\rho}{t}+\epsilon\fpp{\psi_1}{x}=-\frac{\sigma}{2}(2\rho-acT_e^4),\\
    &\frac{\epsilon^2}{c}\fpp{\psi_1}{t}+\frac{2\epsilon}{3}\fpp{\rho}{x}+\frac{2\epsilon}{3}\fpp{\psi_2}{x}=-\sigma\psi_1,\\
    &\frac{\epsilon^2}{c}\fpp{\bbu}{t}+\epsilon\bA\fpp{\bu}{x}+\epsilon\bB\fpp{\bu}{x}=-\sigma\bbu,\\
    &\epsilon^2C_{v_e}\fpp{T_e}{t}=\frac{1}{2}c\kappa (T_i-T_e)+\frac{\sigma}{2}(2\rho-acT_e^4),\\
    &\epsilon^2C_{v_i}\fpp{T_i}{t}=\frac{1}{2}c\kappa (T_e-T_i),
    \end{aligned}
    \right.\\
\label{eq:PN_full_2}
    \text{Second part (macroscopic part): }&\left\{
    \begin{aligned}
    &\frac{2\epsilon^2}{c}\fpp{\rho}{t}=-\frac{\sigma}{2}(2\rho-acT_e^4),\\
    &\epsilon^2C_{v_e}\fpp{T_e}{t}=\epsilon^2\fpp{}{x}\left(D_e\fpp{T_e}{x}\right)+\frac{1}{2}c\kappa (T_i-T_e)+\frac{\sigma}{2}(2\rho-acT_e^4),\\
    &\epsilon^2C_{v_i}\fpp{T_i}{t}=\epsilon^2\fpp{}{x}\left(D_i\fpp{T_i}{x}\right)+\frac{1}{2}c\kappa (T_e-T_i),\\
    &\fpp{\psi_1}{t}=0,\quad \fpp{\bbu}{t}=\boldsymbol{0},\quad \boldsymbol{0}=(0,0,\dots,0)^T\in\mathbb{R}^{M-1}.
    \end{aligned}
    \right.
\end{align}
\end{subequations}

\subsection{Fully temporal discretization and the iterative solver}
\label{sec:semi-dis}

In the split $P_N$ system, the nonlinear terms include the fourth power of temperature and the nonlinear conduction term. Handling systems where both the conduction term and the fourth power term are implicit presents challenges. Explicitly treating the conduction term during the iteration process would impose a time step requirement of $\Delta t=\mathcal{O}(\Delta x^2)$. On the other hand, explicitly treating the fourth power term during the iteration process would result in a time step requirement dependent on $\epsilon$. While the fully implicit scheme is utilized for both systems, designing an efficient iterative solver is crucial to avoid solving a large nonlinear system, especially for the first system \eqref{eq:PN_full_1}.

To iteratively solve these systems, we design an alternating iterative method. It consists of two parts: one implicitly solves a linearized system, and the nonlinear part is considered explicitly. The other part implicitly solves the quartic equation at each space point, while the linearized system part is considered explicitly. The iterations of these two parts are performed alternately until certain convergence criteria are satisfied, and the numerical solution of the implicit system can be obtained. An AP solver for GRTE is developed in \cite{Fu2022, Fuimplicit}, whose computational cost is comparable to an explicit scheme, while its accuracy and stability are independent of the multiscale parameters $\epsilon$. In the subsequent part, we propose an efficient iteration method for \eqref{eq:PN_full_1} based on the AP solver for GRTE in \cite{Fu2022}. 

According to the order analysis in \cite{Fu2022}, the orders of the expansion coefficients are given by 
\begin{equation}
\label{eq:order}
    \psi_l = \mathcal{O}{(\epsilon^l)}, \qquad \fpp{\psi_l}{x} = \mathcal{O}{(\epsilon^l)}.
\end{equation}
For the convenience of readers, we provide the details of order analysis in \ref{app:order_analy}. In Eq. \eqref{eq:PN_full_1}, $\bA \fpp{\bu}{x}$ has a low order of $\epsilon$ and is treated implicitly, while the terms $\bB \fpp{\bu}{x}$ are treated explicitly since their orders are higher. The terms on the right-hand side are treated implicitly, while the coefficient $\sigma$, which may nonlinearly depend on $T_e$, is treated explicitly to simplify the iteration process. 

We now present the details of the iterative method for the microscopic part, as displayed in Eqs. \eqref{eq:iter_full_1} and \eqref{eq:iter_full_1_next}. Eq. \eqref{eq:iter_full_1} is a direct generalization of the IMEX method presented in \cite{Fu2022}. The approach outlined in \cite{Fu2022} requires solving a fourth-order polynomial equation for each spatial grid, offering computational efficiency akin to an explicit scheme while benefiting from a time step at the same order of the mesh size. The advection terms of Eq. \eqref{eq:iter_full_1_next} are updated by the results obtained in \eqref{eq:iter_full_1}, while only a local nonlinear system has to be solved at each grid point. 

{\small\begin{subequations}
\label{eq:iter_full_1}
\begin{align}
\label{eq:iter_full_1_1}
    &  \frac{2\epsilon^2}{c}\frac{\rho^{n+1,\ast, 2k+1} - \rho^n}{\Delta t} + \epsilon \left(\frac{\partial \psi_1}{\partial x}\right)^{n+1, \ast, 2k} = -\frac{\sigma^{n+1, \ast, 2k}}{2}(2\rho - acT_e^4)^{n+1,\ast, 2k+1},\\
    & \frac{\epsilon^2}{c}\frac{\psi_1^{n+1,\ast, 2k+1}-\psi_1^n}{\Delta t} +\frac{2\epsilon}{3}\left(\fpp{\rho}{x}\right)^{n+1,\ast, 2k+1}+\frac{2\epsilon}{3}\left(\fpp{\psi_2}{x}\right)^{n+1,\ast, 2k}=-\sigma^{n+1,\ast, 2k+1}\psi_1^{n+1,\ast, 2k+1},\\
\label{eq:iter_full_1_2}
    &  \frac{\epsilon^2}{c} \frac{\bbu^{n+1,\ast, 2k+1} -\bbu^{n}}{\Delta t} + \epsilon \bA\left(\fpp{\bu}{x}\right)^{n+1,\ast, 2k+1} +\epsilon \bB\left(\fpp{\bu}{x}\right)^{n+1, \ast, 2k} = -\sigma^{n+1,\ast, 2k+1}\bbu^{n+1,\ast, 2k+1} \\
\label{eq:iter_full_1_3}
    & \epsilon^2 C_{v_e}\frac{T_e^{n+1,\ast, 2k+1} - T_e^{n}}{\Delta t}  = \frac{1}{2} c \kappa (T_i - T_e)^{n+1,\ast, 2k+1} + \frac{\sigma^{n+1, \ast, 2k}}{2}(2\rho-ac T_e^4)^{n+1,\ast, 2k+1}, \\
\label{eq:iter_full_1_4}
    &  \epsilon^2 C_{v_i}\frac{T_i^{n+1,\ast, 2k+1} - T_i^{n}}{\Delta t}  = \frac{1}{2} c \kappa (T_e - T_i)^{n+1,\ast, 2k+1},
\end{align}
\end{subequations}}

{\small\begin{subequations}
\label{eq:iter_full_1_next}
\begin{align}
    \label{eq:iter_full_1_next_1}
    &  \frac{2\epsilon^2}{c}\frac{\rho^{n+1,\ast, 2k+2} - \rho^n}{\Delta t} + \epsilon \left(\frac{\partial \psi_1}{\partial x}\right)^{n+1, \ast, 2k+2} = -\frac{\sigma^{n+1, \ast, 2k+1}}{2}\bigg(2\rho^{n+1,\ast, 2k+2} - ac(T_e^{n+1,\ast, 2k+1})^3T_e^{n+1,\ast, 2k+2}\bigg),\\
        \label{eq:iter_full_1_next_2}
    &      \frac{ \epsilon^2}{c} \frac{\psi_1^{n+1, \ast, 2k+2} -\psi_1^{n}}{\Delta t} + \frac{2\epsilon}{3} \left(\fpp{\rho}{x}\right)^{n+1, \ast, 2k+2} +\frac{2\epsilon}{3}\left(\fpp{\psi_2}{x}\right)^{n+1,\ast, 2k+1} = -\sigma^{n+1,\ast, 2k+1}\psi_1^{n+1, \ast, 2k+2} \\
        \label{eq:iter_full_1_next_3}
    &      \frac{ \epsilon^2}{c} \frac{\bbu^{n+1, \ast, 2k+2} -\bbu^{n}}{\Delta t} + \epsilon \bA\left(\fpp{\bu}{x}\right)^{n+1, \ast, 2k+2} +\epsilon \bB\left(\fpp{\bu}{x}\right)^{n+1, \ast, 2k+1} = -\sigma^{n+1,\ast, 2k+2}\bbu^{n+1, \ast, 2k+2} \\
        \label{eq:iter_full_1_next_4}
    & \epsilon^2 C_{v_e}\frac{T_e^{n+1, \ast, 2k+2} - T_e^{n}}{\Delta t}  = \frac{1}{2} c \kappa (T_i - T_e)^{n+1, \ast, 2k+2} + \frac{\sigma^{n+1, \ast, 2k+1}}{2}\bigg(2\rho^{n+1, \ast, 2k+2}-ac (T_e^{n+1, \ast, 2k+1})^3T_e^{n+1, \ast, 2k+2}\bigg), \\
        \label{eq:iter_full_1_next_5}
    &  \epsilon^2 C_{v_i}\frac{T_i^{n+1, \ast, 2k+2} - T_i^{n}}{\Delta t}  = \frac{1}{2} c \kappa (T_e - T_i)^{n+1, \ast, 2k+2},
\end{align}
\end{subequations}}%
where $()^{n+1, \ast, 2k+1}$ is the numerical solution after the step \eqref{eq:iter_full_1} in one iteration. 
When solving \eqref{eq:iter_full_1},  the governing equation of $T_e$ is obtained first. Specifically, by substituting the governing equation of $\rho$ \eqref{eq:iter_full_1_1} and $T_i$ \eqref{eq:iter_full_1_4} into that of $T_e$ \eqref{eq:iter_full_1_3} in \eqref{eq:iter_full_1}, a fourth-order polynomial equation of $T_e^{n+1, \ast, 2k+1}$ is derived as  
\begin{equation}
\label{eq:funTe}
    \mathscr{C}_4 (T_e^{n+1,\ast, 2k+1})^4 + \mathscr{C}_1(T_e^{n+1,\ast, 2k+1}) + \mathscr{C}_0 = 0.
\end{equation}
Here, the coefficients are given by 
\begin{equation}
\label{eq:coeT}
    \begin{gathered}
    \mathscr{C}_4 = \frac{\sigma^{n+1, \ast, 2k}ac\Delta t}{2\epsilon^2+\sigma^{n+1,\ast, 2k}c\Delta t},\qquad
    \mathscr{C}_1 = C_{v_e}+\frac{C_{v_i}c\kappa \Delta t}{2\epsilon^2C_{v_i}+c\kappa\Delta t},\\
    \mathscr{C}_0 = -C_{v_e}T_e^n-\frac{C_{v_i}c\kappa \Delta t}{2\epsilon^2C_{v_i}+c\kappa\Delta t}T_i^n-\frac{2\sigma^{n+1,\ast, 2k}\Delta t}{2\epsilon^2+\sigma^{n+1,\ast, 2k}\Delta t c}\rho^n+\frac{\sigma^{n+1,\ast, 2k}(\Delta t)^2c}{\epsilon(2\epsilon^2+\sigma^{n+1,\ast, 2k}\Delta tc)}\left(\fpp{\psi_1}{x}\right)^{n+1,\ast, 2k}.
    \end{gathered}
\end{equation}
This equation is a key component of the iterative method and allows for the efficient update of the electron temperature $T_e$. The existence and uniqueness of the solution to the fourth-order polynomial equation \eqref{eq:funTe}, with the coefficients provided in \eqref{eq:coeT}, are ensured as long as $\rho^{n+1,\ast, 2k+1}$, $T_e^{n+1, \ast, 2k+1}$, and $T_i^{n+1, \ast, 2k+1}$ exist. In fact, let
\begin{equation}
    f(T) = \mathscr{C}_4 T^4 + \mathscr{C}_1 T + \mathscr{C}_0.    
\end{equation}
Then, it is easy to verify that $f'(T) > 0$ for all $T>0$ and $f(0) < 0$. Thus, the positive solution of \eqref{eq:funTe} is unique. 

Once $T_e^{n+1, \ast, 2k+1}$ is obtained, $\rho^{n+1, \ast, 2k+1}, \sigma^{n+1, \ast, 2k+1}$ and $T_i^{n+1, \ast, 2k+1}$ can be directly calculated using \eqref{eq:iter_full_1_1}, \eqref{eq:iter_full_1_3} and \eqref{eq:iter_full_1_4}. The high-order terms, $\psi_l^{n+1,\ast, 2k+1}, l = 1, \cdots, M$, can be solved successively. Notably, this numerical scheme can be solved with the same computational cost as an explicit scheme, even though the dominant terms are treated implicitly. For example, when updating $\psi_1^{n+1,\ast, 2k+1}$, the implicit term is $\rho^{n+1, \ast, 2k+1}$, which has already obtained by \eqref{eq:iter_full_1_1}.  
When solving the system \eqref{eq:iter_full_1_next}, a linear system of \eqref{eq:iter_full_1_next_1}, \eqref{eq:iter_full_1_next_2}, \eqref{eq:iter_full_1_next_4} and \eqref{eq:iter_full_1_next_5} is first solved to obtain $\rho^{n+1, \ast, 2k+2}, \psi_1^{n+1, \ast, 2k+2}, T_e^{n+1, \ast, 2k+2}$ and $T_i^{n+1, \ast, 2k+2}$. Then, the high-order terms, $\psi_l^{n+1, \ast, 2k+2}, l = 2, \cdots, M$, are solved successively as in \eqref{eq:iter_full_1_2}. 

The alternative iteration method is also adopted for the temporal discretization of the macroscopic part \eqref{eq:PN_full_2}. Assuming the numerical solution after the iteration of the first part is $()^{n+1,\ast}$, the iterative method for the macroscopic part is displayed in Eqs. \eqref{eq:iter_full_2} and \eqref{eq:iter_full_2_next}. In Eq. \eqref{eq:iter_full_2}, the interaction terms on the right side are treated implicitly to allow for larger time step length, while the conduction terms including the coefficients such as the diffusion coefficients $D_s, s = e, i$ and the opacity $\sigma$ are treated explicitly to simplify the iteration process. Eq. \eqref{eq:iter_full_2} solves a local nonlinear system at each spatial grid where the governing equation of $T_e$, which is a fourth-order polynomial equation, is initially solved. Subsequently, in Eq. \eqref{eq:iter_full_2_next}, the conduction coefficients are updated by the results obtained in Eq. \eqref{eq:iter_full_2}, and two linear parabolic equations are solved. 

\begin{subequations}
\label{eq:iter_full_2}
\begin{align}
\label{eq:iter_full_2_1}
    & \frac{2\epsilon^2}{c} \frac{\rho^{n+1,2k+1} - \rho^{n+1, \ast}}{\Delta t} = - \frac{\sigma^{n+1, 2k}}{2} \left(2\rho - acT_e^4 \right)^{n+1,2k+1}, \\
\label{eq:iter_full_2_2}
    &\begin{aligned}
    &  \epsilon^2 C_{v_{e}}\frac{T_e^{n+1,2k+1} - T_e^{n+1, \ast}}{\Delta t} = \epsilon^2 \fpp{}{x} \left(D_e^{n+1, 2k} \left(\fpp{T_e}{x} \right)^{n+1, 2k} \right)+ \frac{1}{2} c \kappa (T_i - T_e)^{n+1,2k+1} \\&\qquad\qquad\qquad\qquad+ \frac{\sigma^{n+1, 2k}}{2}(2\rho - ac T_e^4)^{n+1,2k+1}, 
    \end{aligned}\\
\label{eq:iter_full_2_3}
    &  \epsilon^2 C_{v_{i}} \frac{T_i^{n+1,2k+1} - T_i^{n+1, \ast}}{\Delta t} =\epsilon^2 \fpp{}{x} \left(D_i^{n+1, 2k} \left(\fpp{T_i}{x}\right)^{n+1, 2k} \right)+   \frac{1}{2} c \kappa (T_e - T_i)^{n+1,2k+1},
\end{align}
\end{subequations}

\begin{subequations}
\label{eq:iter_full_2_next}
\begin{align}
    & \frac{2\epsilon^2}{c} \frac{\rho^{n+1, 2k+2} - \rho^{n+1, \ast}}{\Delta t} = - \frac{\sigma^{n+1,2k+1}}{2} \left(2\rho^{n+1,2k+2} - ac(T_e^{n+1,2k+1})^3T_e^{n+1,2k+2} \right), \\
    &\begin{aligned}
    &  \epsilon^2 C_{v_{e}}\frac{ T_e^{n+1, 2k+2} - T_e^{n+1, \ast}}{\Delta t} = \epsilon^2 \fpp{}{x} \left(D_e^{n+1,2k+1} \left(\fpp{T_e}{x} \right)^{n+1, 2k+2} \right)+ \frac{1}{2} c \kappa (T_i - T_e)^{n+1, 2k+2} \\&\qquad\qquad\qquad\qquad+ \frac{\sigma^{n+1,2k+1}}{2}(2\rho^{n+1,2k+2} - ac (T_e^{n+1,2k+1})^3T_e^{n+1,2k+2}), 
    \end{aligned}\\
    &  \epsilon^2 C_{v_{i}} \frac{T_i^{n+1, 2k+2} - T_i^{n+1, \ast}}{\Delta t} =\epsilon^2 \fpp{}{x} \left(D_i^{n+1,2k+1} \left(\fpp{T_i}{x}\right)^{n+1, 2k+2} \right)+   \frac{1}{2} c \kappa (T_e - T_i)^{n+1, 2k+2}.
\end{align}
\end{subequations}
However, using solely \eqref{eq:iter_full_2_next} is not sufficient for the iteration process involving  the macroscopic part \eqref{eq:PN_full_2}. This limitation arises because the linearization of the term $T^4$ in Eq. \eqref{eq:iter_full_2_next}  will no longer match the density $\rho$ as $\epsilon$ diminishes. Consequently, the alternating iteration method is adopted, wherein the numerical solution $()^{n+1,2k+1}$ serves as a precondition. Additionally, numerical observations reveal that in this alternating iteration method, the time step length remains independent of $\epsilon$ and maintains the same order as the spatial step. Besides, when $D_e = D_i = 0$, \eqref{eq:PN_full_2} reduces to an ODEs system, greatly simplifying the numerical solution of the system. 



\section{Spatial discretization and the AP property}
\label{sec:IMEX}

In this section, the spatial discretization of the system \eqref{eq:PN_full_1} and \eqref{eq:PN_full_2} is introduced. The finite volume method is applied to \eqref{eq:PN_full_1}, while the finite difference method is used for \eqref{eq:PN_full_2}. The AP property and the energy conservation properties of the fully discretized scheme are discussed.

\subsection{Spatial discretization}
\label{sec:full-dis}

The spatial discretization for the model \eqref{eq:iter_full_1}, \eqref{eq:iter_full_1_next}, \eqref{eq:iter_full_2} and \eqref{eq:iter_full_2_next} in the slab geometry is introduced, and the two-dimensional extension is straightforward. Let the computational domain be $[0, 1]$ and $\Delta x  = \frac{1}{N}$. $x_{j-1/2} = (j-1/2) \Delta x, j = 1, \cdots, N$  represent the uniform mesh in the Cartesian coordinates. Let $\Gamma_j$ denote the cell $[x_{j-1/2}, x_{j+1/2}]$. $T_{e, j}^n$, $T_{i, j}^n$, and $\bu_j^n$ represent the average of $T_e, T_i$ and $\bu$ at time $t^n$ inside cell $\Gamma_j$, respectively. The finite volume method is utilized to approximate the convection term in \eqref{eq:iter_full_1} and \eqref{eq:iter_full_1_next}. The revised Lax-Friedrichs scheme described in \cite{Fu2022} is adopted to determine the numerical flux. 
\begin{subequations}
\label{eq:flux}
\begin{align}
    &\frac{l+1}{2l+1}\left(\fpp{\psi_{l+1}}{x} \right)_j^{\rm old} \approx \frac{1}{\Delta x} \left(\mF_{l+1}^{\rm old}\left(\bu_{j+1/2}^{L,\rm old},\bu_{j+1/2}^{R,\rm old}\right) - \mF_{l+1}^{\rm old}\left(\bu_{j-1/2}^{L,\rm old},\bu_{j-1/2}^{R,\rm old}\right)\right),\ l = 0, 1, \cdots, M-1,\\
    &\frac{l}{2l+1}\left(\fpp{\psi_{l-1}}{x} \right)_j^{\rm new} \approx \frac{1}{\Delta x} \left(\mF_{l-1}^{\rm new}\left(\bu_{j+1/2}^{L,\rm new},\bu_{j+1/2}^{R,\rm new}\right) - \mF_{l-1}^{\rm new}\left(\bu_{j-1/2}^{L,\rm new},\bu_{j-1/2}^{R,\rm new}\right)\right),\quad l = 1, 2, \cdots, M.
\end{align}
\end{subequations}
where $()^{\rm new}$ is $()^{n+1,\ast,2k+1}$ and $()^{\rm old}$ is $()^{n+1,\ast,2k}$ for the numerical flux in \eqref{eq:iter_full_1}, while $()^{\rm new}$ is $()^{n+1,\ast,2k+2}$ and $()^{\rm old}$ is $()^{n+1,\ast,2k+1}$ for the numerical flux in \eqref{eq:iter_full_1_next}.
The exact forms of $\mF_{l}^{\rm old}$ and $\mF_{l}^{\rm new}$ are given by

{\small
\begin{subequations}
\label{eq:flux_f}
\begin{align}
    &\mF_{l+1}^{\rm old}(\bu^L_{j+\frac12},\bu^R_{j+\frac12})=\frac{l+1}{2l+1}\frac{1}{2}(\psi_{l+1, j+\frac12}^L+\psi_{l+1, j+\frac12}^R)-\frac{\alpha_{j+\frac12}}{2}(\psi_{l,j+\frac12}^R-\psi_{l, j+\frac12}^L),\quad l=0,1,\dots,M-1,\\
    &\mF_{l-1}^{\rm new}(\bu^L_{j+\frac12},\bu^R_{j+\frac12})=\frac{l}{2l+1}\frac{1}{2}(\psi_{l-1,j+\frac12}^L+\psi_{l-1,j+\frac12}^R),\quad l=1,2,\dots,M-1,\\
    &\begin{aligned}
    &\mF_{M-1}^{\rm new}(\bu^{L}_{j+\frac12},\bu^{R}_{j+\frac12})=\frac{M}{2M+1}\frac{1}{2}(\psi_{M-1,j+\frac12}^{L}+\psi_{M-1,j+\frac12}^{R})-\frac{\alpha_{j+\frac12}}{2}(\psi_{M,j+\frac12}^{R}-\psi_{M,j+\frac12}^{L}),
    \end{aligned}
\end{align}
\end{subequations}}%
with
\begin{equation}
\label{eq:flux_coe}
    \alpha_{j+\frac12} = \alpha(\sigma, \epsilon) = \frac{c}{\epsilon^2}\exp\left(-\frac{\sigma_j + \sigma_{j+1}}{2 \epsilon^2}\right),
\end{equation}
and $\bu^L$ and $\bu^R$ are given by reconstruction. The simplest way is
\begin{equation}
    \bu_{j+\frac12}^L=\bu_j,\quad \bu_{j+\frac12}^R=\bu_{j+1}.
\end{equation}
To improve the resolution of the solution, linear reconstruction or WENO reconstruction can be used.

\begin{remark}
The different forms of flux are due to the non-consistent governing equation of $\psi_{M}$. The choice of $\alpha_j$ is necessary to ensuring the AP property of the numerical scheme and we refer to \cite{Fu2022} for more details. 
\end{remark}

For the diffusion terms in \eqref{eq:iter_full_2} and \eqref{eq:iter_full_2_next}, the finite difference method is employed and is expressed as 
\begin{equation}
\label{eq:diff_sch}
    \begin{aligned}
    & \fpp{}{x}\left(D_{s}^{p} \left(\fpp{T_{s}}{x}\right)_j^{q}\right) \approx 
      \frac{1}{\Delta x} \left( \frac{D_{s, j+1}^{p} + D_{s, j}^{p}}{2} \frac{T_{s, j+1}^{q} -T_{s, j}^{q} }{\Delta x} - \frac{D_{s, j}^{p} + D_{s, j-1}^{p}}{2} \frac{T_{s, j}^{q} -T_{s, j-1}^{q} }{\Delta x}\right), \quad s = e, i, 
    \end{aligned}  
\end{equation}
where $p, q$ is the label for the time step.

\subsection{AP property}
\label{sec:ap_analy} 

To demonstrate the AP property, it is necessary to show that for a fixed mesh size and time step, the numerical scheme should automatically converge to a scheme for the diffusion limit \eqref{eq:rtelimit} as the parameter $\epsilon$ approaches zero. Simultaneously, it should also converge to a scheme for the two-temperature limit \eqref{eq:limit_kappa} as the parameter $\kappa$ tends to infinity. In the subsequent part, we illustrate the AP property of the numerical scheme presented in Sec. \ref{sec:semi-dis} and Sec. \ref{sec:full-dis}. 

\paragraph{\bf AP property to the diffusion limit}
{\it As the parameter $\epsilon$ approaches zero, the numerical scheme described in Sec. \ref{sec:full-dis} converges to a consistent scheme for the nonlinear diffusion equation \eqref{eq:rtelimit}.} 

The fully discretized scheme for the $\rho, T_e$ and $T_i$ are 
\begin{equation}
\label{eq:ap0}
    \begin{aligned}
    & \frac{2\epsilon^2}{c}\frac{\rho^{n+1, \ast}_j - \rho^n_j}{\Delta t} + \frac{\mF_{1, j+1/2}^{n+1, \ast} - \mF_{1, j-1/2}^{n+1, \ast}}{\Delta x} = -\frac{\sigma^{n+1, \ast}_j}{2}(2\rho - acT_e^4)_j^{n+1, \ast},\\
    & \epsilon^2 C_{v_e}\frac{T_{e,j}^{n+1, \ast} - T_{e,j}^{n}}{\Delta t}  = \frac{1}{2} c \kappa (T_i - T_e)^{n+1, \ast}_j + {\frac{\sigma_j^{n+1,\ast}}{2}}( 2\rho-ac T_e^4)^{n+1, \ast}_j, \\
    & \epsilon^2 C_{v_i}\frac{T_{i,j}^{n+1, \ast} - T_{i,j}^{n}}{\Delta t}  = \frac{1}{2} c \kappa (T_e - T_i)^{n+1, \ast}_j,
    \end{aligned}
\end{equation}
where
\begin{equation}
\label{eq:ap3}
    \mF_{1,j+1/2}^{n+1, \ast} = \frac{\epsilon}{2}(\psi_{1,j}^{n+1, \ast} + \psi_{1, j+1}^{n+1, \ast}) - \alpha_{j+1/2}^{n+1, \ast}\epsilon (\rho_{j+1}^{n+1,\ast}-\rho_j^{n+1,\ast}).
\end{equation}
Summing up the equations in \eqref{eq:ap0}, we have 
\begin{equation}
\label{eq:ap1}
    \begin{aligned}
     \frac{2\epsilon^2}{c} \frac{\rho_{j}^{n+1,\ast} - \rho_{j}^n}{\Delta t}  + \epsilon^2 C_{v_e} \frac{T_{e,j}^{n+1,  \ast} - T_{e,j}^{n}}{\Delta t} + \epsilon^2C_{v_i} \frac{T_{i,j}^{n+1, \ast} - T_{i,j}^{n}}{\Delta t} = -\frac{\mF_{1,j+1/2}^{n+1, \ast} - \mF_{1,j-1/2}^{n+1, \ast}}{\Delta x}.
    \end{aligned}
\end{equation}
The equation of $\psi_1$ in \eqref{eq:PN_full_1} is 
\begin{equation}
\label{eq:ap2}
    \frac{\epsilon^2}{c} \frac{\psi_{1,j}^{n+1, \ast} - \psi_{1,j}^n}{\Delta t}  + \frac{\mF_{0, j+1/2}^{n+1, \ast} - \mF_{0, j-1/2}^{n+1, \ast} }{\Delta x} + \frac{\mF_{2, j+1/2}^{n+1, \ast} - \mF_{2, j-1/2}^{n+1, \ast} }{\Delta x} = -\sigma_j^{n+1, \ast} \psi_{1,j}^{n+1, \ast},
\end{equation}
where
\begin{equation}
    \mF_{0,j+1/2}^{n+1, \ast} = \frac{2}{3}\frac{\epsilon}{2}(\rho_{j}^{n+1, \ast} + \rho_{j+1}^{n+1, \ast}),\quad \mF_{2,j+1/2}^{n+1, \ast} = \frac{2}{3}\frac{\epsilon}{2}(\psi_{2,j}^{n+1, \ast} + \psi_{2, j+1}^{n+1, \ast})-\alpha_{j+1/2}^{n+1,\ast}\frac{\epsilon}{2}(\psi_{1,j+1}^{n+1,\ast}-\psi_{1,j}^{n+1,\ast}).
\end{equation}
In the limit as $\epsilon$ approaches zero, order analysis (as detailed in \ref{app:order_analy}) yields 
\begin{equation}
\label{eq:ap4}
    \psi_{1,j}^{n+1,\ast} = \mO(\epsilon), \qquad   \psi_{2,j}^{n+1,\ast} = \mO(\epsilon^2). 
\end{equation}
With this, we can rewrite the source term in \eqref{eq:ap2} as 
\begin{equation}
\label{eq:ap5}
    -\sigma^{n+1, \ast}_j \psi_{1, j}^{n+1, \ast} = \frac{2\epsilon}{ 3}\frac{\rho_{j+1}^{n+1, \ast}- \rho_{j-1}^{n+1, \ast}}{2\Delta x} + \mO(\epsilon^2).
\end{equation}
Substituting \eqref{eq:ap5} into \eqref{eq:ap3}, one obtains
\begin{equation}
\label{eq:ap6}
    \mF^{n+1, \ast}_{1, j+1/2} =\frac{ \epsilon^2}{2} \left( \frac{2}{3 \sigma_j^{n+1, \ast}} \frac{\rho_{j+1}^{n+1, \ast}- \rho_{j-1}^{n+1, \ast}}{2\Delta x} + \frac{2}{3 \sigma_{j+1}^{n+1, \ast}} \frac{\rho_{j+2}^{n+1, \ast}- \rho_{j}^{n+1, \ast}}{2\Delta x} \right)  + \mathcal{O}(\epsilon^3).
\end{equation}
With equation \eqref{eq:ap6}, equation \eqref{eq:ap1} is reduced to
\begin{equation}
\label{eq:ap7}
\begin{aligned}
    & \frac{2\epsilon^2}{c} \frac{\rho_{j}^{n+1,\ast} - \rho_{j}^n}{\Delta t}  + \epsilon^2 C_{v_e} \frac{T_{e,j}^{n+1,  \ast} - T_{e,j}^{n}}{\Delta t} + \epsilon^2 C_{v_i} \frac{T_{i,j}^{n+1, \ast} - T_{i,j}^{n}}{\Delta t} \\
    & = \frac{\epsilon^2}{2\Delta x}\left(\frac{2}{3 \sigma_{j+1}^{n+1, \ast}} \frac{\rho_{j+2}^{n+1, \ast}- \rho_{j}^{n+1, \ast}}{2\Delta x}  -  \frac{2}{3 \sigma_{j-1}^{n+1, \ast}} \frac{\rho_{j}^{n+1, \ast}- \rho_{j-2}^{n+1, \ast}}{2\Delta x}\right)+\mathcal{O}(\epsilon^3).
\end{aligned}
\end{equation}
As $\epsilon$ tends to zero, the dominance of the sources in \eqref{eq:ap0} leads to  
\begin{equation}
\label{eq:ap8}
    T_{e,j}^{n+1,\ast} = (T^{(0)})_j^{n+1,\ast}+\mathcal{O}(\epsilon), \qquad T_{i,j}^{n+1,\ast} = (T^{(0)})_j^{n+1,\ast}+\mathcal{O}(\epsilon), \qquad  2\rho_j^{n+1,\ast} = ac ((T^{(0)})^4)_j^{n+1,\ast}+\mathcal{O}(\epsilon).
\end{equation}
By substituting \eqref{eq:ap8} into \eqref{eq:ap7}, it holds that, as $\epsilon$ approaches zero, 
\begin{equation}
\label{eq:ap9}
    \begin{aligned}
    &  \frac{a ((T^{(0)})^4)_j^{n+1,\ast} - a((T^{(0)})^4)_{j}^n}{\Delta t} +  \frac{C_{v_e} (T^{(0)}_j)^{n+1,  \ast} - C_{v_e} (T^{(0)}_j)^{n}}{\Delta t} + \frac{C_{v_i} (T^{(0)}_j)^{n+1, \ast} - C_{v_i} (T^{(0)}_j)^{n}}{\Delta t} \\
    & = \frac{1}{2\Delta x}\left(\frac{ac}{3 \sigma((T^{(0)})_{j+1}^{n+1, \ast})} \frac{((T^{(0)})^4)_{j+2}^{n+1, \ast}- ((T^{(0)})^4)_{j}^{n+1, \ast}}{2\Delta x}  - \frac{ac}{3 \sigma((T^{(0)})_{j-1}^{n+1, \ast})} \frac{((T^{(0)})^4)_{j}^{n+1, \ast}- ((T^{(0)})^4)_{j-2}^{n+1, \ast}}{2\Delta x}\right),
\end{aligned}
\end{equation}
which is a numerical scheme for the nonlinear diffusion equation \eqref{eq:rtelimit}. We have completed the proof of the first part. 

In the second part, as $\epsilon$ tends to zero in \eqref{eq:PN_full_2}, it holds that
\begin{equation}
\label{eq:ap11}
    T_{e,j}^{n+1} = (T^{(0)})_j^{n+1}+\mathcal{O}(\epsilon),\qquad  T_{i,j}^{n+1} = (T^{(0)})_j^{n+1}+\mathcal{O}(\epsilon), \qquad 2\rho_j^{n+1}=(ac(T^{(0)})^4)_j^{n+1}+\mathcal{O}(\epsilon),
\end{equation}
Summing up the equations for $\rho$, $T_e$, and  $T_i$ in \eqref{eq:PN_full_2} and employing \eqref{eq:diff_sch}, we can derive that 
{\small \begin{equation}
\label{eq:ap10}
    \begin{aligned}
    &   \frac{a ((T_j^{(0)})^4)_j^{n+1} - a ((T_j^{(0)})^4)_{j}^{n+1, \ast}}{\Delta t} +  \frac{C_{v_e} (T^{(0)}_{j})^{n+1} - C_{v_e} (T^{(0)}_{j})^{n+1, \ast}}{\Delta t} + \frac{C_{v_i} (T^{(0)}_{j})^{n+1} - C_{v_i}(T^{(0)}_{j})^{n+1, \ast}}{\Delta t}\\
    &  \qquad = \frac{1}{\Delta x} \left( \frac{D_e((T^{(0)})_{j+1}^{n+1}) + D_e((T^{(0)})_{j}^{n+1})}{2} \frac{(T^{(0)})_{j+1}^{n+1} -(T^{(0)})_{j}^{n+1} }{\Delta x} - \frac{D_e((T^{(0)})_{j}^{n+1}) + D_e((T^{(0)})_{j-1}^{n+1})}{2} \frac{(T^{(0)})_{j}^{n+1} -(T^{(0)})_{j-1}^{n+1} }{\Delta x}  \right) \\
    & \qquad +\frac{1}{\Delta x} \left( \frac{D_i((T^{(0)})_{j+1}^{n+1}) + D_i((T^{(0)})_{j}^{n+1})}{2} \frac{(T^{(0)})_{j+1}^{n+1} -(T^{(0)})_{j}^{n+1} }{\Delta x} - \frac{D_i((T^{(0)})_{j}^{n+1}) + D_i((T^{(0)})_{j-1}^{n+1})}{2} \frac{(T^{(0)})_{j}^{n+1} -(T^{(0)})_{j-1}^{n+1} }{\Delta x}
    \right). 
    \end{aligned}
\end{equation}}%
Considering \eqref{eq:ap9}, \eqref{eq:ap11}, and \eqref{eq:ap10}, we observe that the numerical scheme \eqref{eq:PN_full}, with the spatial discretization \eqref{eq:flux} and \eqref{eq:diff_sch}, reduces to an implicit scheme of the diffusion limit \eqref{eq:rtelimit}. Furthermore, \eqref{eq:PN_full_1} and \eqref{eq:PN_full_2} converge to the numerical schemes corresponding to the 1D form of the split systems \eqref{eq:MM_limit_split_2_1} and \eqref{eq:MM_limit_split_2_2}. This convergence completes the proof. 

\paragraph{\bf AP property to the two temperature limit}
{\it As the parameter $\kappa$ tends to infinity, the numerical scheme proposed in Sec. \ref{sec:full-dis} converges to a consistent numerical scheme for the two-temperature model given by \eqref{eq:limit_kappa}. }

As $\kappa$ approaches infinity, it holds that 
\begin{equation}
    T_{e,j}^{n+1, \ast}= T_j^{n+1,\ast} + \mathcal{O}(\kappa^{-1}),\quad T_{i,j}^{n+1,\ast}= T_j^{n+1,\ast} + \mathcal{O}(\kappa^{-1}).
\end{equation}
With $C_v = C_{v_e} + C_{v_i}$, by summing up the governing equation for $T_e$ and $T_i$ in \eqref{eq:PN_full_1}, we can derive 
\begin{equation}
\label{eq:kappa_ap_1} 
    \begin{aligned}
    &   \frac{2\epsilon^2}{c}\frac{\rho_j^{n+1, \ast} - \rho_j^n}{\Delta t} + \frac{\mF^{n+1, \ast}_{1, j+1/2} - \mF^{n+1, \ast}_{1, j-1/2}}{\Delta x} = -\frac{\sigma_j^{n+1, \ast}}{2}(2\rho - acT^4)_j^{n+1, \ast},\\
    & \frac{\epsilon^2}{c}\frac{\psi_{1,j}^{n+1,\ast}-\psi_{1,j}^n}{\Delta t}+\frac{\mathcal{F}_{0,j+1/2}^{n+1,\ast}-\mathcal{F}_{0,j-1/2}^{n+1,\ast}}{\Delta x}+\frac{\mathcal{F}_{2,j+1/2}^{n+1,\ast}-\mathcal{F}_{2,j-1/2}^{n+1,\ast}}{\Delta x}=-\sigma_j^{n+1,\ast}\psi_{1,j},\\
    & \frac{\epsilon^2}{c} \frac{\bbu_j^{n+1, \ast} -\bbu_j^{n}}{\Delta t} +  \frac{{\bf F}^{n+1, \ast}_{1, j+1/2} - {\bf F}^{n+1, \ast}_{1, j-1/2}}{\Delta x} = -\sigma_j^{n+1,\ast}\bbu_j^{n+1, \ast}, \\
    &  \epsilon^2 C_v \frac{T_j^{n+1, \ast} - T_j^{n}}{\Delta t} = \frac{\sigma_j^{n+1, \ast}}{2} (2\rho - ac T^4)_j^{n+1, \ast}.
    \end{aligned}
\end{equation}
Here, $\bf F$ is a vector with $M-1$ entries, and the $l$-th entry represents the corresponding numerical flux of $\psi_{l+1}, l = 1, \cdots, M-1$. When $\kappa$ approach infinity in \eqref{eq:PN_full_2}, one obtains
\begin{equation}
    T_{e,j}^{n+1}= T_j^{n+1} + \mathcal{O}(\kappa^{-1}),\quad T_{i,j}^{n+1}= T_j^{n+1} + \mathcal{O}(\kappa^{-1}).
\end{equation}
By summing up the governing equations for $T_e$ and $T_i$ in \eqref{eq:PN_full_2}, we have
\begin{equation}
\label{eq:kappa_ap_3}
    \begin{aligned}
    & \frac{2\epsilon^2}{c} \frac{\rho_j^{n+1} - \rho_j^{n+1, \ast}}{\Delta t} = - \frac{\sigma_j^{n+1}}{2} \left(2\rho - acT^4 \right)_j^{n+1}, \\
    &  \epsilon^2 C_{v}\frac{ T_j^{n+1} - T_j^{n+1, \ast}}{\Delta t} = 
    \frac{\epsilon^2}{\Delta x} \left( \frac{D_{T, j+1}^{n+1} + D_{T, j}^{n+1}}{2} \frac{T_{j+1}^{n+1} -T_{j}^{n+1} }{\Delta x} \right. \\
    & \qquad\qquad \qquad \qquad  \left. - \frac{D_{T, j}^{n+1} + D_{T, j-1}^{n+1}}{2} \frac{T_{j}^{n+1} -T_{j-1}^{n+1} }{\Delta x} \right)  +\frac{\sigma^{n+1}}{2} \left(2\rho - acT^4 \right)_j^{n+1},
    \end{aligned}
\end{equation}
with $ C_v = C_{v_e} + C_{v_i}$, and $D_T = D_e(T) + D_i(T)$. It is evident that \eqref{eq:kappa_ap_1} and \eqref{eq:kappa_ap_3} together constitute a splitting scheme of the two-temperature model \eqref{eq:limit_kappa}, thus completing the proof.

\subsection{Conservation of total energy}
\label{sec:con_energy}

In this section, we will discuss the energy-conserving property of this numerical scheme. For the 3TRTE model \eqref{eq:rte_1D}, it conserves the total energy, as demonstrated by the theorem below.

\paragraph{\bf Energy-conserving property}
{\it For the 3TRTE system \eqref{eq:rte_1D}, with the periodic boundary conditions in the spatial space and the absence of external source terms, the system conserves total energy. This conservation is precisely expressed by 
\begin{equation}
\label{eq:conser_E}
    \fpp{}{t} E_{\rm total}(t) = 0,
\end{equation}
where the total energy $E_{\rm total}(t)$ is defined as 
\begin{equation}
\label{eq:total_E} 
    E_{\rm total}(t) = \int   E_e(t, x) + E_i(t, x) + E_r(t, x) \dd x.
\end{equation} }

The proof is straightforward by integrating \eqref{eq:Pn_energy} in the spatial domain, and we will omit it.  For the numerical scheme \eqref{eq:PN_full}, utilizing the spatial discretization methods presented in \eqref{eq:flux} and \eqref{eq:diff_sch}, it effectively conserves the discretized total energy. This conservation is expressed as follows. 
\paragraph{\bf Conservation of the discretized total energy}
{\it The numerical scheme \eqref{eq:PN_full}, with the spatial discretization given by \eqref{eq:flux} and \eqref{eq:diff_sch}, conserves the total energy 
\begin{equation}
\label{eq:dis_con_E}
    \mE_{\rm total}^{n+1} = \mE_{\rm total}^{n}
\end{equation}
for the 3TRTE system \eqref{eq:rte_1D} with periodic boundary condition. The discretized total energy at time $t= t^n$ is defined as 
\begin{equation}
\label{eq:dis_total_E}
    \mE_{\rm total}^n =  \sum_{j=1}^N( E_{e,j}^n + E_{i,j}^n + E_{r,j}^n), 
\end{equation}
with $N$ the total number of mesh points in the spatial domain.   }

To prove this, first, by summing up the governing equation for $\rho$,  $T_e$, and $T_i$ in \eqref{eq:PN_full}, with the spatial discretization \eqref{eq:flux} and \eqref{eq:diff_sch}, it holds that 
\begin{equation}
\label{eq:dis_total_E_1}
    \begin{aligned}
    & \frac{2\epsilon^2}{c}\frac{\rho_j^{n+1} - \rho_j^n}{\Delta t} + \epsilon^2 C_{v_e} \frac{T_{e,j}^{n+1} - T_{e,j}^{n}}{\Delta t} + \epsilon^2 C_{v_i} \frac{T_{i,j}^{n+1} - T_{i,j}^n}{\Delta t} \\
    & \qquad \qquad   = -\frac{\mF_{1,j+1/2}^{n+1, \ast} - \mF_{1,j-1/2}^{n+1, \ast}}{\Delta x} \\
    & \qquad \qquad + \frac{\epsilon^2}{\Delta x}\left(D^{n+1}_{e, j+1/2} \frac{T^{n+1}_{e, j+1} - T^{e, n+1}_{e,j}}{\Delta x} - D^{n+1}_{e, j-1/2} \frac{T^{n+1}_{e, j} - T^{e, n+1}_{e,j-1}}{\Delta x} \right) \\
    & \qquad \qquad + \frac{\epsilon^2}{\Delta x}\left(D^{n+1}_{i, j+1/2} \frac{T^{n+1}_{i, j+1} - T^{i, n+1}_{i,j}}{\Delta x} - D^{n+1}_{i, j+1/2} \frac{T^{n+1}_{i, j+1} - T^{i, n+1}_{i,j}}{\Delta x} \right),    
    \end{aligned}
\end{equation}
with 
\begin{equation}
    D^{n+1}_{s, j+1/2} =\frac{1}{2}\left( D^{n+1}_{s, j+1} + D^{n+1}_{s, j}\right), \qquad s = e, i.
\end{equation}
Summing up over the spatial space, we obtain 
\begin{equation}
\label{eq:dis_total_E_2}
    \begin{aligned}
    & \sum_{j = 1}^N \left(\frac{2\epsilon^2}{c}\frac{\rho_{j}^{n+1} - \rho_{j}^n}{\Delta t} + \epsilon^2 C_{v_e} \frac{T_{e, j}^{n+1} - T_{e, j}^{n}}{\Delta t} + \epsilon^2 C_{v_i} \frac{T_{i, j}^{n+1} - T_{i, j}^n}{\Delta t}\right)   \\
    & \qquad \qquad = \frac{\mF_{1,1/2}^{n+1, \ast} - \mF_{1,N+1/2}^{n+1, \ast}}{\Delta x} \\
    & \qquad \qquad +  \frac{\epsilon^2}{\Delta x}\left(D^{n+1}_{e, N+1/2} \frac{T^{n+1}_{e, N+1} - T^{n+1}_{e, N}}{\Delta x} - D^{n+1}_{e, 1/2} \frac{T^{n+1}_{e, 1} - T^{n+1}_{e, 0}}{\Delta x} \right) \\
    & \qquad \qquad + \frac{\epsilon^2}{\Delta x}\left(D^{n+1}_{i, N+1/2} \frac{T^{n+1}_{i, N+1} - T^{n+1}_{i, N}}{\Delta x} - D^{n+1}_{i, 1/2} \frac{T^{n+1}_{i, 1} - T^{n+1}_{i, 0}}{\Delta x} \right).     
    \end{aligned}
\end{equation}
With the periodic boundary condition, we have  
\begin{equation}
\label{eq:per}
    \mF_{1,1/2}^{n+1, \ast} = \mF_{1, N+1/2}^{n+1, \ast}, \qquad D^{n+1}_{s, 1/2} = D^{n+1}_{s, N+1/2}, \qquad T^{n+1}_{s, 0} = T^{n+1}_{s, N}, \qquad T^{n+1}_{s, 1} = T^{n+1}_{s, N+1}, \qquad s = e, i.
\end{equation}
Substituting \eqref{eq:per} into \eqref{eq:dis_total_E_2} yields 
\begin{equation}
\label{eq:dis_total_E_3}
    \sum_{j = 1}^N \left(\frac{2\epsilon^2}{c}\frac{\rho_{j}^{n+1} - \rho_{j}^n}{\Delta t} + \epsilon^2 C_{v_e} \frac{T_{e, j}^{n+1} - T_{e, j}^{n}}{\Delta t} 
+ \epsilon^2 C_{v_i} \frac{T_{i, j}^{n+1} - T_{i, j}^n}{\Delta t}\right)  = 0.
\end{equation}
Thus, we have completed the proof.

\subsection{The whole algorithm}
\label{sec:alg}

In this section, for the completeness of the numerical scheme, we provide an outline of the numerical algorithm as follows.

\begin{enumerate}
\item Let $t = 0$ and $n=0$, and calculate the initial value of $\rho_{j}^0$ ,$\psi_j^0$, $\bbu^0_j$, $T_{e,j}^0$, $T_{i, j}^0$;
\item Set time step $\Delta t$ with the CFL condition \eqref{eq:CFL};
\item Solve \eqref{eq:PN_full} to obtain $\rho_{j}^{n+1},\psi_j^{n+1},\bbu^{n+1}_j$, $T_{e, j}^{n+1}$, $T_{i, j}^{n+1}$:
    \begin{enumerate}
    \item Let $k = 0$; 
    \item Set $\rho_{j}^{n+1,\ast, 2k} = \rho_{j}^n, \psi_{1,j}^{n+1,\ast, 2k}=\psi_{1,j}^n, \bbu^{n+1, \ast, 2k}_j = \bbu_j^n$, $T_{e, j}^{n+1,\ast, 2k} = T_{e,j}^{n},T_{i, j}^{n+1,\ast, 2k} =T_{i, j}^{n}$ and $\sigma_j^{n+1,\ast, 2k} = \sigma_j^{n}$;
    \item Obtain $\rho_{j}^{n+1,\ast, 2k+1}, T_{e,j}^{n+1,\ast, 2k+1}, T_{i, j}^{n+1, \ast, 2k+1}$ using \eqref{eq:iter_full_1} and \eqref{eq:flux}; 
    \item Update $\sigma_j^{n+1, \ast, 2k+1}$ using $T_{e, j}^{n+1, \ast, 2k+1}$;
    \item Obtain $\psi_{1,j}^{n+1,\ast, 2k+1},\bbu^{n+1, \ast, 2k+1}_j$ using \eqref{eq:iter_full_1} and \eqref{eq:flux}; 
    \item Obtain $\rho_{j}^{n+1,\ast, 2k+2},\psi_{1,j}^{n+1,\ast, 2k+2}, T_{e,j}^{n+1, \ast, 2k+2}, T_{i, j}^{n+1, \ast, 2k+2}$ using \eqref{eq:iter_full_1_next} and \eqref{eq:flux}; 
    \item Update $\sigma_j^{n+1, \ast, 2k+2}$ using $T_{e, j}^{n+1,\ast, 2k+2}$;
    \item Obtain $\bbu^{n+1, \ast, 2k+2}_j$ using \eqref{eq:iter_full_1_next} and \eqref{eq:flux}; 
    \item Calculate the error $\mathcal{E}$ using \eqref{eq:error_1};
    \item If $\mathcal{E}^{n+1, \ast, 2k+2} < \bar{\epsilon}$, set $\rho_{j}^{n+1,\ast} = \rho_{j}^{n+1, \ast, 2k+2},\psi_{1,j}^{n+1,\ast}=\psi_{1,j}^{n+1,\ast, 2k+2}, \bbu^{n+1,\ast}_j = \bbu_j^{n+1, \ast, 2k+2}$, $T_{e, j}^{n+1,\ast} = T_{e,j}^{n+1, \ast, 2k+2}$ and $T_{i, j}^{n+1,\ast} =T_{i, j}^{n+1, \ast, 2k+2}$, else let $k \leftarrow k+1$, and return to step (c);
    \item Set $\psi_{1,j}^{n+1}=\psi_{1,j}^{n+1,*},\bbu^{n+1}_j = \bbu_j^{n+1, *}$.
    \item Let $k = 0$; 
    \item Set $\rho_{j}^{n+1, 2k} = \rho_{j}^{n+1,\ast},  \bbu^{n+1, 2k}_j = \bbu_j^{n+1,\ast}$, $T_{e, j}^{n+1, 2k} = T_{e,j}^{n+1,\ast},T_{i, j}^{n+1, 2k} =T_{i, j}^{n+1,\ast}$ and $\sigma_j^{n+1, 2k} = \sigma_j^{n+1,\ast}$;
    \item Obtain $\rho_{j}^{n+1,2k+1}, \bbu^{n+1, 2k+1}_j, T_{e,j}^{n+1, 2k+1}, T_{i, j}^{n+1, 2k+1}$ using \eqref{eq:iter_full_2} and \eqref{eq:diff_sch};
    \item Update $\sigma_j^{n+1, 2k+1}$ and $D_{s}^{n+1, 2k+1}, s = e, i$ using $T_{e, j}^{n+1, 2k+1}$ and $T_{i,j}^{n+1, 2k+1}$;
    \item Obtain $\rho_{j}^{n+1,2k+2}, \bbu^{n+1, 2k+2}_j, T_{e,j}^{n+1, 2k+2}, T_{i, j}^{n+1, 2k+2}$ using \eqref{eq:iter_full_2_next} and \eqref{eq:diff_sch};
    \item Update $\sigma_j^{n+1, 2k+2}$ and $D_{s}^{n+1, 2k+2}, s = e, i$ using $T_{e, j}^{n+1, 2k+2}$ and $T_{i,j}^{n+1, 2k+2}$;
    \item Calculate the error $\mathcal{E}$ using \eqref{eq:error_2};
    \item If $\mathcal{E}^{n+1, 2k+2} < \bar{\epsilon}$, set $\rho_{j}^{n+1} = \rho_{j}^{n+1, 2k+2}$, $T_{e, j}^{n+1} = T_{e,j}^{n+1, 2k+2}$ and $T_{i, j}^{n+1} =T_{i, j}^{n+1, 2k+2}$, else let $k \leftarrow k+1$, and return to step (n);
    \end{enumerate}
\item Let $n \leftarrow n+1$, and return to step 2. 
\end{enumerate}

\section{Numerical examples}
\label{sec:num}
In this section, the validation of the numerical method is conducted through a series of tests. First, the AP property of this splitting scheme is verified. Then, numerical examples including spatially homogeneous conditions, the 1D Marshak wave problem, and the 2D Riemann problem are investigated. In the first test, the periodic condition is applied, while the inflow boundary condition in \cite{Fu2022} is applied for the last three tests.

\subsection{Stopping criteria}
\label{sec:iter_stop}

We discuss the stopping criteria for each iteration to solve \eqref{eq:PN_full}. The key stopping criterion is based on the error measure, denoted as $\mathcal{E}$, which ensures the accuracy of the solution. Two forms of this error, $\mathcal{E}^{n+1, \ast. 2k+2}$ and $\mathcal{E}^{n+1, 2k+2}$, are adopted to monitor the convergence of the iterative process. The detailed form is defined as 
\begin{subequations}
\label{eq:error}
\begin{equation}
    \begin{aligned}
    \mathcal{E}^{n+1, \ast, 2k+2} &= \frac{1}{\Delta t}\left\{\Delta x\sum_{j=1}^N \left[\sum_{l = 1}^M \frac{2}{2l+1}\left(\psi_{l,j}^{n+1, \ast, 2k+2} - \psi_{l,j}^{n+1, \ast, 2k} \right)^2 + 2\left(2\rho_{j}^{n+1,\ast, 2k+2} - 2\rho_{j}^{n+1, \ast, 2k} \right)^2 \right.\right. \\
    &\qquad \left.\left.  + \left(T_{e, j}^{n+1, \ast, 2k+2} - T_{e, j}^{n+1, \ast, 2k}\right)^2 + \left(T_{i, j}^{n+1, \ast, 2k+2} - T_{i, j}^{n+1,\ast,2k}\right)^2\right]\right\}^{1/2};
    \end{aligned}
\label{eq:error_1}
\end{equation}
\begin{equation}
    \mathcal{E}^{n+1, 2k+2} = \frac{1}{\Delta t}\left\{\Delta x\sum_{j=1}^N \left[\left(T_{e, j}^{n+1, 2k+2} - T_{e, j}^{n+1, 2k}\right)^2 + \left(T_{i, j}^{n+1, 2k+2} - T_{i, j}^{n+1, 2k}\right)^2\right]\right\}^{1/2}.\label{eq:error_2}
\end{equation}
\end{subequations}
The parameter $\bar{\epsilon}$ is introduced as a stopping threshold for these errors. It is a problem-specific value and set to $10^{-10}$ for both 1D and 2D problems in our simulation. The  CFL (Courant-Friedrichs-Lewy) condition  is utilized to determine the time step length for this splitting method, which is calculated as 
\begin{equation}
\label{eq:CFL}
    \Delta t =  C  \Delta x / c,
\end{equation}
where $\Delta x$ is the spatial mesh size, $c$ is the speed of light, and $C$ is the CFL number.




\subsection{Tests of the AP property}
\label{sec:AP}

We first investigate the AP property of this numerical scheme. The initial conditions and parameters are 
similar to examples in \cite{Fu2022}, which is an equilibrium as  
\begin{equation}
    T_r=T_e=T_i=(3+\sin(\pi x))/4,\qquad x\in[0,2].
\end{equation}
Periodic condition is applied, with other parameters set as 
\begin{equation}
\label{eq:para_ex0}
    a = 1, \quad c = 1, \quad C_{v_e} = 0.1, \quad C_{v_i} = 0.2, \quad \sigma = 10, \quad K_e = 0.01, \quad K_i = 0.02. 
\end{equation}

\begin{figure}[!hptb]
\centering
\subfloat[$T_r$]{
\includegraphics[width=0.32\linewidth]{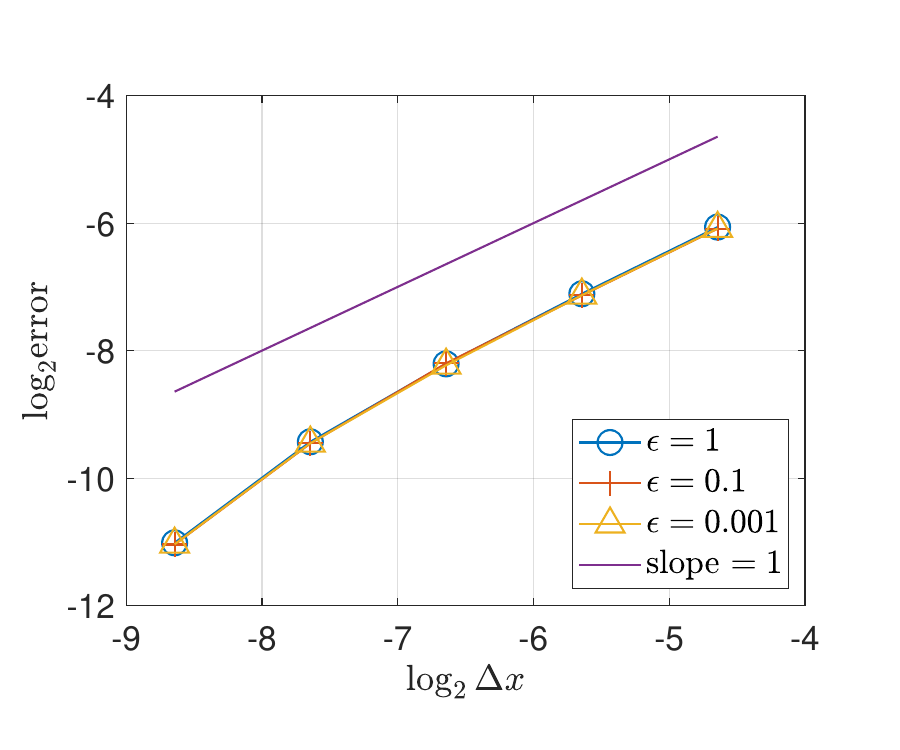}} \hfill 
\subfloat[$T_e$]{
\includegraphics[width=0.32\linewidth]{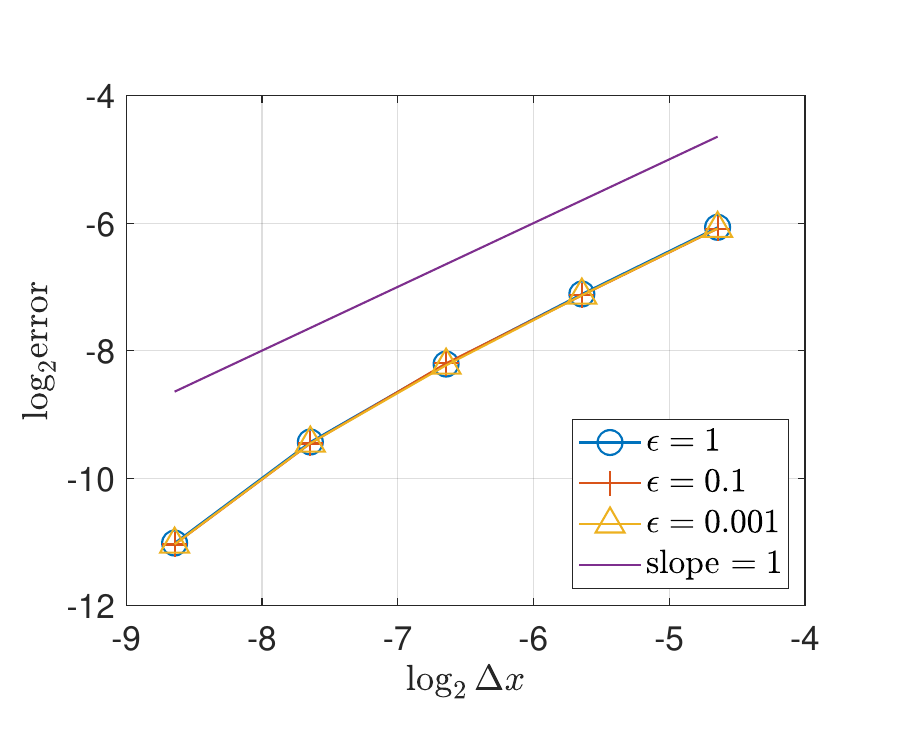} } \hfill
\subfloat[$T_i$]{
\includegraphics[width=0.32\linewidth]{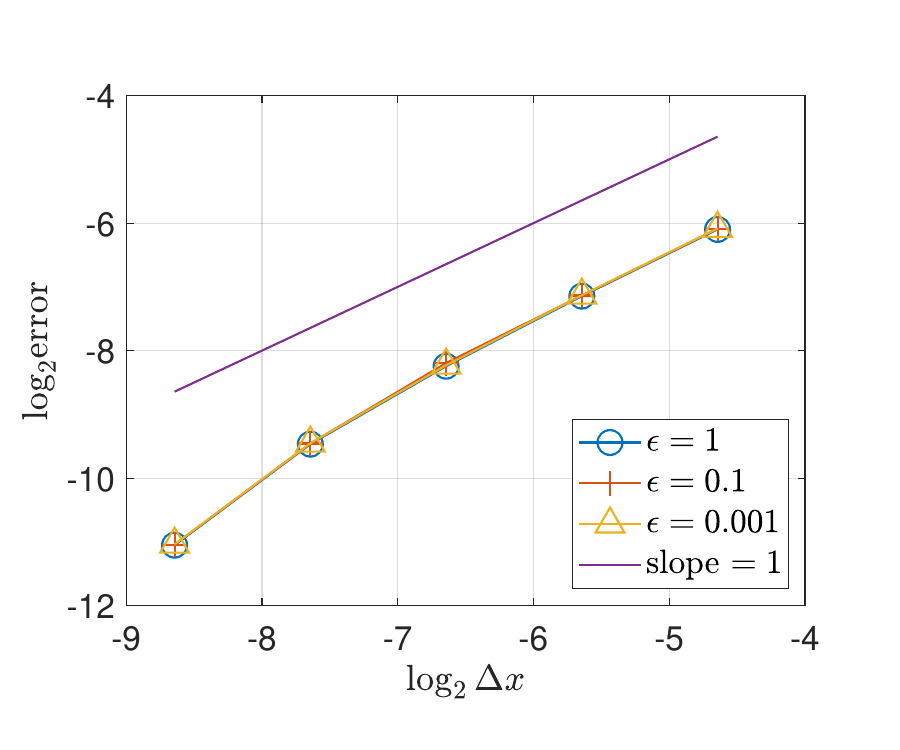}}
\caption{(Tests  of the AP property in Sec. \ref{sec:AP}) The $l_2$ error between the numerical solution with mesh size $N=50,100,200,400$ and $800$ and the reference solution. Here, the reference solution is obtained by the same numerical method with $N=1600$. The parameter $\epsilon$ is set as $\epsilon = 1, 0.1$ and $0.001$ and $\kappa = 1$.}
\label{fig:ex0_eps}
\end{figure}

\begin{figure}[!hptb]
\centering
\subfloat[$T_r$]{
\includegraphics[width=0.32\linewidth]{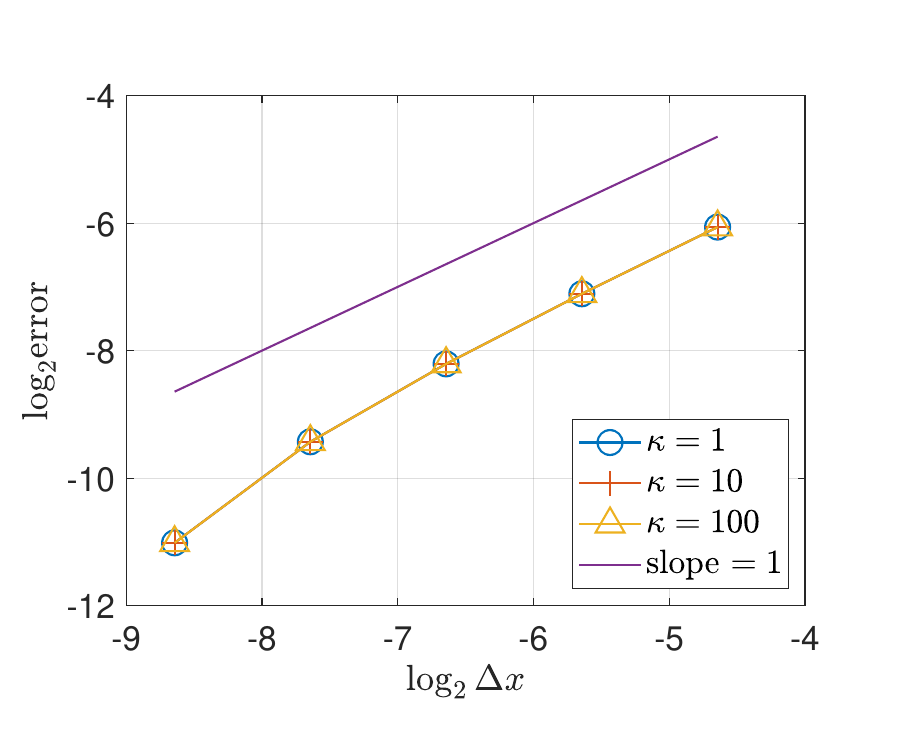}} \hfill
\subfloat[$T_e$]{
\includegraphics[width=0.32\linewidth]{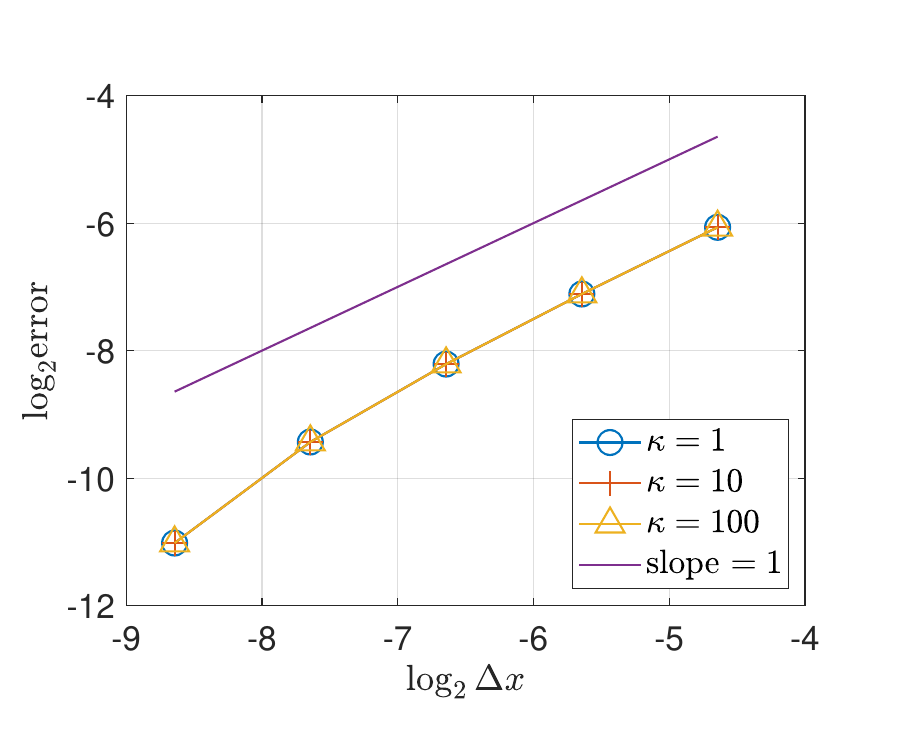}} \hfill
\subfloat[$T_i$]{
\includegraphics[width=0.32\linewidth]{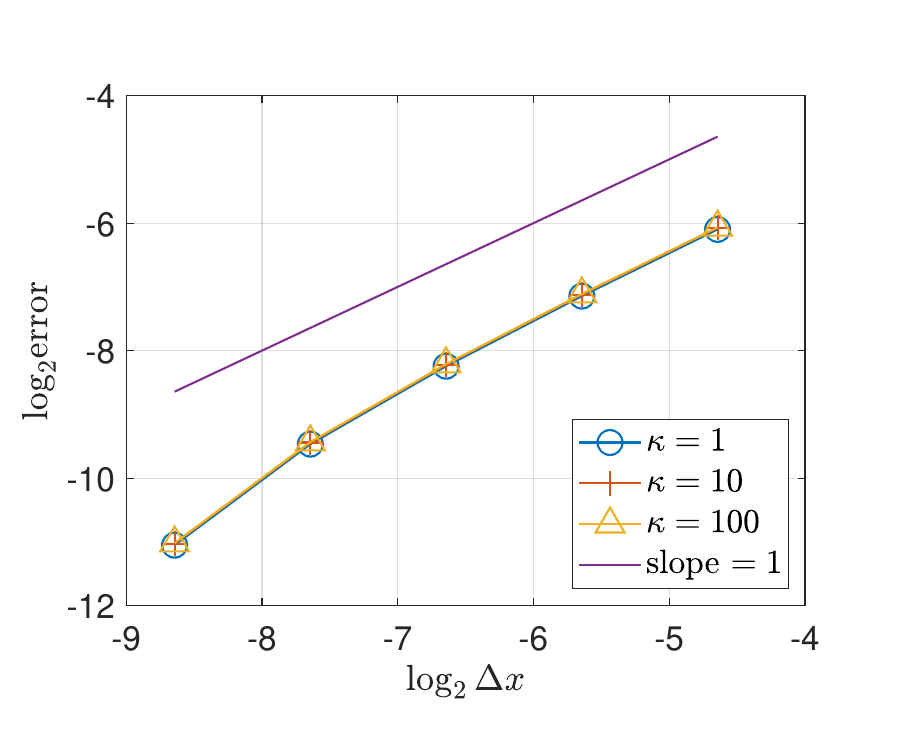}}
\caption{(Tests of the AP property in Sec. \ref{sec:AP}) The $l_2$ error between the numerical solution with grid size $N=50,100,200,400$ and $800$ and the reference solution. Here, the reference solution is obtained by the same numerical method with $N = 1600$. The parameter $\kappa$ is set as $\kappa = 1, 10$ and $100$ and $\epsilon = 1$. }    
\label{fig:ex0_kappa}
\end{figure}

The AP property in the radiation diffusion limit with the parameter $\kappa$ fixed at $\kappa = 1$ is first studied. The expansion order in the $P_N$ method is set to $M = 7$. During the simulation, the mesh size $N$ is varied, with values set as $50$, $100$, $200$, $400$, and $800$. The reference solution is obtained using the same method with a mesh size of $N = 1600$. The CFL number is set to $C = 0.1$, and the time step length is determined as \eqref{eq:CFL}. The $l_2$ error between the numerical solution and the reference solution at $t = 0.1$ is depicted in Fig. \ref{fig:ex0_eps}, where $\epsilon$ is set to $1, 0.1$, and $0.001$. The results in Fig. \ref{fig:ex0_eps} demonstrate that, for all three temperatures, the $l_2$ error converges at the first order for all values of $\epsilon$, confirming the AP property of this splitting scheme concerning $\epsilon$. In Fig \ref{fig:ex0_kappa}, the AP property related to $\kappa$ is investigated, with $\epsilon$ fixed at $\epsilon = 1$, and the other parameters identical as \eqref{eq:para_ex0}. Fig. \ref{fig:ex0_kappa} displays the $l_2$ error between the numerical solution and the reference solution at $t = 0.1$, with $\kappa$ taking values of $\kappa = 1, 10$, and $100$. The results in Fig. \ref{fig:ex0_kappa} demonstrate that, for all three temperatures, the $l_2$ error converges at the first order for all values of $\kappa$, providing evidence for the AP property of this splitting scheme concerning $\kappa$.

To investigate the impact of the CFL number on the efficiency of the splitting method \eqref{eq:PN_full}, we examine the iteration numbers for both steps in the numerical scheme under various CFL numbers. Initially, with $\kappa = 1$, we test three different values of $\epsilon$ as $\epsilon = 1, 0.1$, and $0.001$. Tab. \ref{fig:ex0_itereps} displays the iteration numbers for both steps of the numerical scheme \eqref{eq:PN_full} when the CFL number $C = 0.01, 0.05$, and $0.1$, with a fixed mesh size of $N = 200$. The results reveal that, for the first part of the numerical scheme, iteration numbers are consistently small across different CFL numbers and all three values of $\epsilon$. In the second part of the scheme, the iteration numbers are consistently small with the CFL number for all tested $\epsilon$ values too.

Subsequently, with $\epsilon = 1$ held constant, we explore the effect of three different $\kappa$ values as $\kappa = 1, 10$, and $100$. Tab. \ref{fig:ex0_iterkappa} illustrates the iteration numbers for both steps of the numerical scheme \eqref{eq:PN_full} at the CFL number $C = 0.01, 0.05$, and $0.1$, with a fixed mesh size of $N = 200$. Similar to the $\epsilon$ scenario, the iteration numbers for the first part of the numerical scheme consistently stay small across various CFL numbers and all three $\kappa$ values. In the second part of the scheme, iteration numbers also consistently remain small.

In summary, across all tested regimes of $\epsilon$ and $\kappa$, the iteration numbers for both steps in the numerical scheme \eqref{eq:PN_full} are consistently small. Importantly, the CFL number shows minimal influence on the iteration numbers, indicating the high efficiency of this asymptotic-preserving splitting numerical method.

\begin{table}[!hptb]
\centering
\def\arraystretch{1.5}
{\footnotesize
\begin{tabular}{c|ccc|ccc}
\multirow{2}{*}{iteration   number} & \multicolumn{3}{c|}{first iteration}                                    & \multicolumn{3}{c}{second iteration}                                   \\ \cline{2-7} 
& \multicolumn{1}{c|}{$C = 0.01$} & \multicolumn{1}{c|}{$C = 0.05$} & $C = 0.1$ & \multicolumn{1}{c|}{ $C = 0.01$} & \multicolumn{1}{c|}{$C = 0.05$} & $C = 0.1$ \\ \hline
$\epsilon=1$                        & \multicolumn{1}{c|}{6}        & \multicolumn{1}{c|}{6}        & 6       & \multicolumn{1}{c|}{6}        & \multicolumn{1}{c|}{8}        & 8       \\ \hline
$\epsilon=0.1$                      & \multicolumn{1}{c|}{6}        & \multicolumn{1}{c|}{8}        & 8       & \multicolumn{1}{c|}{8}        & \multicolumn{1}{c|}{8}        & 8       \\ \hline
$\epsilon=0.001$                    & \multicolumn{1}{c|}{8}        & \multicolumn{1}{c|}{8}        & 8       & \multicolumn{1}{c|}{8}        & \multicolumn{1}{c|}{8}        & 8       \\ 
\end{tabular}
}
\caption{(Tests of the AP property in Sec. \ref{sec:AP}) Iteration numbers for the both steps in the numerical scheme \eqref{eq:PN_full} with CFL number $C = 0.01, 0.05$ and $0.1$ for $\epsilon = 1, 0.1$, and $0.001$ and $\kappa = 1$. }
\label{fig:ex0_itereps}
\end{table}

\begin{table}[!hptb]
\centering
\def\arraystretch{1.5}
{\footnotesize
\begin{tabular}{c|ccc|ccc}
\multirow{2}{*}{iteration   number} & \multicolumn{3}{c|}{first iteration}                                    & \multicolumn{3}{c}{second iteration}                                   \\ \cline{2-7} 
& \multicolumn{1}{c|}{$C = 0.01$} & \multicolumn{1}{c|}{$C = 0.05$} & $C = 0.1$ & \multicolumn{1}{c|}{ $C = 0.01$} & \multicolumn{1}{c|}{$C = 0.05$} & $C = 0.1$\\ \hline
$\kappa=1$                          & \multicolumn{1}{c|}{6}        & \multicolumn{1}{c|}{6}        & 6       & \multicolumn{1}{c|}{6}        & \multicolumn{1}{c|}{8}        & 8       \\ \hline
$\kappa=10$                         & \multicolumn{1}{c|}{6}        & \multicolumn{1}{c|}{6}        & 6       & \multicolumn{1}{c|}{6}        & \multicolumn{1}{c|}{8}        & 8       \\ \hline
$\kappa=100$                        & \multicolumn{1}{c|}{6}        & \multicolumn{1}{c|}{6}        & 6       & \multicolumn{1}{c|}{6}        & \multicolumn{1}{c|}{8}        & 8       \\ 
\end{tabular}
}
\caption{(Tests of the AP property in Sec. \ref{sec:AP}) Iteration numbers for the both steps in the numerical scheme \eqref{eq:PN_full} with CFL number $C = 0.01, 0.05$ and $0.1$ for $\kappa = 1, 10$, and $100$ and $\epsilon = 1$.  }  
\label{fig:ex0_iterkappa}
\end{table}

\subsection{Homogeneous model problem}
\label{sec:ex1}

In this section, we examine two homogeneous model problems, each involving a source affecting either the ions or the photons, respectively. In this case, the system becomes an ODE system with three equations with strong stiffness. We can verify that the splitting scheme can provide an accurate solution even if the ODE is stiff. The parameters for these tests are consistent with those in references like \cite{3T2022, enaux2020numerical, EVANS20071695}. Detailed parameters for both tests are provided in Tab. \ref{tab:ex1}. The source term is defined as
\begin{equation}
\label{eq:ex1_para}
    Q = \frac{\bar{\rho} \mathcal{C}}{\sqrt{2 \pi} t_w} \exp\left(\frac{(t - t_c)^2}{2t_w^2}  \right) [{\rm GJ cm^{-3} ns^{-1}}],
\end{equation}
and the initial conditions for both problems are as follows.
\begin{equation}
\label{eq:ex1_ini}
    T_i = T_e = T_r  =  2.52487 \times 10^{-5} {\rm keV}.
\end{equation}

\begin{table}[!hptb]
\centering
\def\arraystrech{1.5}
{\footnotesize
\begin{tabular}{c||c|c}
&  Problem I & Problem II \\
\hline 
$Q_i$ & \eqref{eq:ex1_para} & $0$\\
$Q_r$ & $0$ & \eqref{eq:ex1_para} \\
$C_{v_e}$ & $  0.1\bar{\rho}$ &  $0.1T_e \bar{\rho}$  \\ 
$C_{v_i}$ & $0.05 \bar{\rho}$ & $0.05 \bar{\rho}$ \\
$\kappa$ & $C_{v_e}/ (c\tau)$  & $0.01379 T_e^{-0.5}$  \\[1mm] \hline \hline
$c~({\rm cm ns}^{-1})$ &  $29.979$ &  $29.979$ \\  
$a ~({\rm GJ cm}^{-3} {\rm keV}^4)$ & $0.01372$ & $0.01372$ \\
$\tau ~({\rm ns})$ & $0.1$ & $0.1$  \\
$\epsilon$  & $1$ & $1$ \\
$\mathcal{C}$& 25.06628 & 25.06628 \\
$t_w~ ({\rm ns})$ & 1.0 & 1.0\\
$t_c~ ({\rm ns})$ & 1.0 & 1.0 \\
$\bar{\rho}~ ({\rm g cm}^{-3})$  & 3.0 & 3.0 \\
$\sigma~({\rm cm}^{-1})$ & $ \sigma_0/T_e^2$ &$ \sigma_0/T_e^2$   \\
$\sigma_0$ & 0.5 & 0.5\\
$M$ & 7 & 7\\
$\Delta t ({\rm ns})$ & $0.0025$  & $0.0025$
\end{tabular}
}
\caption{(Homogeneous model problem in Sec. \ref{sec:ex1}) Parameters for the homogeneous problem with source terms.}
\label{tab:ex1}
\end{table}

\begin{figure}[!hptb]
\centering
\subfloat[Problem I]{
\includegraphics[width=0.45\linewidth]{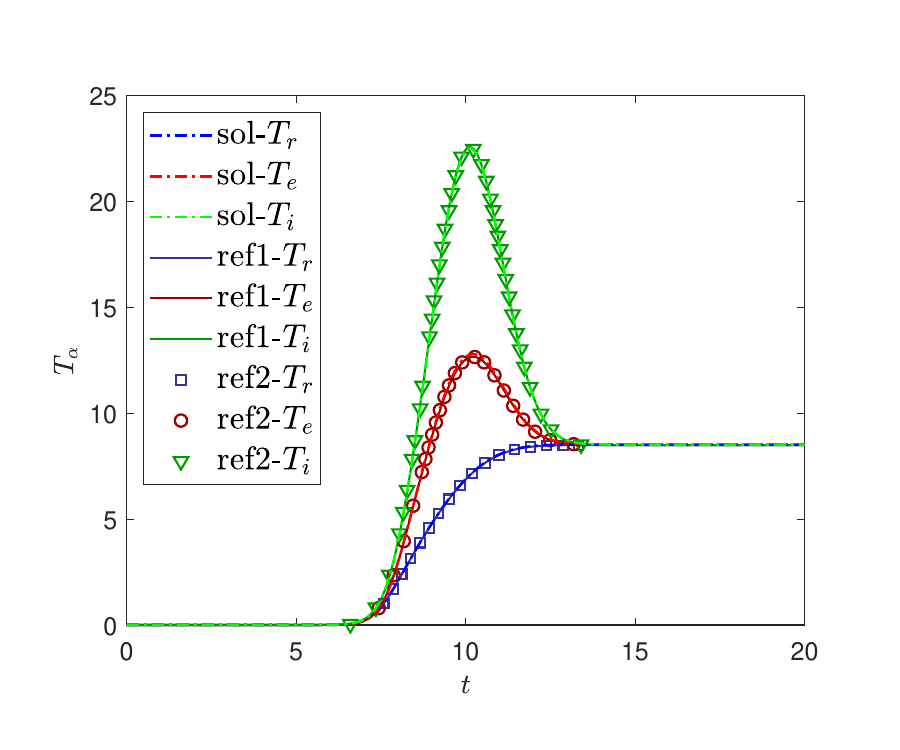}
\label{fig:ex1_sol}
}
\subfloat[Problem II]{
\includegraphics[width=0.45\linewidth]{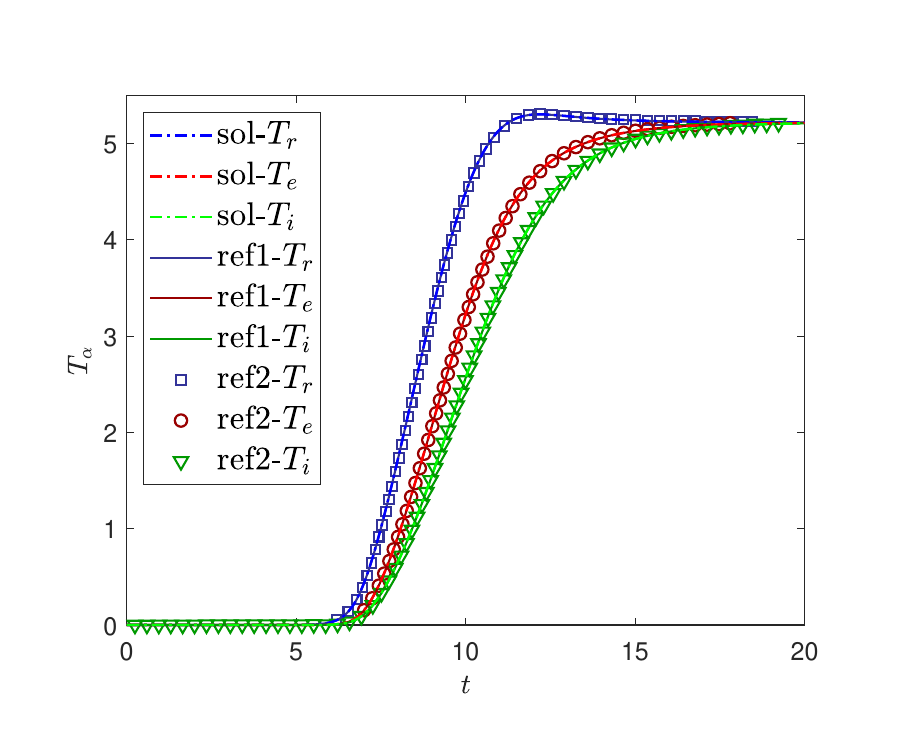}
\label{fig:ex1_test2sol}
}
\caption{(Homogeneous model problem in Sec. \ref{sec:ex1}) The time evolution of the numerical solution, the reference solution by the same method with a time step length of $\Delta t=0.01\times 2^{-10}$ns, and the reference solution from \cite{EVANS20071695}. The dot-dashed lines are the numerical solution by the AP $P_N$ method. The solid lines are the reference solution by the same method with a time step length of $\Delta t=0.01\times 2^{-10}$ns, and the symbols are the reference solutions from \cite{EVANS20071695}.  (a) Solutions of Problem I, (b) Solutions of Problem II.}
\label{fig:ex1}
\end{figure}

\paragraph{Problem I: Ion-source homogeneous problem}
In Problem I, the source term of the ion temperature is introduced. The source term is specified in \eqref{eq:ex1_para}, and the initial condition is presented in \eqref{eq:ex1_ini}. Additionally, the opacity depends on the temperature of electron $T_e$ and follows the relation
$\sigma = \frac{1}{2T_e^2}.$

\begin{figure}[!hptb]
    \centering
    \subfloat[First order]{
    \includegraphics[width=0.45\linewidth]{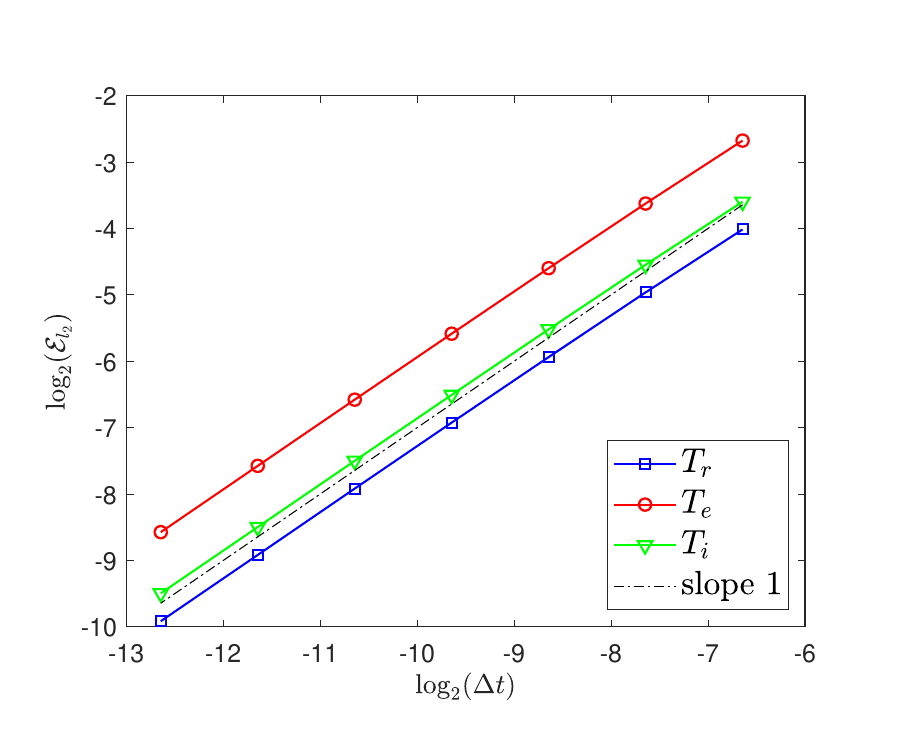}   
      \label{fig:ex1_err}
      } \hfill
    \subfloat[Second order]{
    \includegraphics[width=0.45\linewidth]{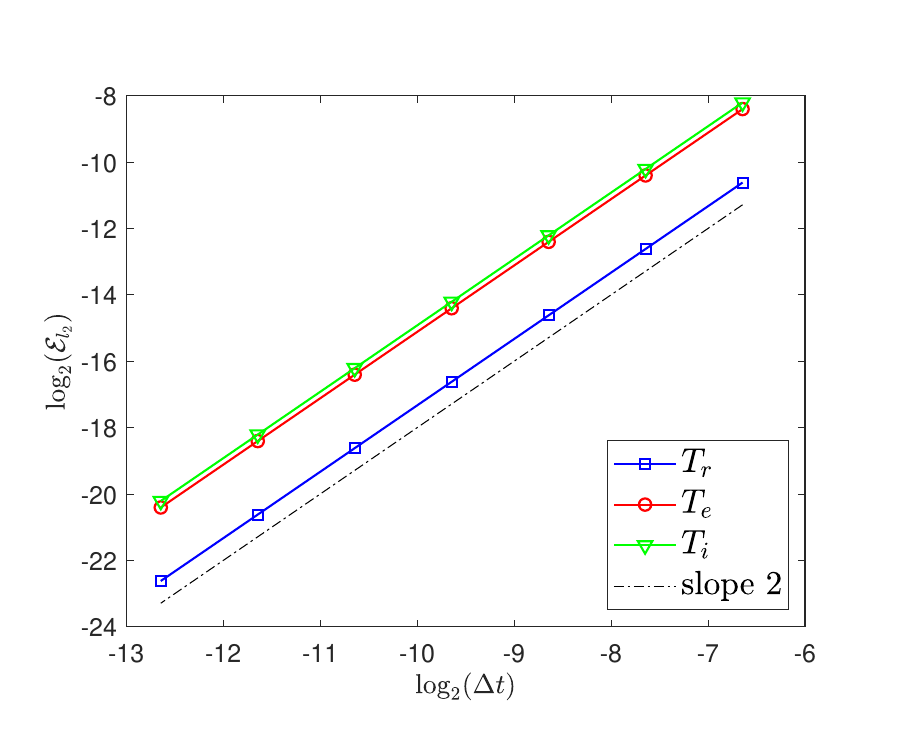}   
    \label{fig:ex1_errorder2}
    }
    \caption{(Homogeneous model problem in Sec. \ref{sec:ex1}) The $l_2$ error of the numerical solution by the AP $P_N$ method and the reference solution by the same method with a time step length of $\Delta t=0.01\times 2^{-10}$ns. (a) The error for the first-order $\Pn$ method, (b) The error for the second-order $\Pn$ method.     
    }
\end{figure}

During the simulation, the expansion number for the $P_N$ method is set to $M = 7$, and the time step length is defined as $\Delta t = 0.0025$ns. The simulation is run until the final time $t = 20$ns. The time evolution of $T_r$, $T_e$, and $T_i$ is illustrated in Fig. \ref{fig:ex1_sol}, along with the reference solution obtained using the same method with a time step length of $\Delta t=0.01\times 2^{-10}$ ns and reference solutions from \cite{EVANS20071695}. The simulation results show close agreement with the reference solutions. Initially, all three temperatures are identical, but they gradually diverge over time. Due to the source term affecting the ions, the ion temperature exhibits the highest value. The $l_2$ error between the numerical solution and the reference solution is presented in Fig. \ref{fig:ex1_err}, demonstrating first-order convergence in time. When the implicit midpoint scheme is employed for temporal discretization, the convergence rate improves to the second order as shown in Fig. \ref{fig:ex1_errorder2}. Here, the implicit midpoint method with a time step length of $\Delta t=0.01\times 2^{-10}$ ns is adopted as the reference solution.

\paragraph{Problem II: Photons-source homogeneous problem}
In Problem II, the source term in \eqref{eq:ex1_para} is introduced for the photons, and the identical initial condition \eqref{eq:ex1_ini} is employed. The specific parameters are detailed in Tab. \ref{tab:ex1}, where the specific heat of electrons $C_{v_e}$ depends on the electron temperature and is defined as
\begin{equation}
\label{eq:ex1_cv}
    C_{v_e} = 0.1 \bar{\rho} T_e, \qquad \bar{\rho} = 3.0,
\end{equation}
and the coefficient $\kappa$ depends on $T_e$ as well 
\begin{equation}
\label{eq:ex1_kappa} 
    \kappa = 0.01379 T_e^{-0.5}.
\end{equation}
In the simulation, the expansion number for the $P_N$ method is set to $M = 7$, and the time step length is configured as $\Delta t = 0.0025$ns. The time evolution of the radiative temperature $T_r$, electron temperature $T_e$, and ion temperature $T_i$ is illustrated in Fig. \ref{fig:ex1_test2sol}, alongside the reference solutions obtained using the same method with a time step length of $\Delta t=0.01\times 2^{-10}$ns and the reference solution from \cite{3T2022}, respectively. The simulation is conducted until the final time $t = 20$ns. The results indicate that the numerical solution matches the reference solution well. In this problem, all three temperatures are nearly identical initially, with the radiative temperature $T_r$ gradually becoming the highest over time due to the influence of the source term.

\subsection{1D Marshak wave problem without conduction terms}
\label{sec:marshak_ex2}

In this section, we investigate the Marshak problem without conduction terms. This problem is a classical benchmark extensively studied in \cite{enaux2020numerical, Fu2022, sun2015asymptotic1}. For this specific example, the computational domain is $[0, 0.5]$ and the equilibrium initial condition is chosen as
\begin{equation}
\label{eq:ex2_ini}
    T_e=T_i=10^{-6}, \quad \psi=\frac{1}{2}a c T_e^{4}.
\end{equation}
The inflow boundary conditions are defined as  
\begin{equation}
\label{eq:ex2_boun}
    {\rm left~boundary:}~T_e=T_i=1,\quad \psi=\frac{1}{2}ac T_e^4, \qquad {\rm right~boundary:}~T_e=T_i=10^{-6}, \quad \psi=\frac{1}{2}acT_e^4,
\end{equation}
where the parameters are set as follows.
\begin{equation}
\label{eq:ex2_param}
    \epsilon=1,\qquad a=0.01372, \qquad c=299.79, \qquad C_{v_e}=0.03, \qquad C_{v_i}=0.27, \qquad \sigma=300T_e^{-3}.
\end{equation}
Detailed application of the inflow boundary condition within the framework of the $P_N$ method is proposed in \cite{Fu2022}. 

\begin{figure}[!hptb]
\centering
\subfloat[$T_r, \kappa = 1$]{
\includegraphics[width=0.3\linewidth]{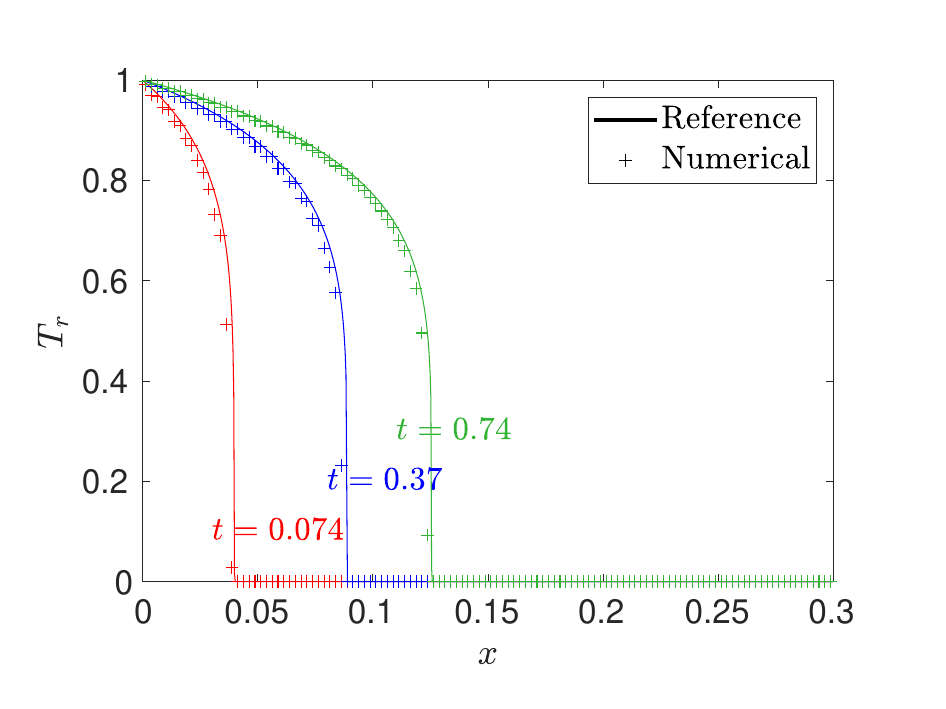}
}\hfill
\subfloat[$T_e, \kappa = 1$]{
\includegraphics[width=0.3\linewidth]{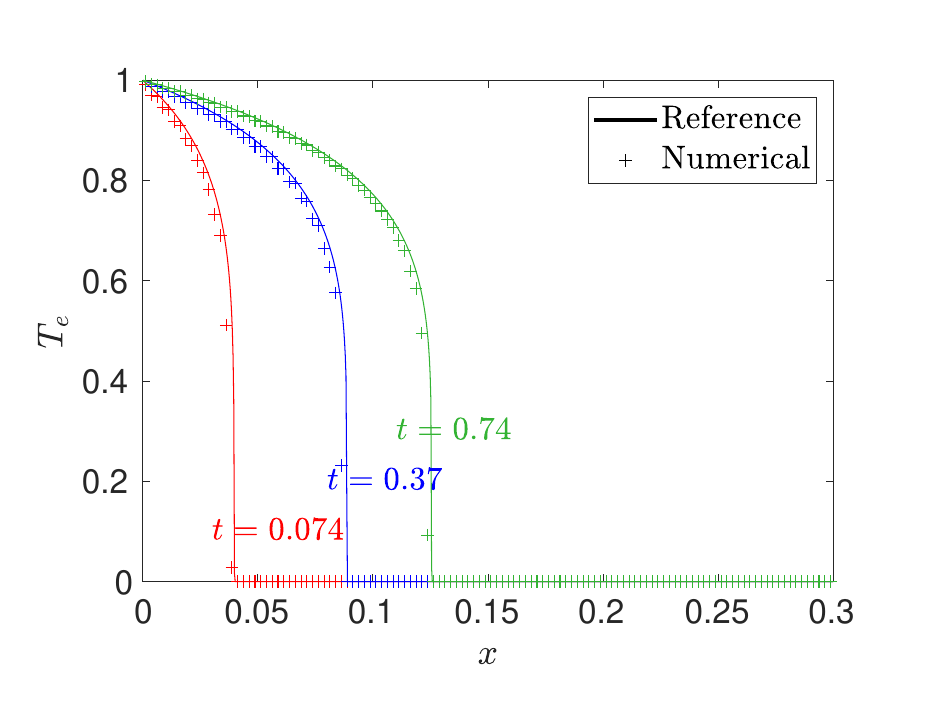}
}\hfill
\subfloat[$T_i, \kappa = 1$]{
\includegraphics[width=0.3\linewidth]{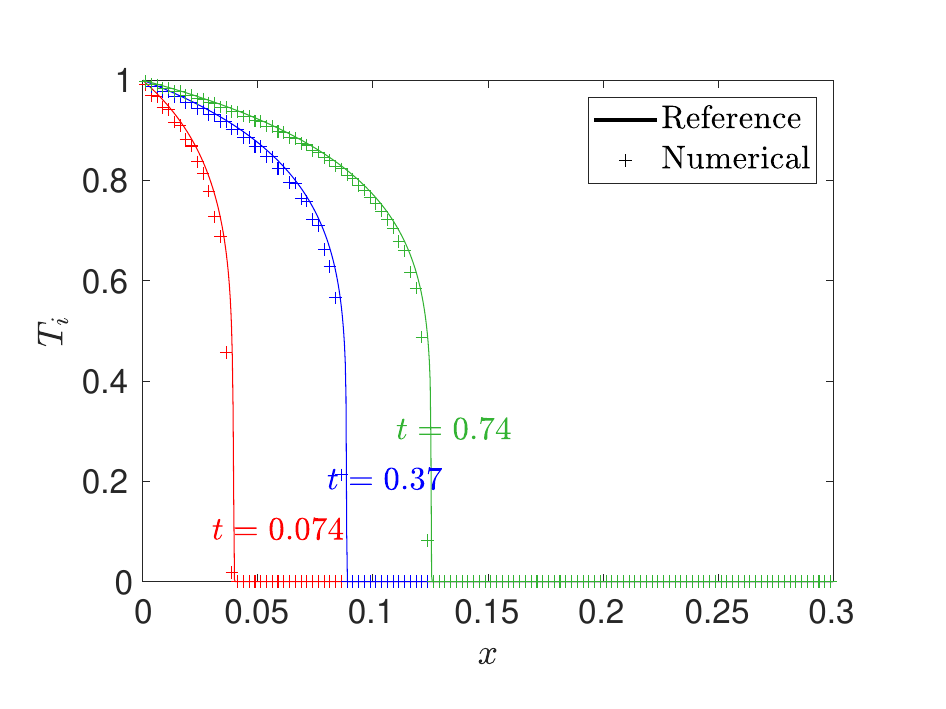}
}
\caption{(1D Marshak wave problem without conduction terms in Sec. \ref{sec:marshak_ex2}) The numerical solution of $T_r, T_e$ and $T_i$ for $\kappa = 1$ at different times. The symbols are the numerical solution by the AP $P_N$ method, and the solid lines are the reference solution obtained from the method in \cite{enaux2020numerical}.}
\label{fig:ex2_kappa1e0}
\end{figure}

In the absence of the conduction terms, both $D_{e}$ and $D_{i}$ are set to zero. A grid size of $N=200$ is employed, along with linear reconstruction. We use $M=7$ in the $P_N$ method. The numerical solution for the radiative temperature $T_r$, as well as the temperature of the electron and ion $T_e, T_i$ at time $t = 0.074, 0.37$, and $0.74$ for $\kappa = 1$, is depicted in Fig. \ref{fig:ex2_kappa1e0}. Additionally, the reference solution obtained using the method from \cite{enaux2020numerical} is included. The results illustrate that all three temperatures propagate forward and closely align with the reference solution. 

\begin{figure}[!hptb]
\centering
\subfloat[$T_r, \kappa = 0.001$]{
\includegraphics[width=0.3\linewidth]{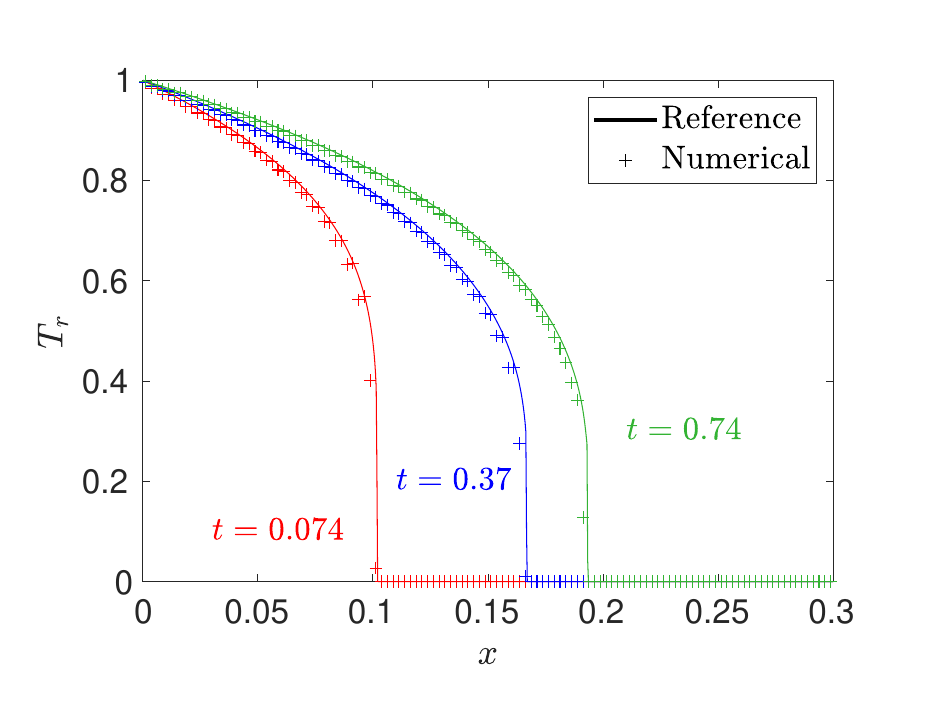}
}\hfill
\subfloat[$T_e, \kappa = 0.001$]{
\includegraphics[width=0.3\linewidth]{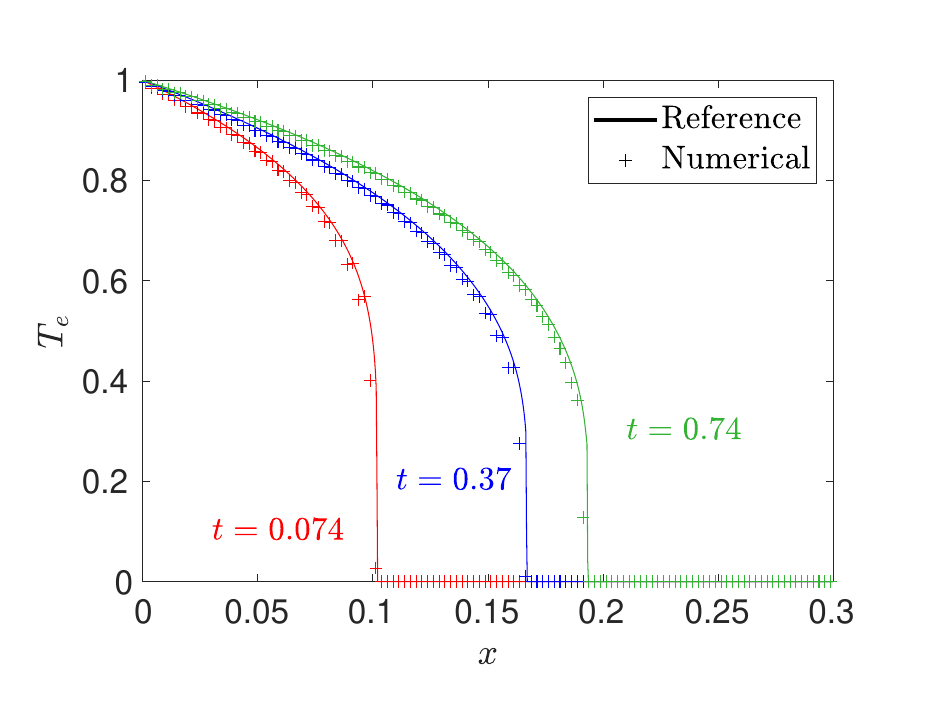}
}\hfill
\subfloat[$T_i, \kappa = 0.001$]{
\includegraphics[width=0.3\linewidth]{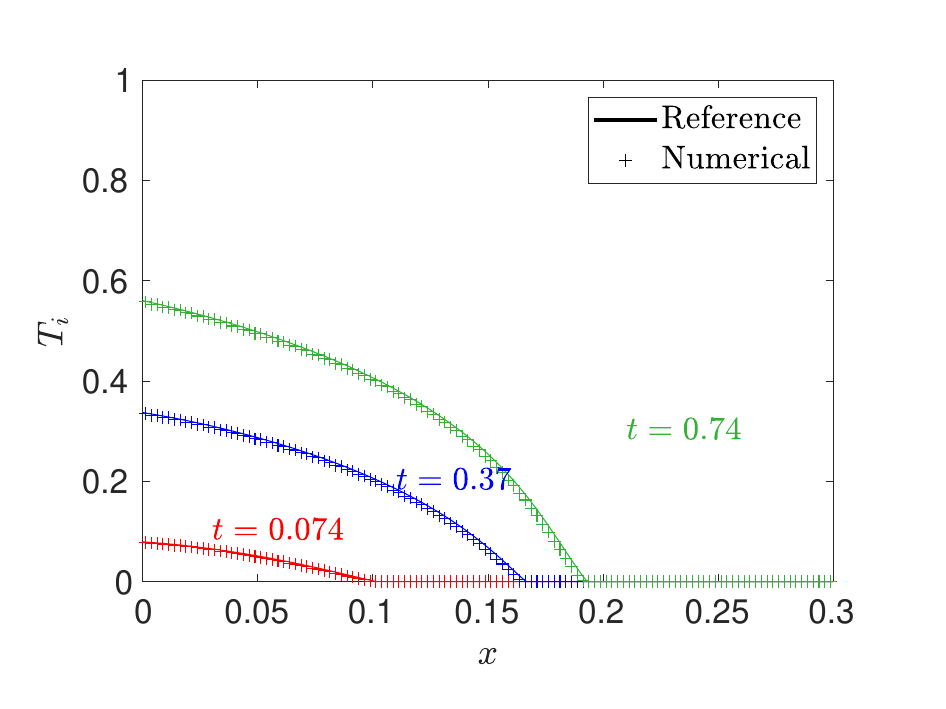}
\label{fig:ex2_kappa_0p001_Ti}
} \\
\subfloat[$T_r, \kappa = 100$]{
\includegraphics[width=0.3\linewidth]{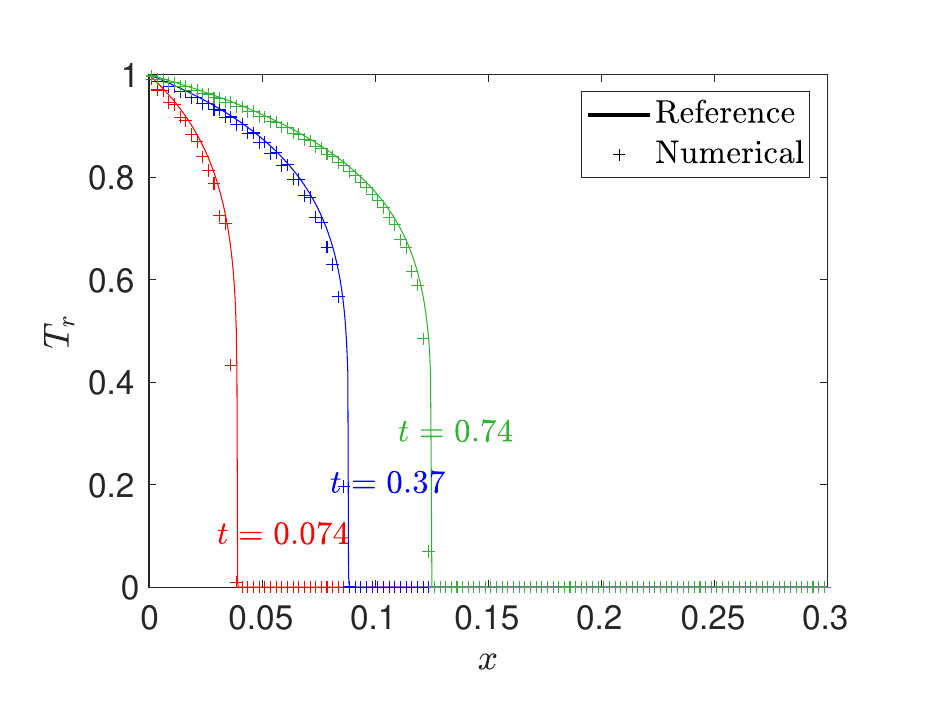}
}\hfill
\subfloat[$T_e, \kappa = 100$]{
\includegraphics[width=0.3\linewidth]{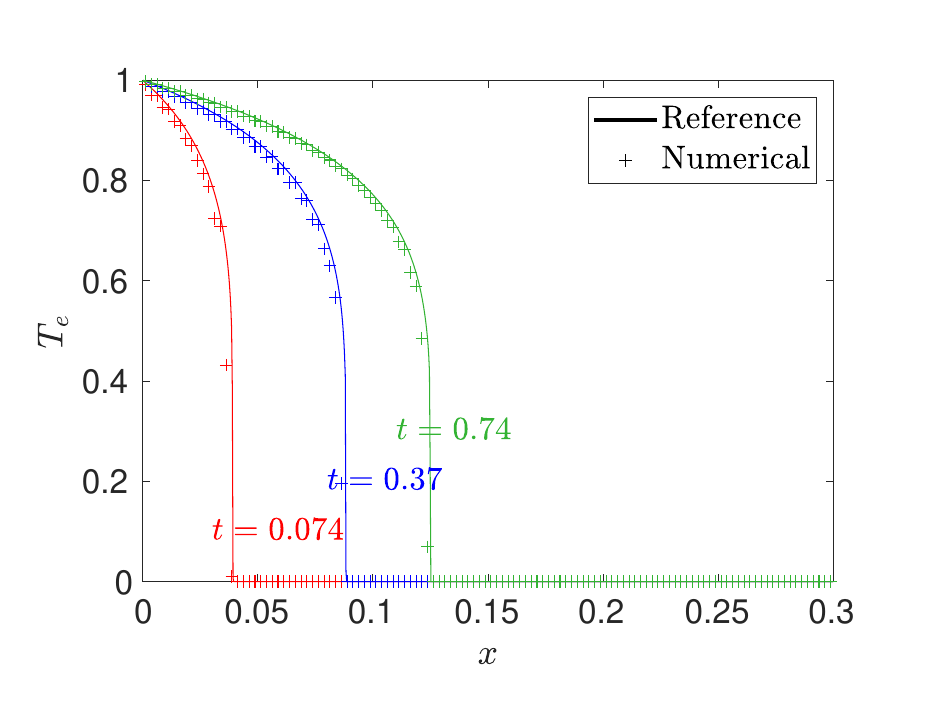}
}\hfill
\subfloat[$T_i, \kappa = 100$]{
\includegraphics[width=0.3\linewidth]{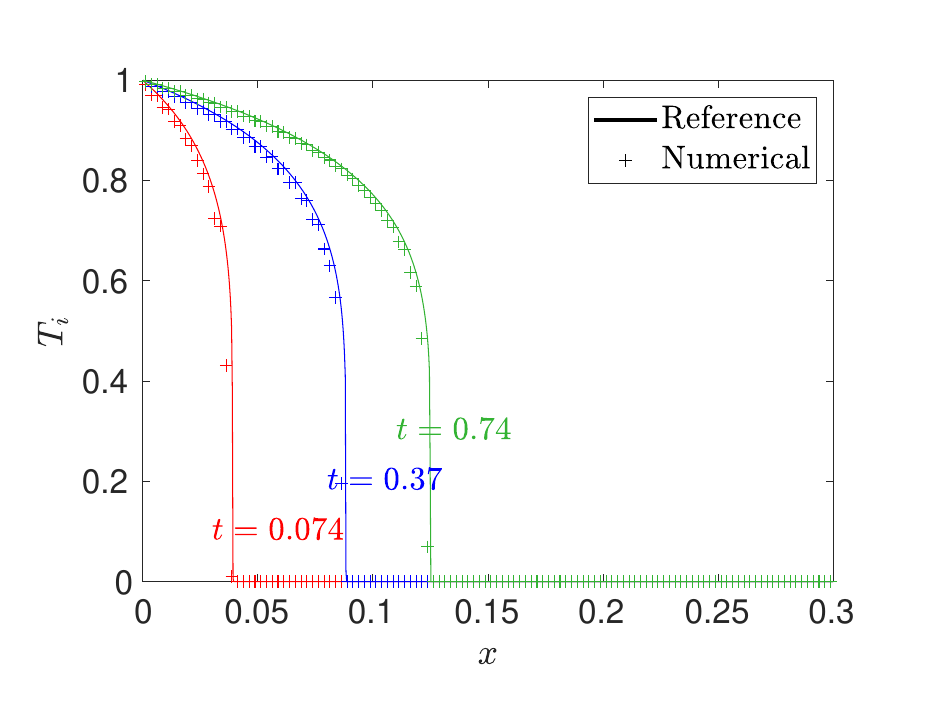}
\label{fig:ex2_kappa_100_Ti}
}
\caption{(1D Marshak wave problem without conduction terms in Sec. \ref{sec:marshak_ex2}) The numerical solution of $T_r, T_e$ and $T_i$ for $\kappa = 0.001$ and $100$ at different times. The top column is for $\kappa =0.001$ and the bottom column is for $\kappa = 100$. The symbols are the numerical solution by the AP $P_N$ method, and the solid lines are the reference solution obtained from the method in \cite{enaux2020numerical}.}
\label{fig:ex2_kappa}
\end{figure}

\begin{figure}
\centering
\subfloat[$\kappa = 0.001$]{
\includegraphics[width=0.3\linewidth]{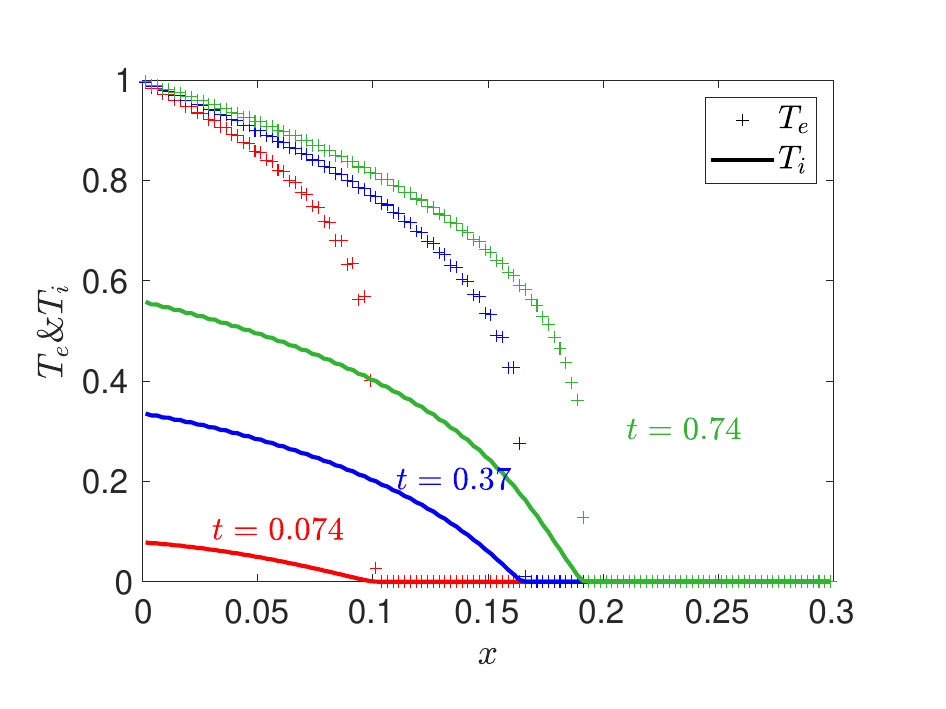}
\label{fig:ex2_kappa_compare_0p001}
}\hfill
\subfloat[$\kappa = 1$]{
\includegraphics[width=0.3\linewidth]{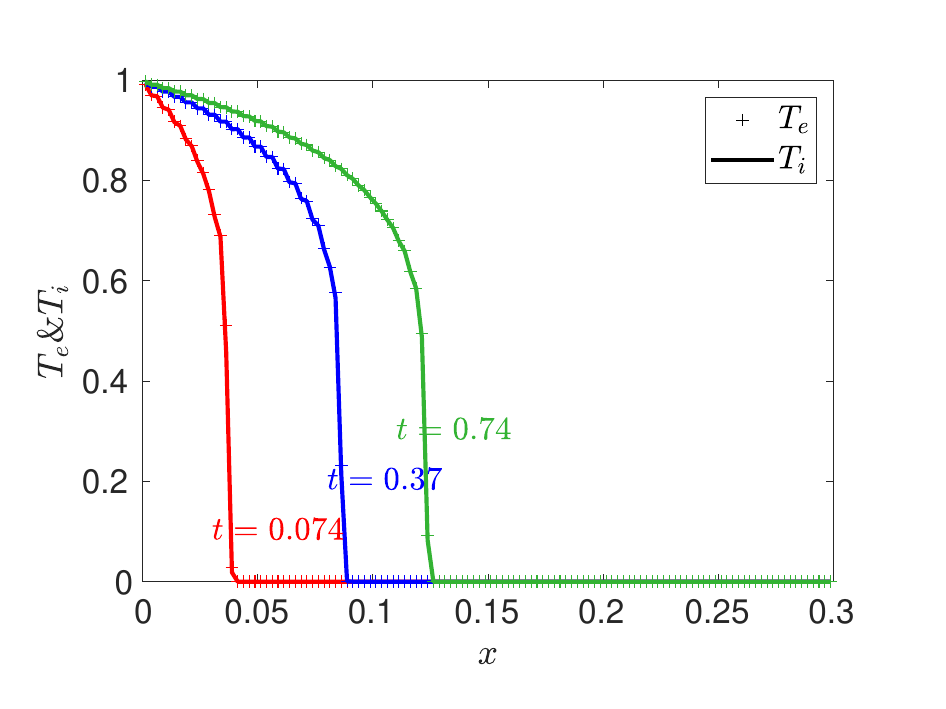}
\label{fig:ex2_kappa_compare_1}
}\hfill
\subfloat[$\kappa = 100$]{
\includegraphics[width=0.3\linewidth]{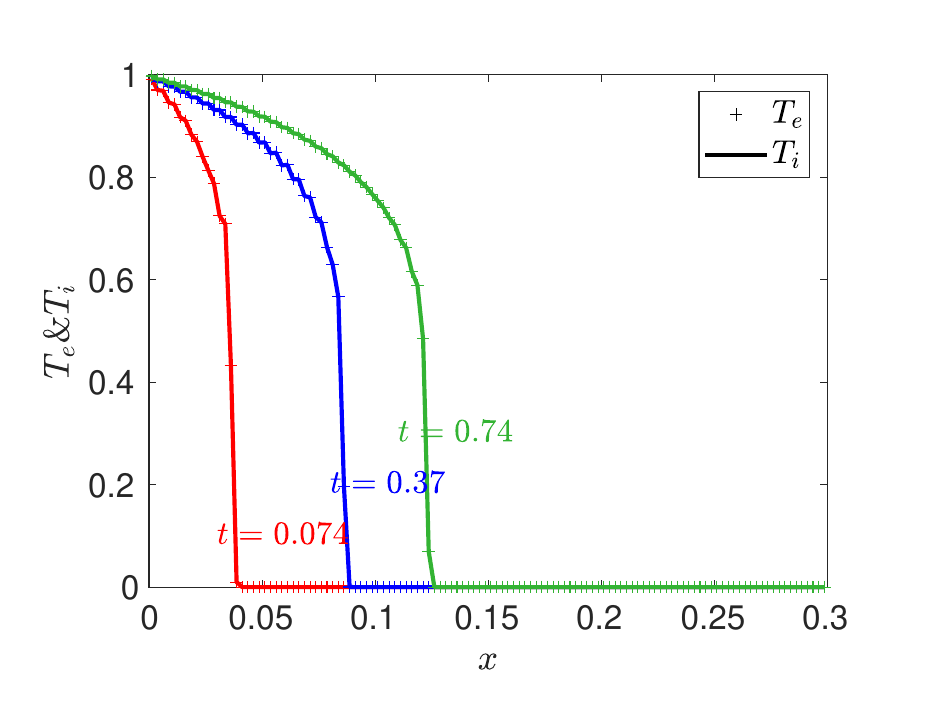}
\label{fig:ex2_kappa_compare_100}
}
\caption{(1D Marshak wave problem without conduction terms in Sec. \ref{sec:marshak_ex2}) The comparison of $T_e$ and $T_i$ for different $\kappa$ at different times. The symbols are the numerical solution of $T_e$ and the solid lines are the numerical solution of $T_i$.  (a) $\kappa = 0.001$, (b) $\kappa = 1$, (c) $\kappa = 100$. }
\label{fig:ex2_kappa_compare}
\end{figure}

To demonstrate the AP property of the scheme proposed in Sec. \ref{sec:alg}, the behavior of the numerical solution is examined for two additional values of $\kappa$: $\kappa = 0.001$ and $100$, as depicted in Fig. \ref{fig:ex2_kappa}. For both small $\kappa$ ($\kappa = 0.001$) and large $\kappa$ ($\kappa = 100$), Fig. \ref{fig:ex2_kappa} reveals that the numerical solution closely aligns with the reference solution. Importantly, it is clearly illustrated that the behavior of the ion temperature $T_i$ differs from that of $T_e$ at $\kappa = 0.001$, while the two behave similarly for $\kappa = 100$. For a more precise comparison of the behavior of $T_e$ and $T_i$ for different $\kappa$, the numerical solutions of $T_e$ and $T_i$ for various $\kappa$ at $t = 0.74$ are presented in Fig. \ref{fig:ex2_kappa_compare}. It is observed that when $\kappa = 0.001$ in Fig. \ref{fig:ex2_kappa_compare_0p001}, the two temperatures exhibit significant differences. However, as $\kappa$ increases to $\kappa = 1$, the numerical solutions of $T_e$ and $T_i$ become nearly identical. Furthermore, the numerical solutions for these two temperatures completely overlap when $\kappa$ increases to $100$. Additionally, the same time step length, which is independent of $\kappa$, is employed for all three cases, indicating the AP property of this new method concerning $\kappa$.

\subsection{Marshak wave problem with conduction terms}
\label{sec:marshak_ex3}

In this section, the Marshak wave problem with conduction terms is examined. The same initial and boundary conditions as \eqref{eq:ex2_ini} and \eqref{eq:ex2_boun} are utilized. The parameters are set as 
\begin{equation}
    \epsilon=1,\quad a=0.01372,\quad c=29.979,\quad C_{v_e}=0.3,\quad C_{v_i}=0.27,\quad \sigma=300T_e^{-3}.
\end{equation}
Additionally, the conduction coefficients are set to $K_{e} = 1.0$ and $K_{i} = 0.1$.

\begin{figure}[!hptb]
\centering
\subfloat[$T_r, \kappa = 1$]{
\includegraphics[width=0.3\linewidth]{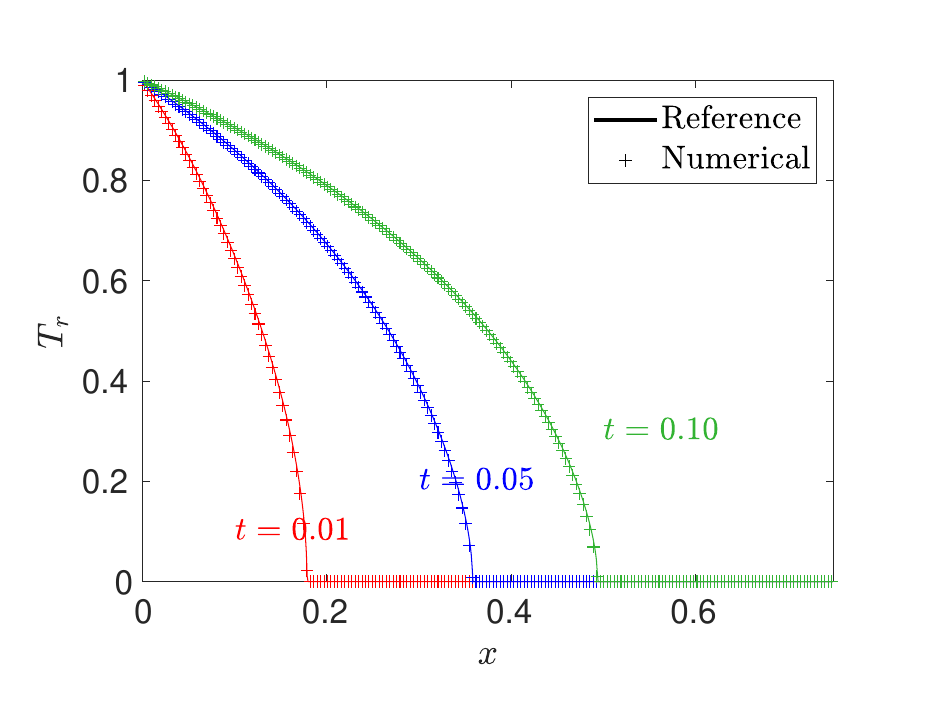}
}\hfill
\subfloat[$T_e, \kappa = 1$]{
\includegraphics[width=0.3\linewidth]{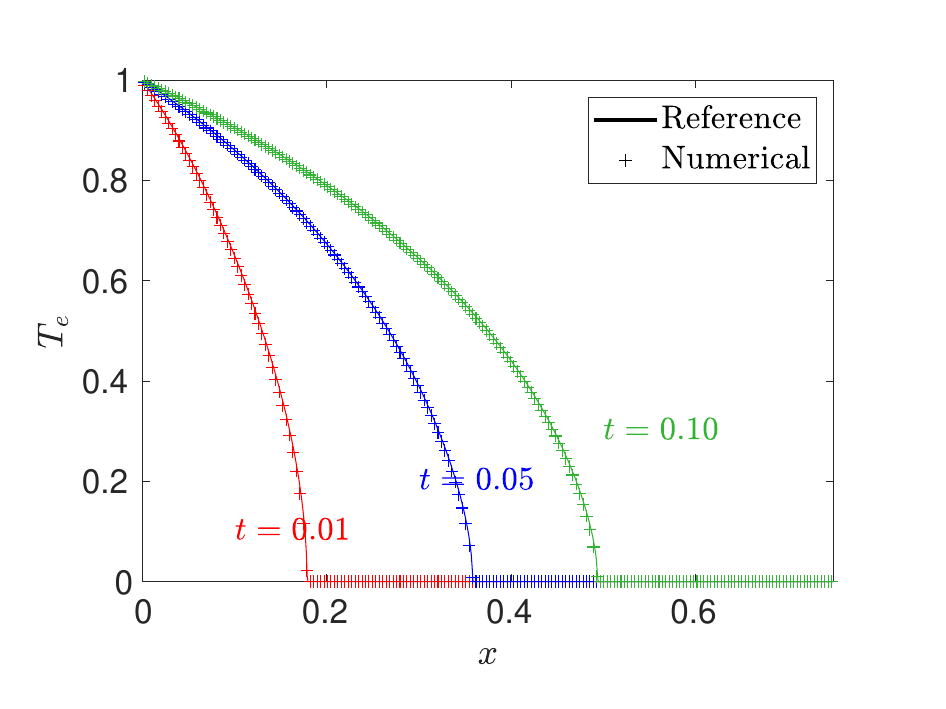}
}\hfill
\subfloat[$T_i, \kappa = 1$]{
\includegraphics[width=0.3\linewidth]{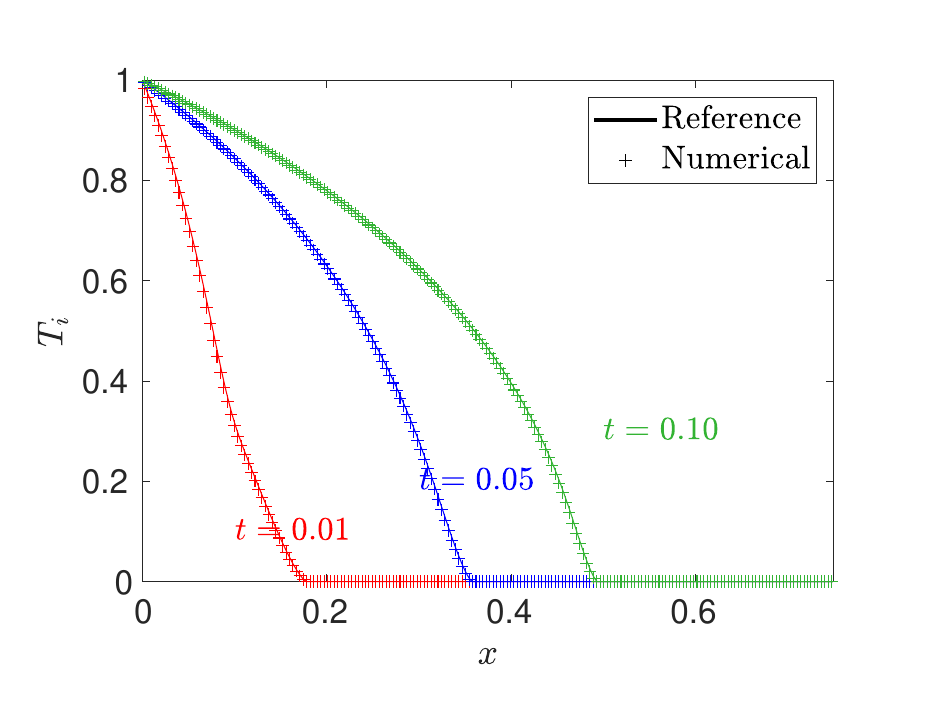}
}
\caption{(1D Marshak wave problem with conduction terms in Sec. \ref{sec:marshak_ex3}) The numerical solution of $T_r, T_e$ and $T_i$ for $\kappa = 1$ at different times. The symbols are the numerical solution by the AP $P_N$ method, and the solid lines are the reference solution obtained from the $S_N$ method.}
\label{fig:ex3_kappa1e0}
\end{figure}

In the simulation, the same grid size, expansion order, and time step length as in Sec. \ref{sec:marshak_ex2} are adopted. The numerical solution of the radiative temperature $T_r$, and the temperature of electron and ion $T_e, T_i$ for $\kappa = 1$ at $t = 0.01, 0.05$ and $0.1$ is plotted in Fig. \ref{fig:ex3_kappa1e0}. It shows that the numerical solution and the reference solution are almost on top of each other. Here the reference solution is obtained by the $S_N$ method. For the  $S_N$ method, the fully implicit scheme is applied with the mesh size in the spacial space $N = 1600$, and the CFL number $C = 0.8$. The upwind scheme is applied to solve the convection term, with $8$ Gaussian nodes utilized to discrete the angular space. Comparing Fig. \ref{fig:ex2_kappa1e0} and \ref{fig:ex3_kappa1e0}, the evolution of these three temperatures is quite different with and without the conduction terms. With the conduction terms, the evolution of these three temperatures is much smoother, aligning with the analytical properties of the system \eqref{eq:rte}. Furthermore, significant differences between $T_e$ and $T_i$ are noticed when $\kappa = 1$. 

\begin{figure}[!hptb]
\centering
\subfloat[$T_r, \kappa = 0.001$]{
\includegraphics[width=0.3\linewidth]{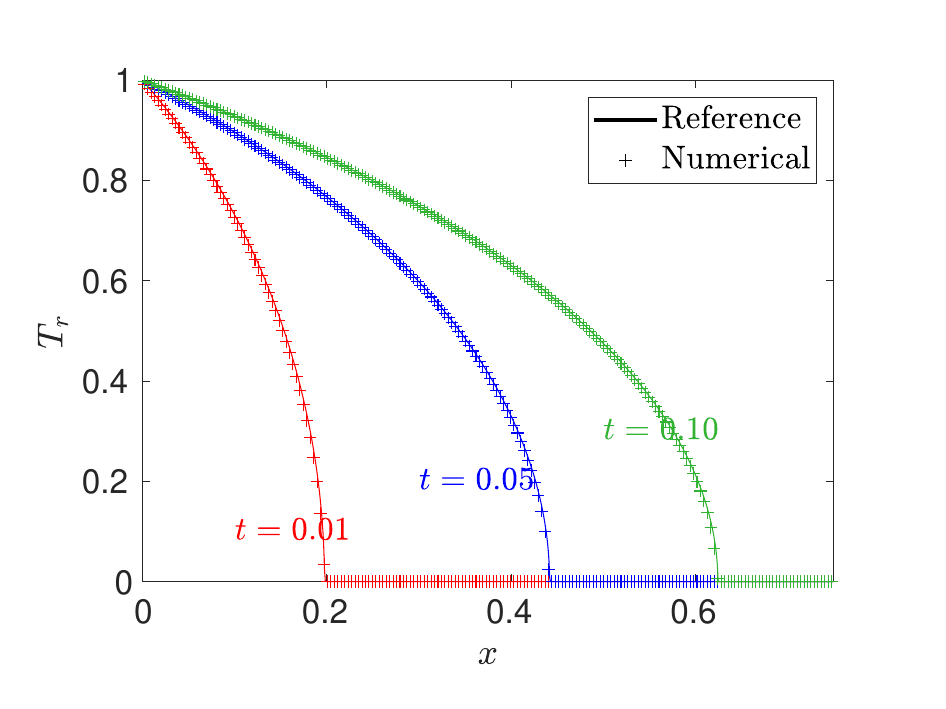}
}\hfill
\subfloat[$T_e, \kappa = 0.001$]{
\includegraphics[width=0.3\linewidth]{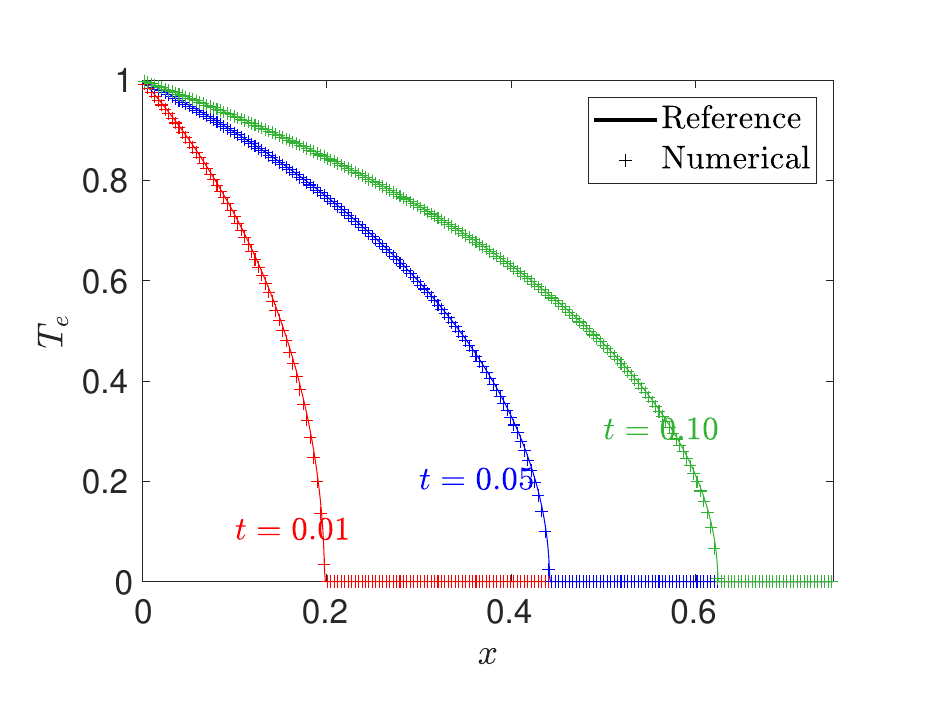}
\label{fig:ex3_kappa_0p001_Te}
}\hfill
\subfloat[$T_i, \kappa = 0.001$]{
\includegraphics[width=0.3\linewidth]{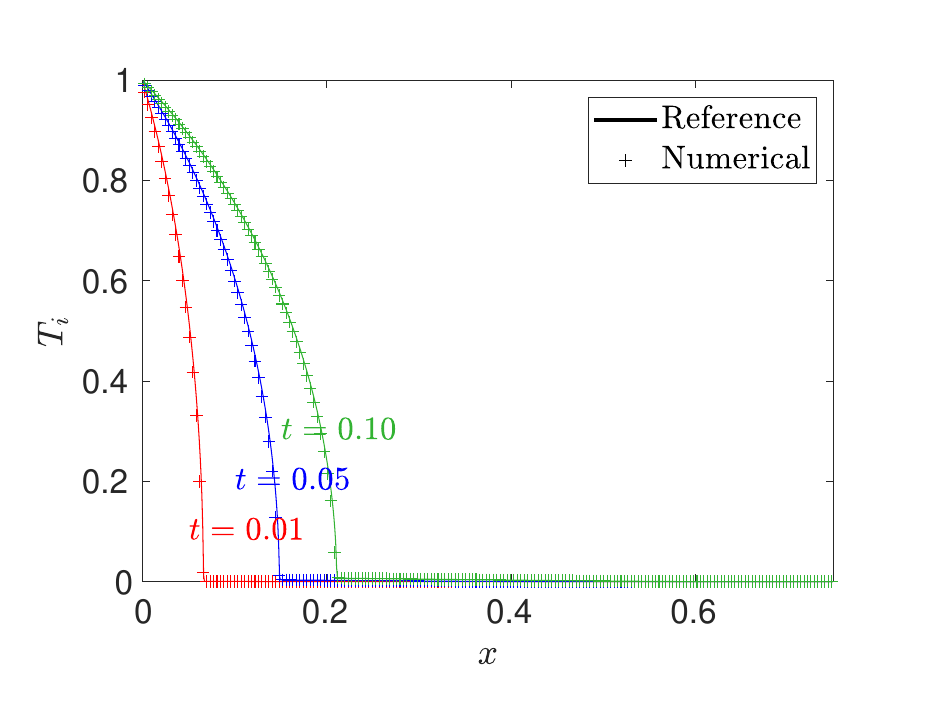}
\label{fig:ex3_kappa_0p001_Ti}
} \\
\subfloat[$T_r, \kappa = 100$]{
\includegraphics[width=0.3\linewidth]{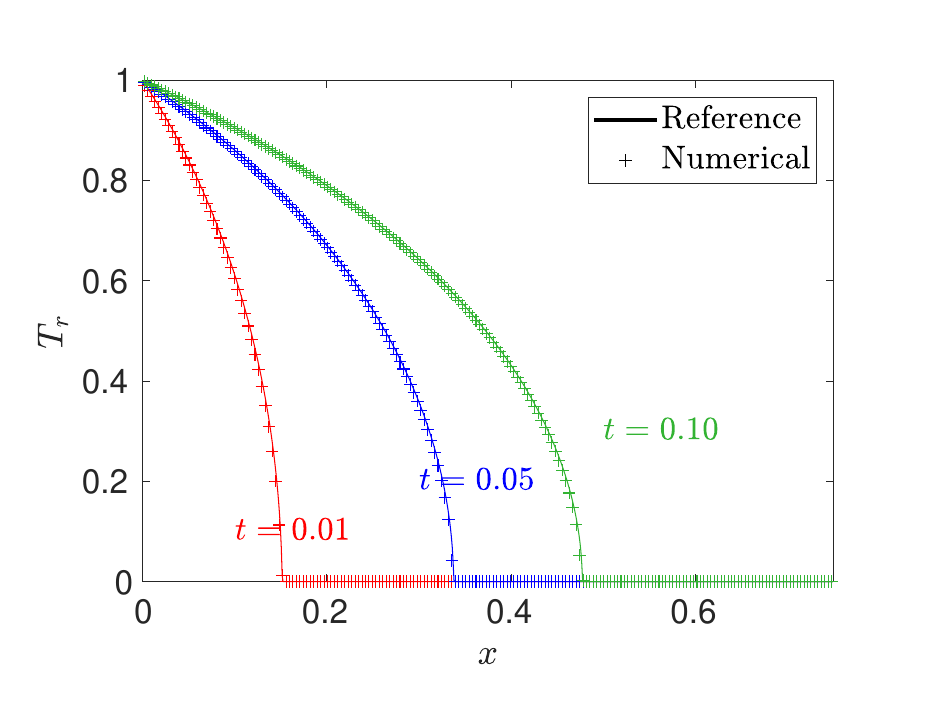}
}\hfill
\subfloat[$T_e, \kappa = 100$]{
\includegraphics[width=0.3\linewidth]{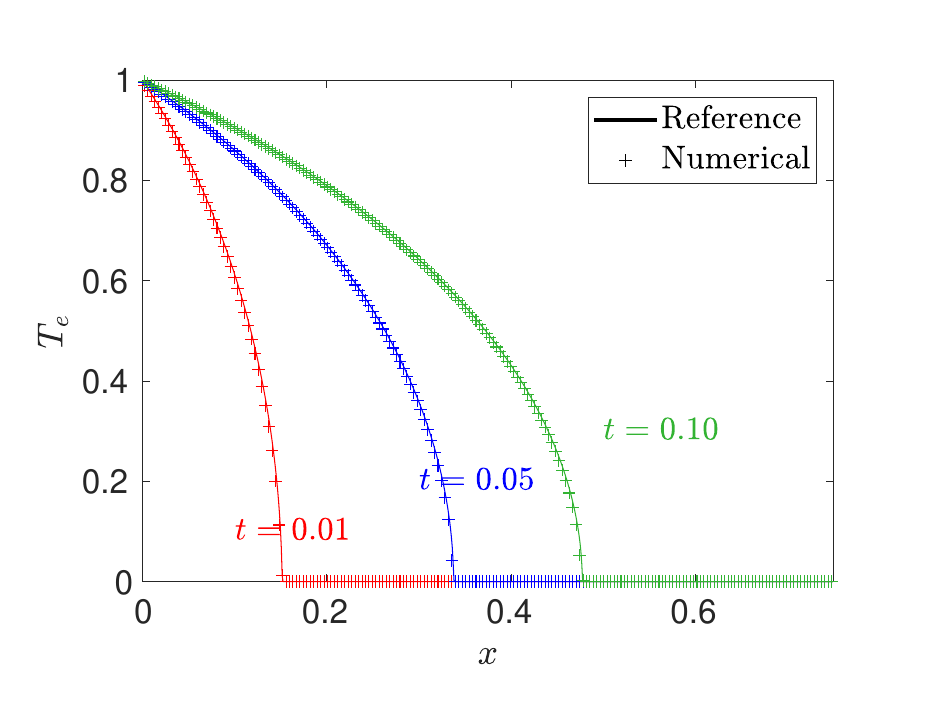}
}\hfill
\subfloat[$T_i, \kappa = 100$]{
\includegraphics[width=0.3\linewidth]{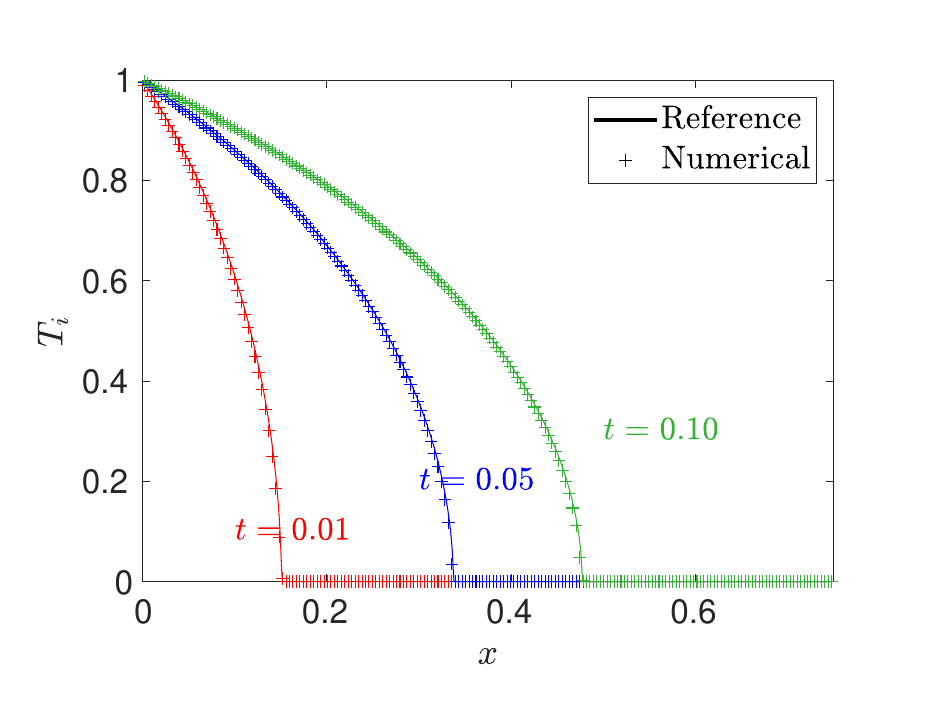}
\label{fig:ex3_kappa_100_Ti}
}
\caption{(1D Marshak wave problem with conduction terms in Sec. \ref{sec:marshak_ex3}) The numerical solution of $T_r, T_e$ and $T_i$ for $\kappa = 0.001$ and $100$ at different times. The top column is for $\kappa =0.001$ and the bottom column is for $\kappa = 100$. The symbols are the numerical solution by the AP $P_N$ method, and the solid lines are the reference solution obtained from the $S_N$ method.}
\label{fig:ex3_kappa}
\end{figure}

\begin{figure}[!hptb]
\centering
\subfloat[$\kappa = 0.001$]{
\includegraphics[width=0.3\linewidth]{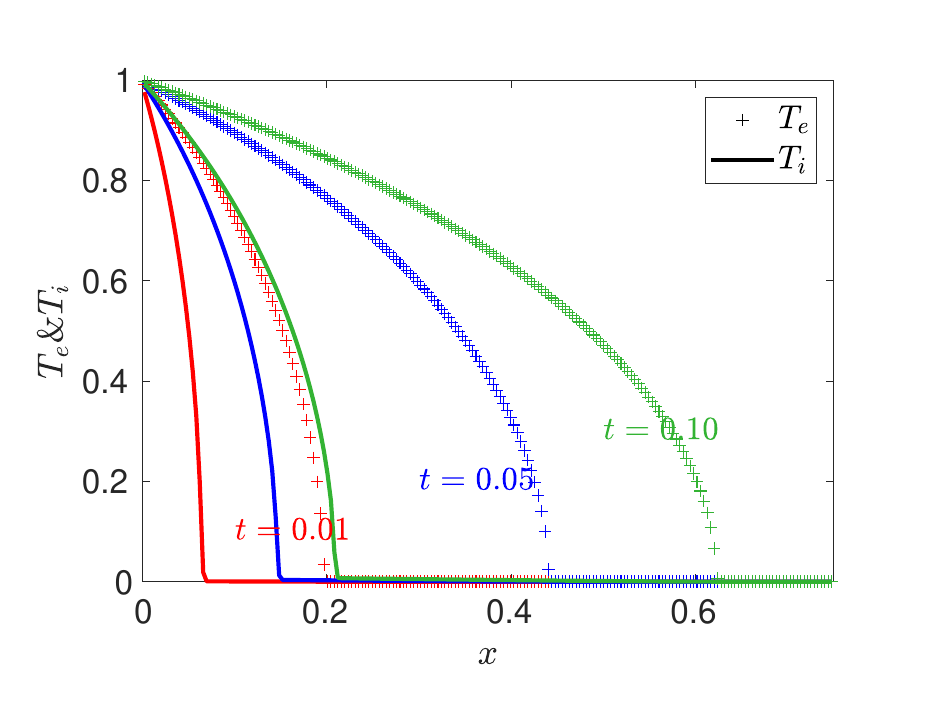}
\label{fig:ex3_kappa_compare_0p001}
}\hfill
\subfloat[$\kappa = 1$]{
\includegraphics[width=0.3\linewidth]{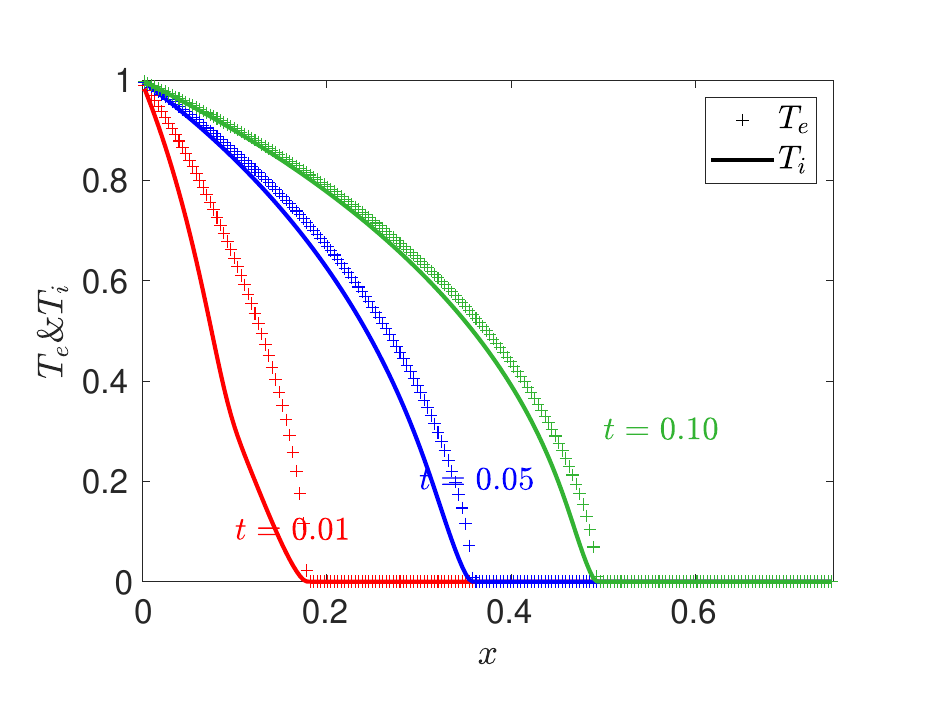}
\label{fig:ex3_kappa_compare_1}
}\hfill
\subfloat[$\kappa = 100$]{
\includegraphics[width=0.3\linewidth]{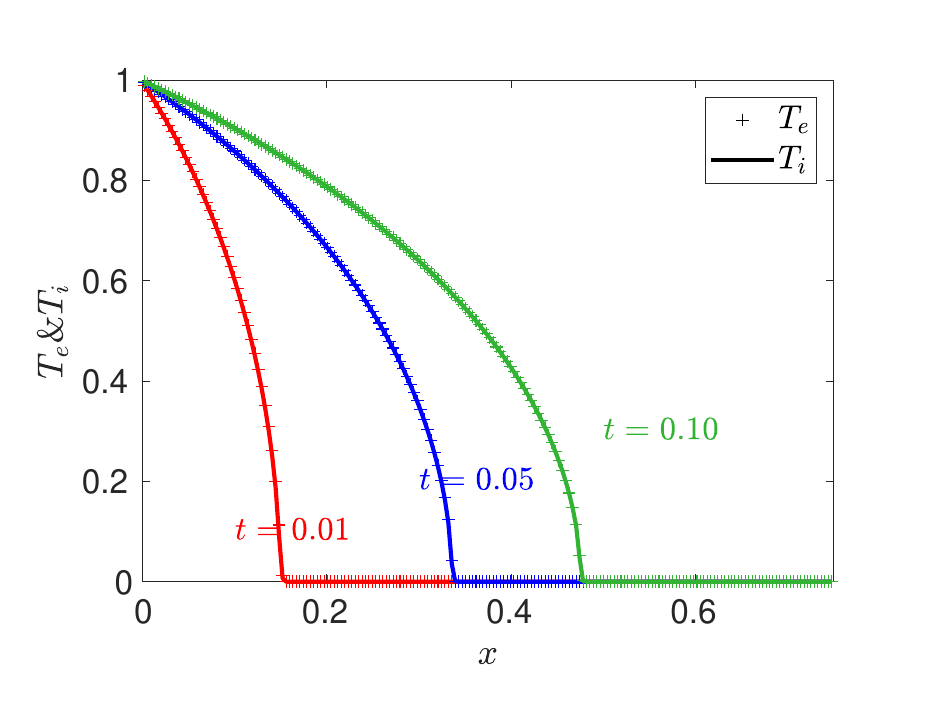}
\label{fig:ex3_kappa_compare_100}
}    
\caption{(1D Marshak wave problem with conduction terms in Sec. \ref{sec:marshak_ex3}) The comparison of $T_e$ and $T_i$ for different $\kappa$ at different times. The symbols are the numerical solution of $T_e$ and the solid lines are the numerical solution of $T_i$. (a) $\kappa = 0.001$, (b) $\kappa = 1$, (c) $\kappa = 100$.
}
\label{fig:ex3_kappa_compare}
\end{figure}

To demonstrate the AP property of this scheme, we have conducted tests for two additional values of $\kappa$ as $\kappa = 0.001$ and $100$. The numerical solutions for the three temperatures at $t = 0.01, 0.05$, and $0.1$ are presented in Fig. \ref{fig:ex3_kappa}, alongside the reference solution obtained by the $S_N$ method. For both values of $\kappa$, the numerical solutions closely align with the reference solutions, demonstrating the AP property of the numerical scheme introduced in Sec. \ref{sec:alg}. Furthermore, Fig. \ref{fig:ex3_kappa_0p001_Te} and \ref{fig:ex3_kappa_0p001_Ti} provide a clear depiction of the different behaviors between electrons and ions for small $\kappa$. However, as $\kappa$ increases, the behaviors of electrons and ions become more similar. To facilitate a precise comparison of their behavior, Fig. \ref{fig:ex3_kappa_compare} displays the temperature profiles of electrons and ions for different values of $\kappa$. It is evident that when $\kappa = 0.001$, the behaviors of the two temperatures differ significantly. In contrast, when $\kappa = 1$, the differences diminish considerably, and as $\kappa$ increases to $100$, the behaviors of these two temperatures nearly coincide. Notably, the same time step length, independent of $\kappa$, is employed for all three cases, indicating the AP property of this new method with respect to $\kappa$. 

\subsection{2D Riemann problem}
\label{sec:2D}
\begin{figure}[!htpb]
\centering
\subfloat[Problem setting for the 2D Riemann problem]{
\includegraphics[width=0.45\linewidth]{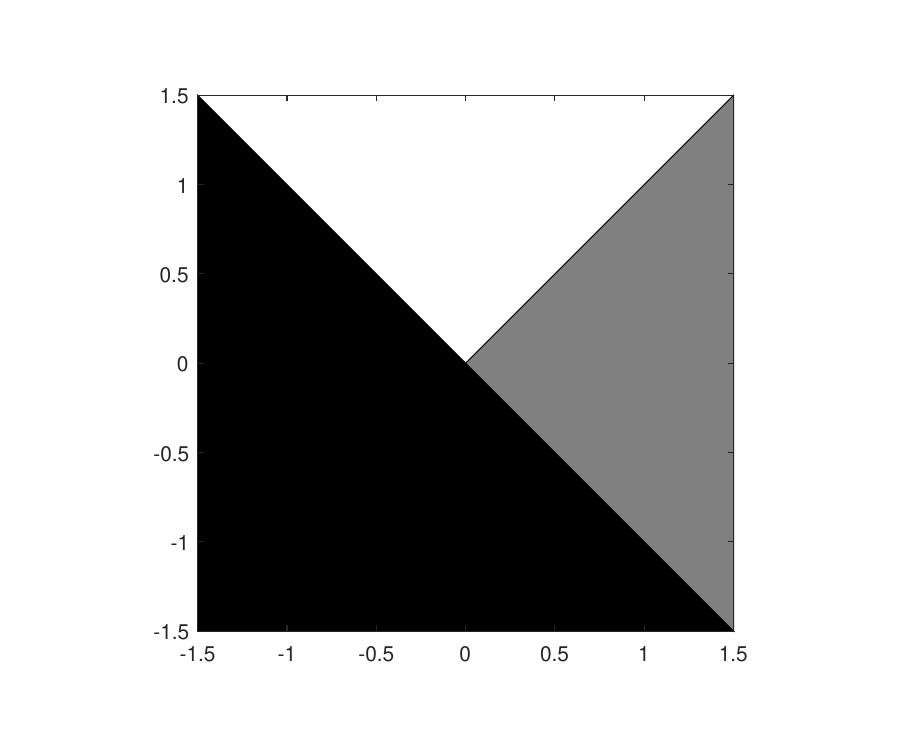}
\label{fig:2D_1}
}\hfill
\subfloat[Problem setting for the 1D Riemann problem]{
\includegraphics[width=0.45\linewidth]{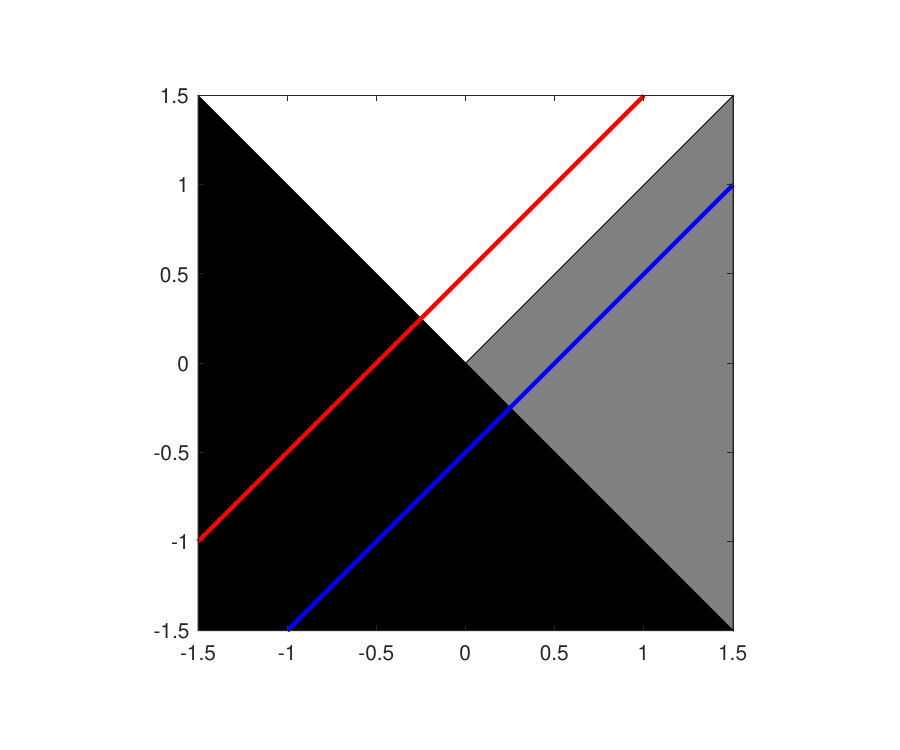}
\label{fig:2D_2}
}
\caption{(2D Riemann problem in Sec. \ref{sec:2D}) Computational regime and problem setting for the 2D Riemann problem. Here the blue and red line in Fig. \ref{fig:2D_2} is the computational regime in the corresponding 1D problem. }
\label{fig:2darea}
\end{figure}
In this section, a 2D Riemann problem is examined using a similar setup to the one described in \cite{XIONG2022111308}. This problem is defined on a square spatial domain $[-1.5,1.5]\times [-1.5,1.5]$ as shown in Fig. \ref{fig:2D_1}. In this example, the constant parameters are set as follows.
\begin{equation}
\label{eq:2D_para}
    a=1,\quad c=1,\quad C_{v_e}=1,\quad C_{v_i}=0.1.
\end{equation}
As in Fig. \ref{fig:2darea}, the opacity $\sigma$ is chosen based on the spatial region:  
\begin{equation}
    \sigma(x,y)=\left\{ 
    \begin{aligned}
    10,\quad &x+y>0 \text{ and } x<y,\quad\text{ (white region)},\\
    1,\quad &x+y<0 \text{ or } x>y, \quad\text{(black and gray regions)}.
    \end{aligned}
    \right.
\end{equation}
Moreover, the initial conditions for the material temperatures are given as
\begin{equation}
\label{eq:2dini}
    T_e(0,x,y)=T_i(0,x,y)=\left\{
    \begin{aligned}
    1,\quad & x+y<0, \quad\text{(black region)},\\
    0.1,\quad & x+y>0, \quad\text{(gray and white regions)}.
    \end{aligned}
    \right.
\end{equation}
The boundary conditions are set as
\begin{equation}
    T_e(t,x,y)=T_i(t,x,y)=\left\{
    \begin{aligned}
    1,\quad &x=-1.5\text{ or }y=-1.5,\\
    0.1,\quad &x=1.5\text{ or }y=1.5,
    \end{aligned}\qquad \forall t>0.
    \right.
\end{equation}
The initial and boundary radiative intensity is in thermal equilibrium
\begin{equation}
    \begin{gathered}
    \psi(0,x,y,\Omega)=\frac{1}{4\pi}acT_r^4(0,x,y),\\ \psi(t,x,y,\Omega)=\frac{1}{4\pi}acT_r^4(t,x,y),\quad  (x,y)\in\partial [-1.5,1.5]^2,
    \end{gathered}
\end{equation}
with $T_r(0,x,y)=T_e(0,x,y)$ given in \eqref{eq:2dini}, and $T_r(t,x,y)=T_e(t,x,y)$ for boundary condition. 

\begin{figure}[!htpb]
\centering
\subfloat[$T_r, P_N$]{
\includegraphics[width=0.3\linewidth]{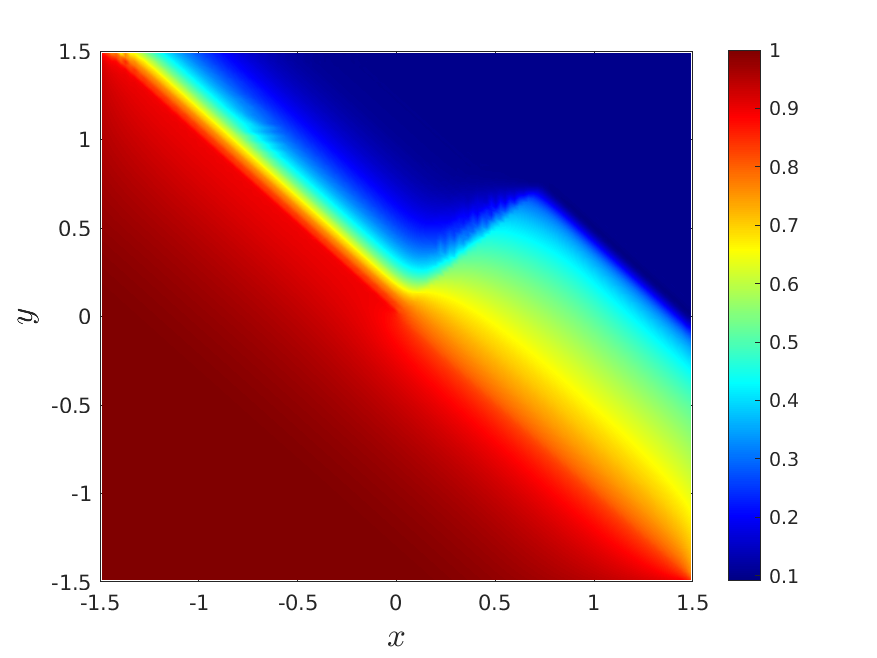}
} \hfill 
\subfloat[$T_e, P_N$]{
\includegraphics[width=0.3\linewidth]{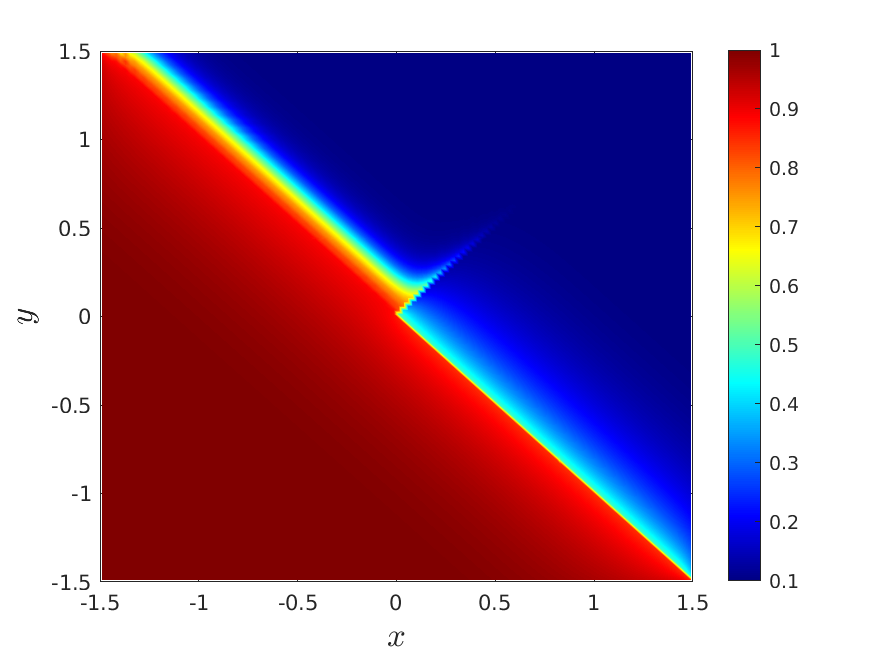}
}\hfill
\subfloat[$T_i, P_N$]{
\includegraphics[width=0.3\linewidth]{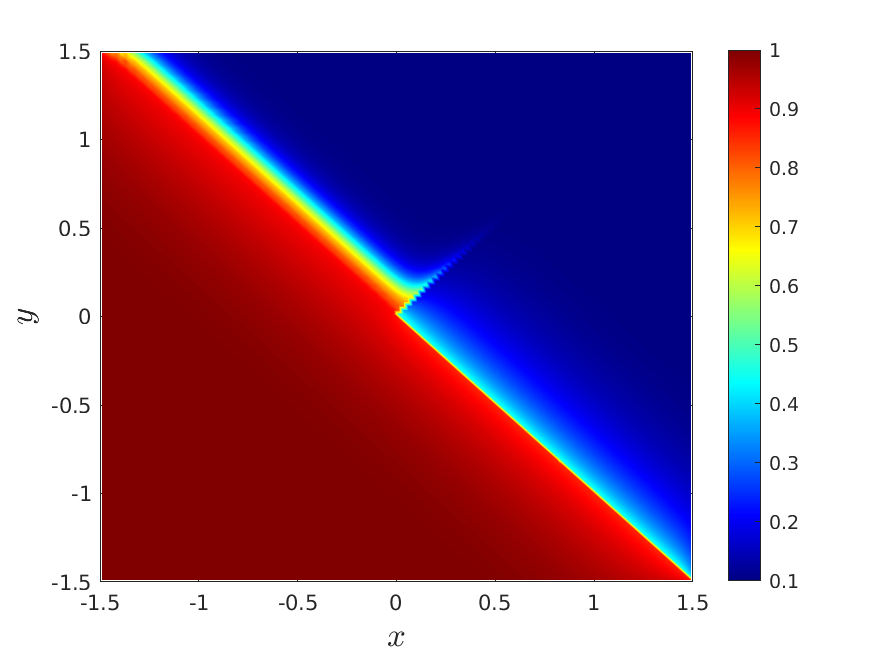}
}\\
\subfloat[$T_r, S_N$]{
\includegraphics[width=0.3\linewidth]{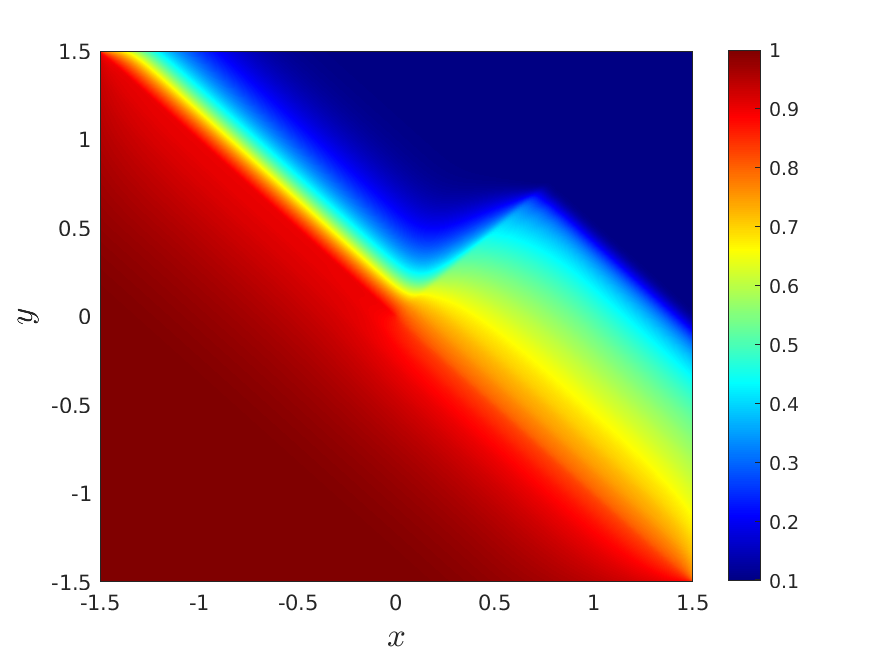}
} \hfill 
\subfloat[$T_e, S_N$]{
\includegraphics[width=0.3\linewidth]{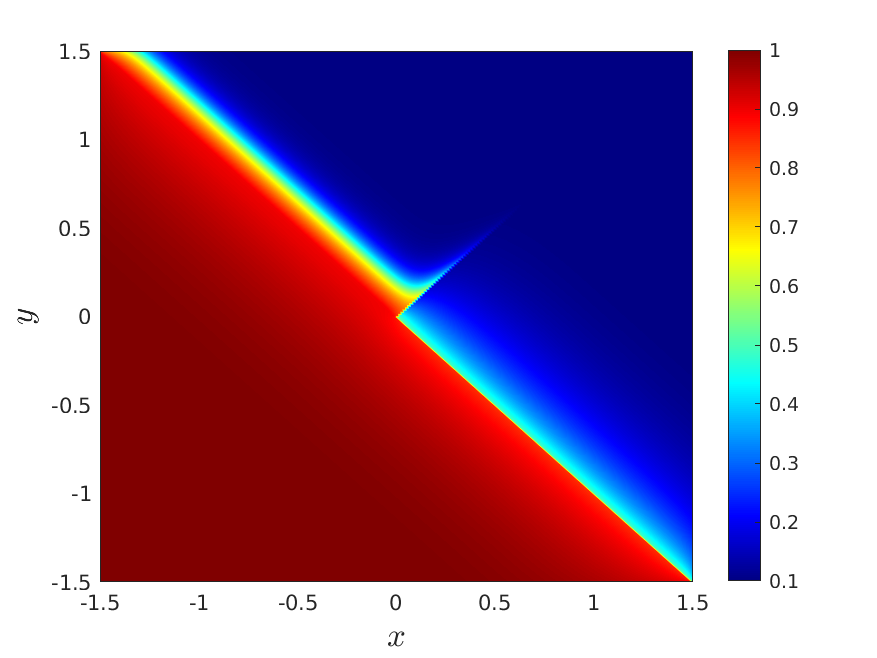}
}\hfill
\subfloat[$T_i, S_N$]{
\includegraphics[width=0.3\linewidth]{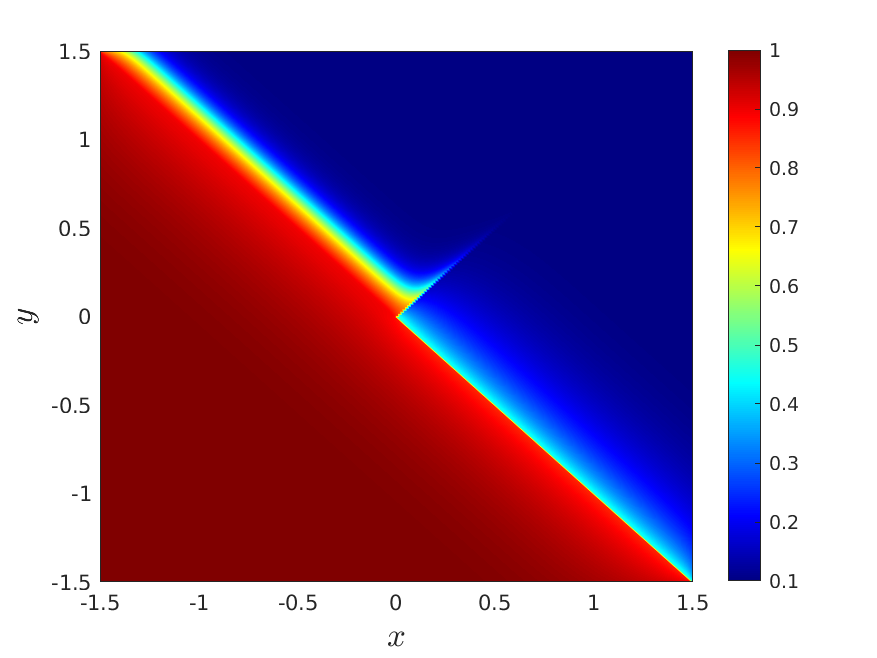}
}
\caption{(2D Riemann problem in Sec. \ref{sec:2D}) The radiative temperature $T_r$ and the temperature of electron and ion $T_e,T_i$ at $t = 1$ for $\kappa = 1$. The first column is the numerical solution obtained by the AP $P_N$ method, and the second column is the reference solution obtained by the $S_N$ method.}
\label{fig:2D_sol}
\end{figure}

\begin{figure}[!htpb]
\centering
\subfloat[$T_r, \kappa = 1$]{
\includegraphics[width=0.3\linewidth]{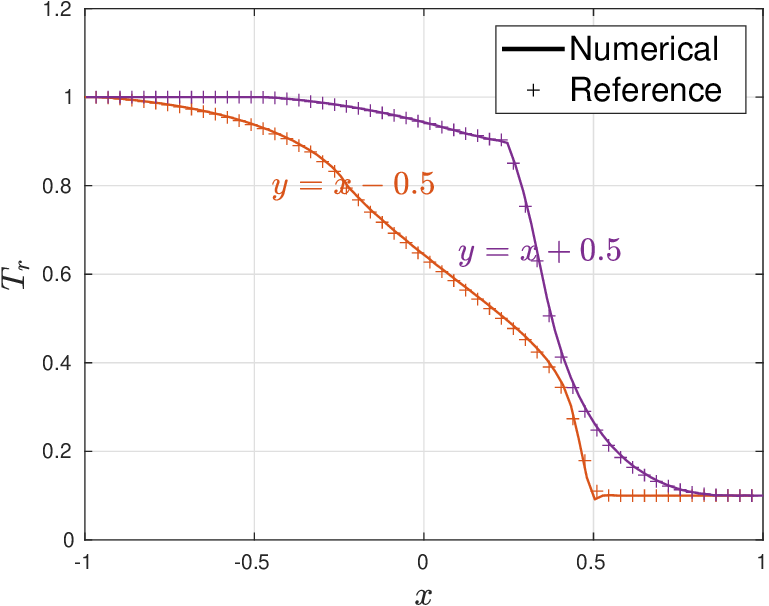}
} \hfill 
\subfloat[$T_e, \kappa = 1$]{
\includegraphics[width=0.3\linewidth]{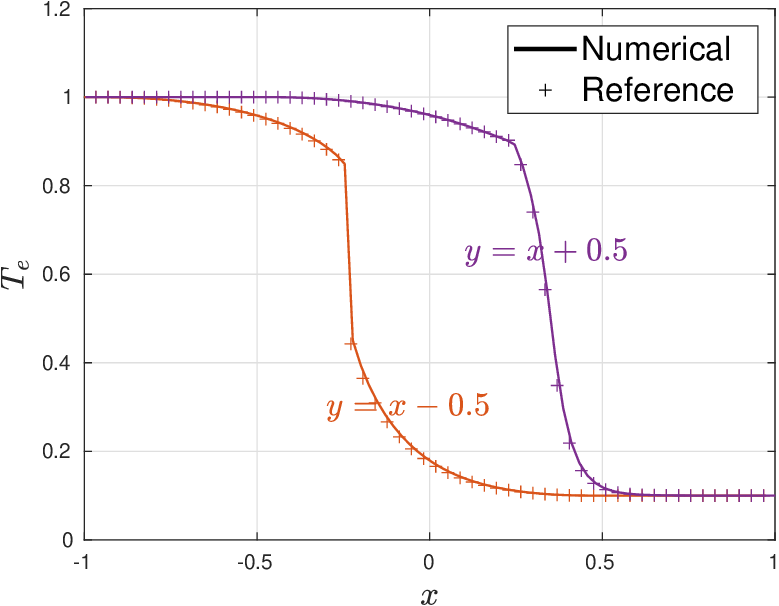}
}\hfill
\subfloat[$T_i, \kappa = 1$]{
\includegraphics[width=0.3\linewidth]{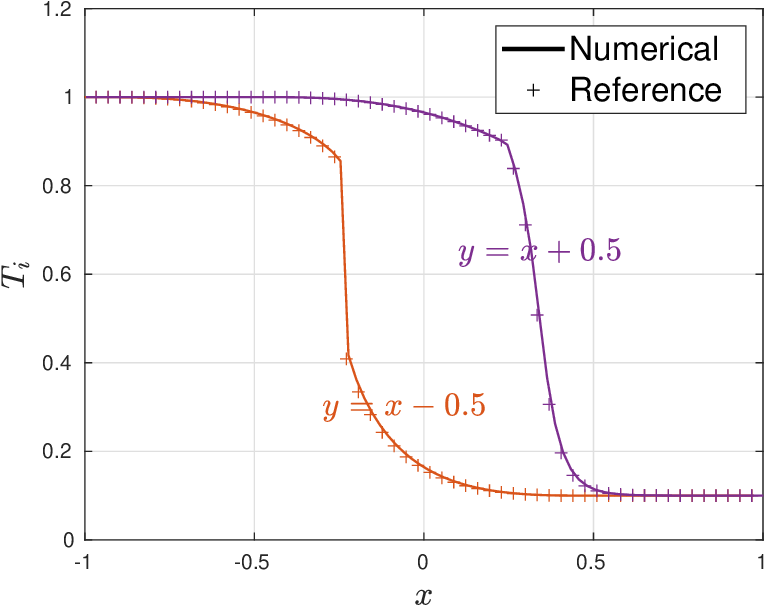}
} 
\caption{(2D Riemann problem in Sec. \ref{sec:2D}) The radiative temperature $T_r$ and the temperature of electron and ion $T_e, T_i$ at $t = 1$ along $y = x - 0.5$ and $y = x+0.5$ for $\kappa = 1$. The solid line is the numerical solution obtained by the AP $P_N$ method, and the symbol is the reference solution obtained by the $S_N$ method.}
\label{fig:2D_sol_local}
\end{figure}

\begin{figure}[!htpb]
\centering
\subfloat[$T_r, \kappa = 1$]{
\includegraphics[width=0.3\linewidth]{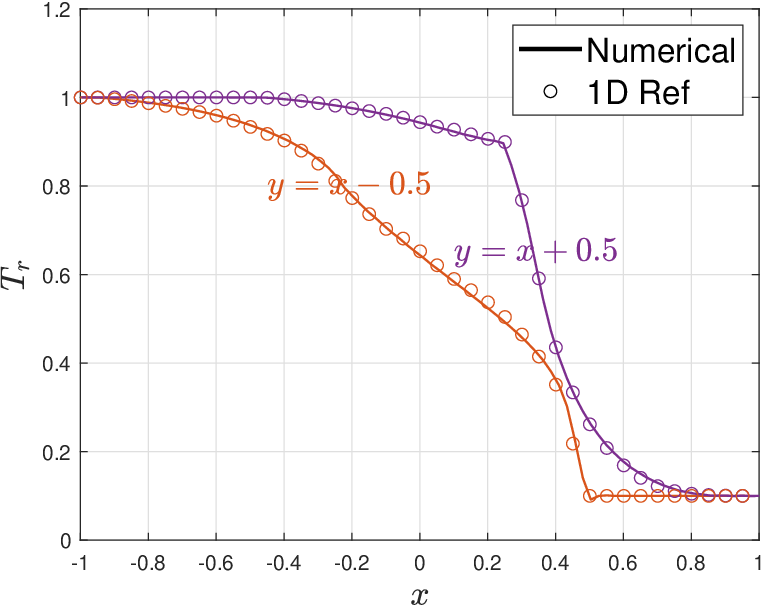}
} \hfill 
\subfloat[$T_e, \kappa = 1$]{
\includegraphics[width=0.3\linewidth]{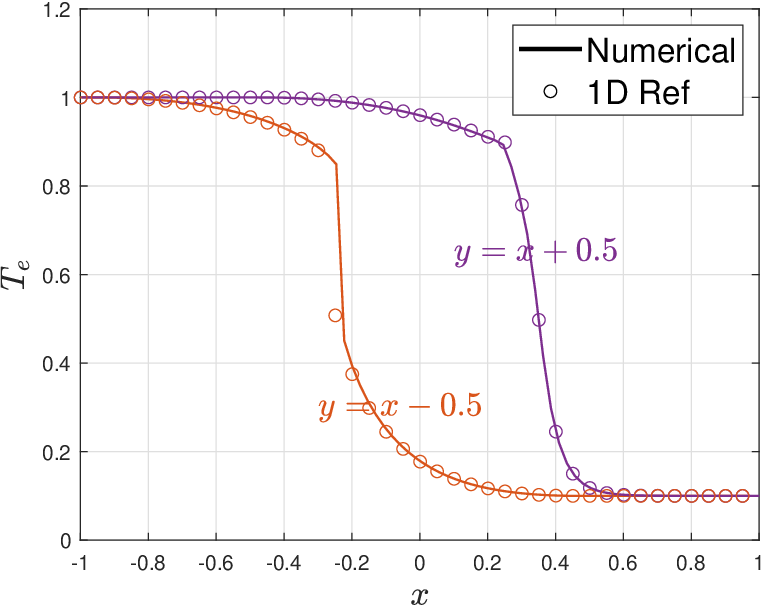}
}\hfill
\subfloat[$T_i, \kappa = 1$]{
\includegraphics[width=0.3\linewidth]{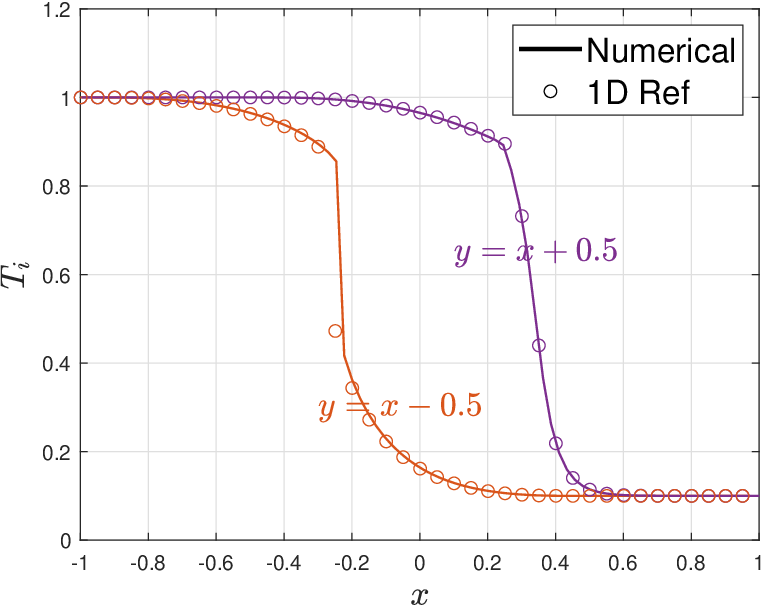}
} 
\caption{(2D Riemann problem in Sec. \ref{sec:2D}) Comparision of the numerical solution with the 1D problem. The radiative temperature $T_r$ and the temperature of electron and ion $T_e, T_i$ at $t = 1$ for $\kappa = 1$. The solid line is the numerical solution obtained by the AP $P_N$ method along $y = x - 0.5$ and $y = x + 0.5$, and the symbol is the numerical solution obtained by the 1D problem.}
\label{fig:2D_1D_sol_local}
\end{figure}

In the simulation, the $\Pn$ method employs an expansion order of $M = 30$, and the grid size is set to $N_x = N_y = 128$ with WENO reconstruction \cite{WENO} in \ref{app:weno}. First, we consider the case $\kappa = 1$. The radiative temperature $T_r$, electron temperature $T_e$, and ion temperature $T_i$ at $t = 1$ are presented in Fig. \ref{fig:2D_sol}, where the reference solution is obtained using the $S_N$ method with a mesh size of $256 \times 256$. In the $S_N$ method, the angular space is discreted with $1202$ Lebedev points (order $59$ polynomial), and the upwind scheme is utilized for the convection term with the time step length \eqref{eq:CFL} with $\rm CFL = 0.1$. It is evident that the numerical solution closely matches the reference solution. Moreover, in the white region, the radiative temperature is lower than in the gray region, consistent with the opacity value $\sigma$. To provide a more detailed comparison between the numerical and reference solutions, we plot the numerical solution along the lines $y = x - 0.5$ and $y = x + 0.5$ in Fig. \ref{fig:2D_sol_local}, demonstrating a close agreement. Additionally, to verify the numerical solution more carefully, two one-dimensional problems along the regime shown in Fig. \ref{fig:2D_2} are tested with the initial and boundary conditions the same as the 2D problem, where the same $P_N$ method is utilized with the expansion number $M = 30$, and the grid size $N = 1600$. The numerical results are shown in Fig. \ref{fig:2D_1D_sol_local}, where the behavior of the numerical solution for the 1D problem is almost the same as the 2D numerical solution along $y = x - 0.5$ and $y = x + 0.5$. This indicates that a smaller opacity $\sigma$ accelerates the propagation of radiation, yet decelerates the heating of the material by the radiation. Note that in the problems considered in this section, the conduction term is absent, meaning the change in material temperature is solely influenced by radiation. Therefore, for a small opacity $\sigma$, the material temperature rises over a large spatial range, but the increase is relatively small. This shows that the numerical results are reasonable. 

\begin{remark}
In this example, the interface for the one-dimensional test is located at the cell boundary, and the interface for the two-dimensional test passes through the center of the cell. 
\end{remark}

\begin{figure}[!htpb]
\centering
\subfloat[$\kappa = 0$]{
\includegraphics[width=0.3\linewidth]{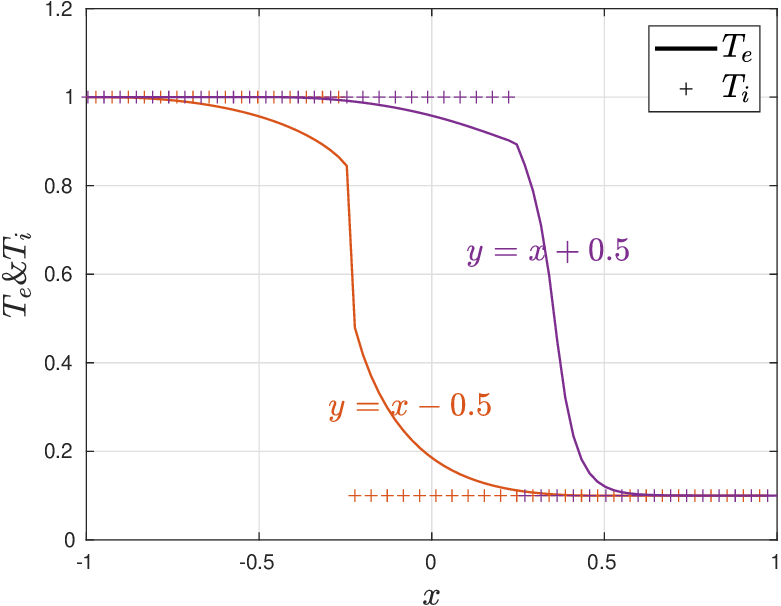}
\label{fig:2D_sol_compare_1}
} \hfill 
\subfloat[$\kappa = 0.01$]{
\includegraphics[width=0.3\linewidth]{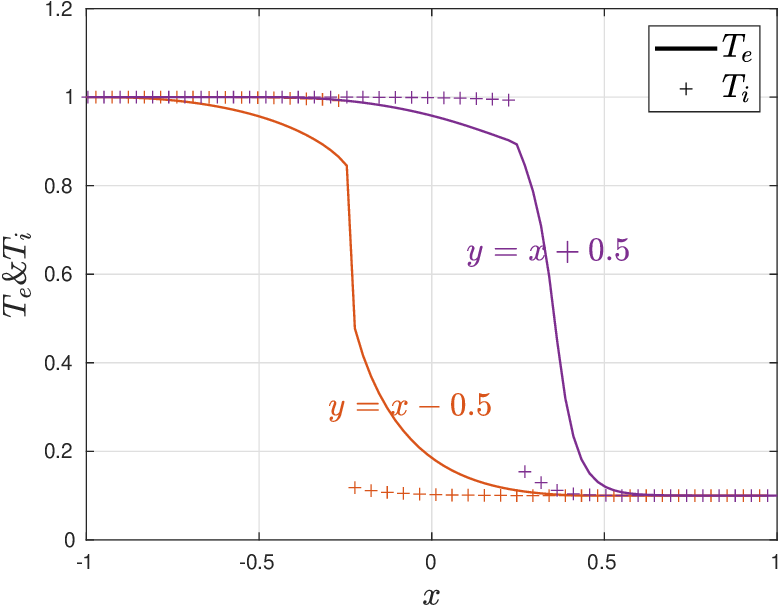}
\label{fig:2D_sol_compare_2}
}\hfill
\subfloat[$\kappa = 0.1$]{
\includegraphics[width=0.3\linewidth]{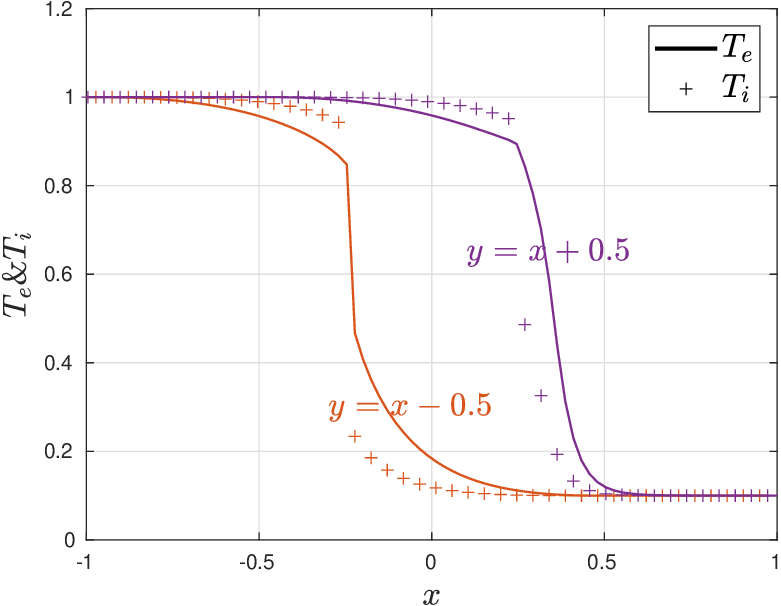}
\label{fig:2D_sol_compare_3}
} \\
\subfloat[$\kappa = 1$]{
\includegraphics[width=0.3\linewidth]{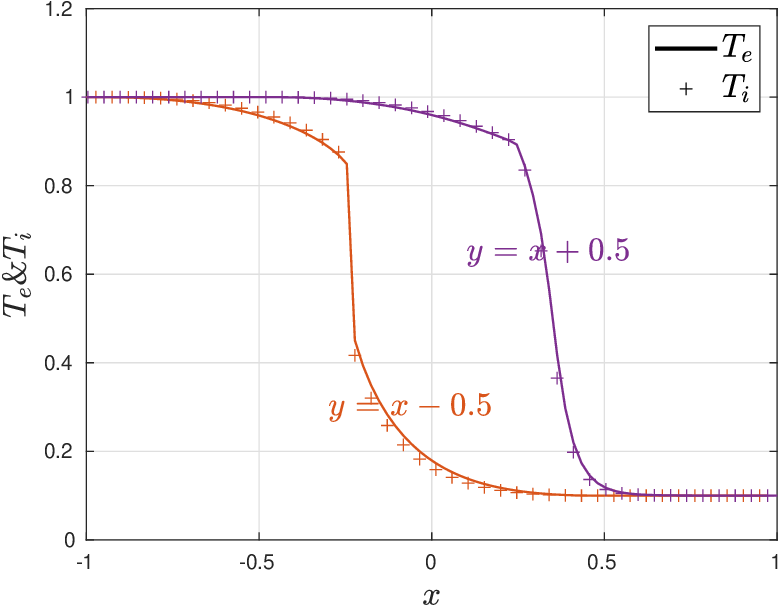}
\label{fig:2D_sol_compare_4}
} \hfill 
\subfloat[$\kappa = 10$]{
\includegraphics[width=0.3\linewidth]{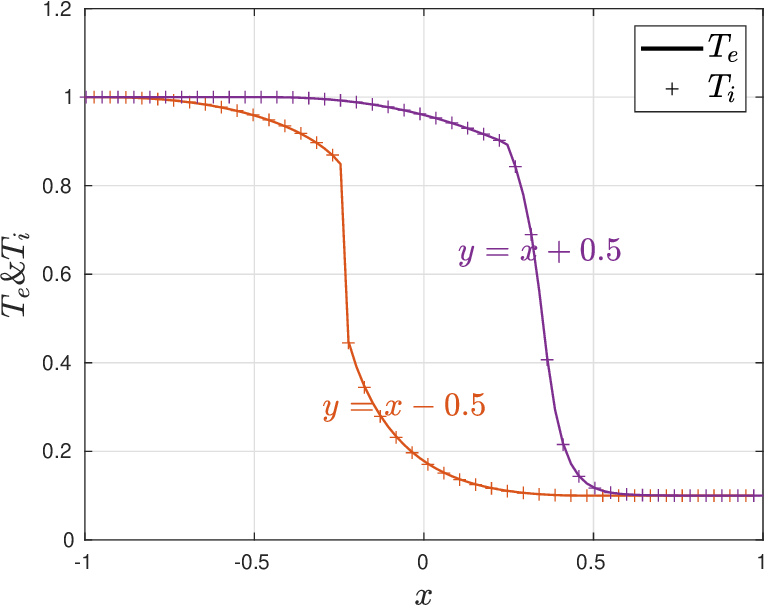}
\label{fig:2D_sol_compare_5}
}\hfill
\subfloat[$\kappa = 100$]{
\includegraphics[width=0.3\linewidth]{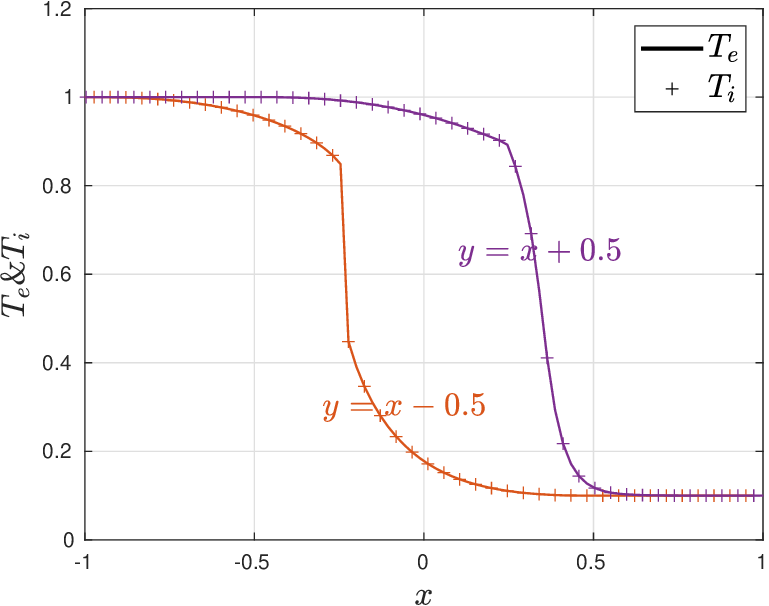}
\label{fig:2D_sol_compare_6}
} 
\caption{(2D Riemann problem in Sec. \ref{sec:2D}) The temperature $T_e, T_i$ at $t = 1$ along $y = x - 0.5$ and $y = x+0.5$ for $\kappa = 0, 0.01, 0.1, 1, 10$ and $100$. The solid line is the numerical solution of the temperature $T_e$, and the symbol is the numerical solution of the temperature $T_i$.}
\label{fig:2D_sol_compare}
\end{figure}

To investigate the impact of $\kappa$ on the temperatures $T_e$ and $T_i$, the numerical solutions for different $\kappa$ along the same lines $y = x - 0.5$ and $y = x + 0.5$ are presented in Fig. \ref{fig:2D_sol_compare}. When $\kappa = 0$, the system reduces to a two-temperature limit, where $T_i$ remains the same as the initial condition and significantly differs from $T_e$ at $t = 1$, as shown in Fig. \ref{fig:2D_sol_compare_1}. As $\kappa$ increases, the variation in $T_i$ becomes more pronounced. As illustrated in Fig. \ref{fig:2D_sol_compare_2} and \ref{fig:2D_sol_compare_3}, the numerical solution of $T_i$ gradually converges towards that of $T_e$. When $\kappa$ reaches $\kappa = 1$, the discrepancy between these two numerical solutions is small. Moreover, at $\kappa = 10$ and $100$, these two solutions nearly coincide. Fig. \ref{fig:2D_sol_compare} demonstrates that the behavior of temperatures $T_e$ and $T_i$ aligns with the AP property of the system \eqref{eq:rte} as $\kappa$ increases.

\begin{figure}[!htpb]
\centering
\subfloat[$T_r, P_N$]{
\includegraphics[width=0.3\linewidth]{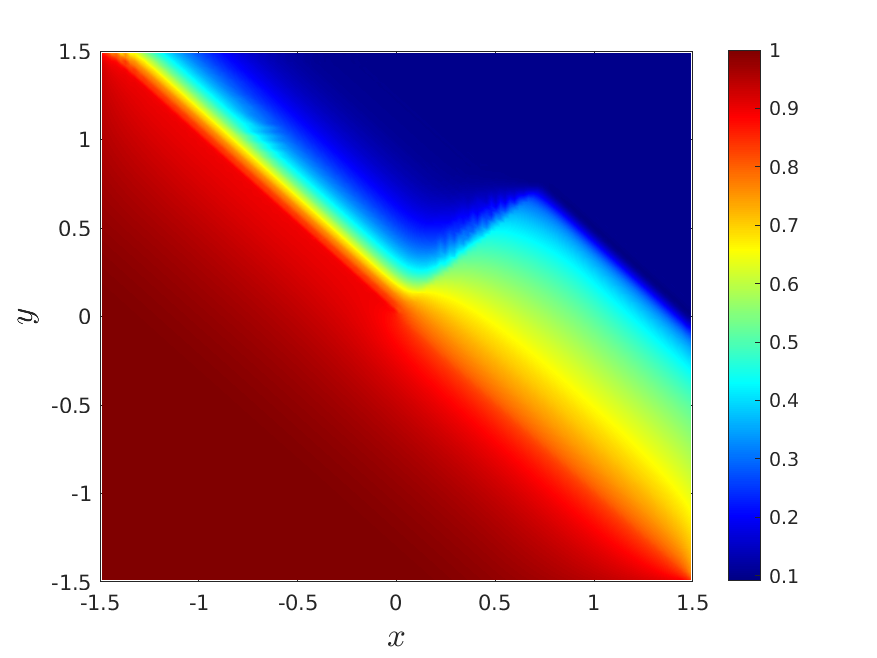}
} \hfill 
\subfloat[$T_e, P_N$]{
\includegraphics[width=0.3\linewidth]{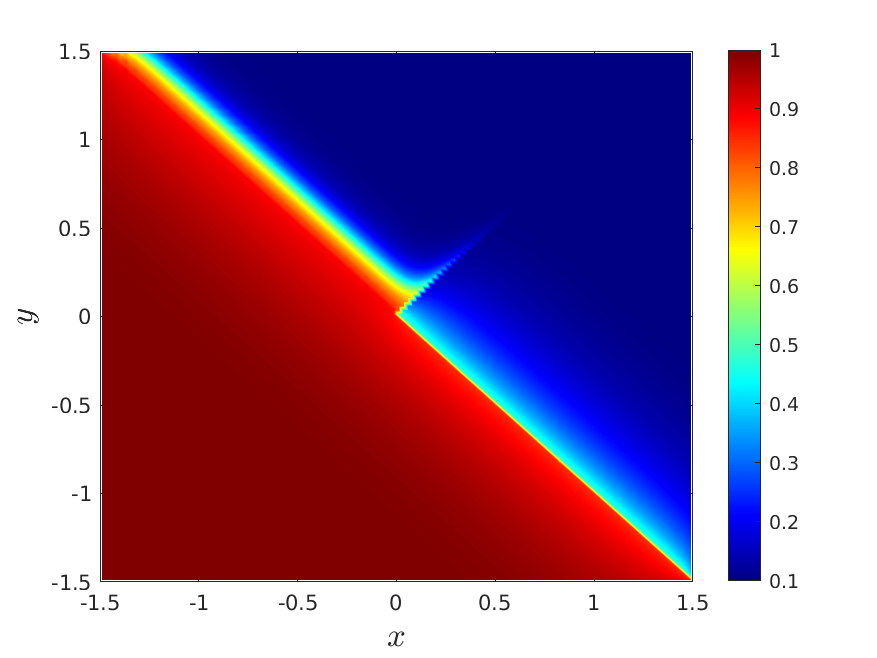}
}\hfill
\subfloat[$T_i, P_N$]{
\includegraphics[width=0.3\linewidth]{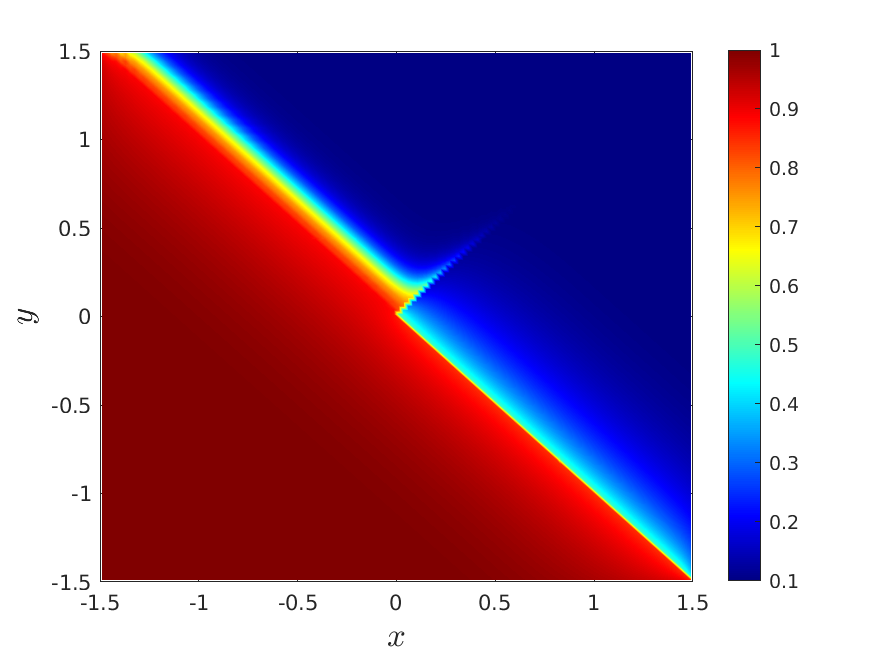}
}\\
\subfloat[$T_r, S_N$]{
\includegraphics[width=0.3\linewidth]{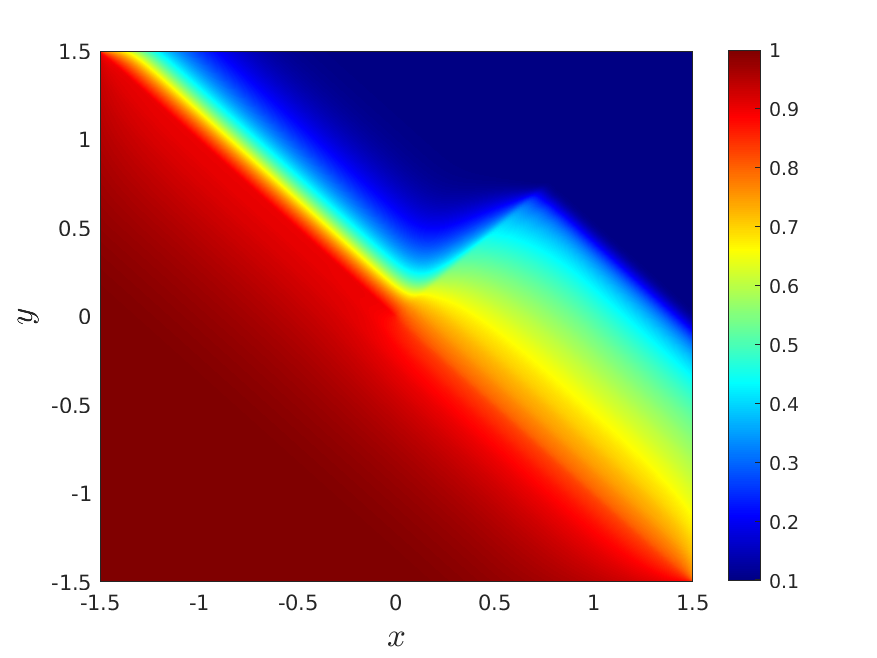}
} \hfill 
\subfloat[$T_e, S_N$]{
\includegraphics[width=0.3\linewidth]{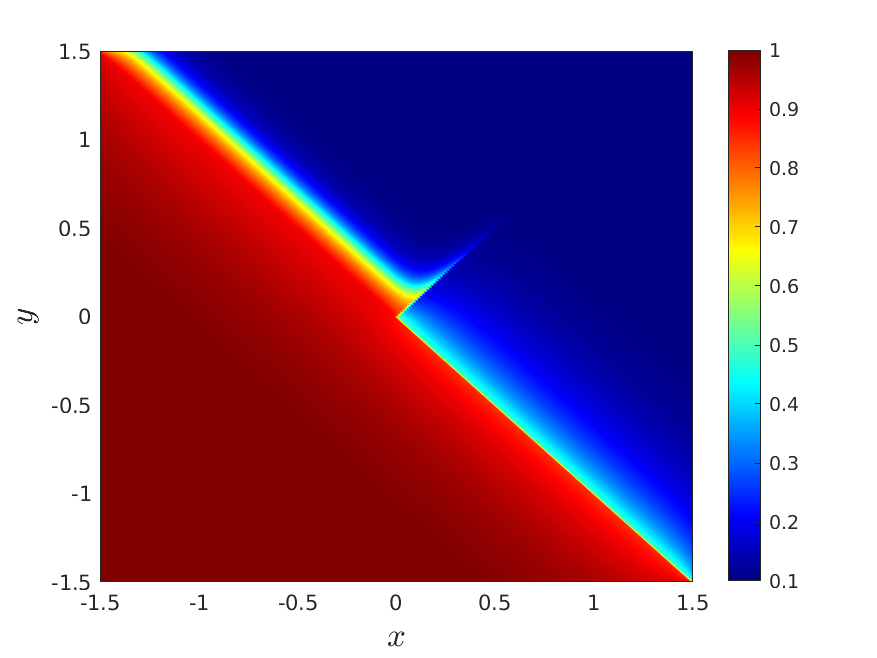}
}\hfill
\subfloat[$T_i, S_N$]{
\includegraphics[width=0.3\linewidth]{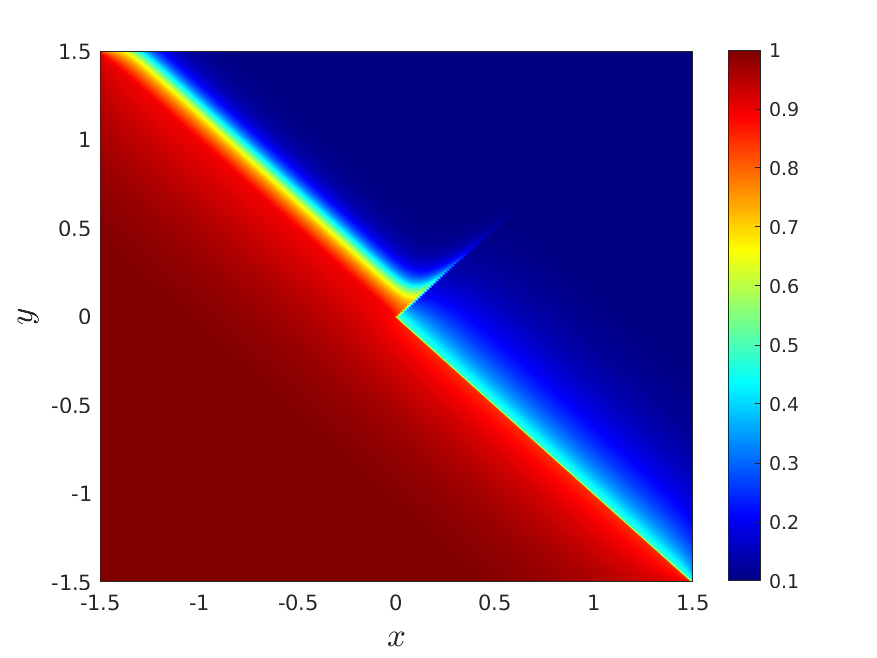}
}
\caption{(2D Riemann problem in Sec. \ref{sec:2D}) The radiative temperature $T_r$ and the temperature of electron and ion $T_e, T_i$ at $t = 1$ for $\kappa = 1000$. The first column is the numerical solution obtained by the AP $P_N$ method, and the second column is the reference solution obtained by the $S_N$ method.}
\label{fig:2D_sol_1000}
\end{figure}

\begin{figure}[!htpb]
\centering
\subfloat[$T_r, \kappa = 1000$]{
\includegraphics[width=0.3\linewidth]{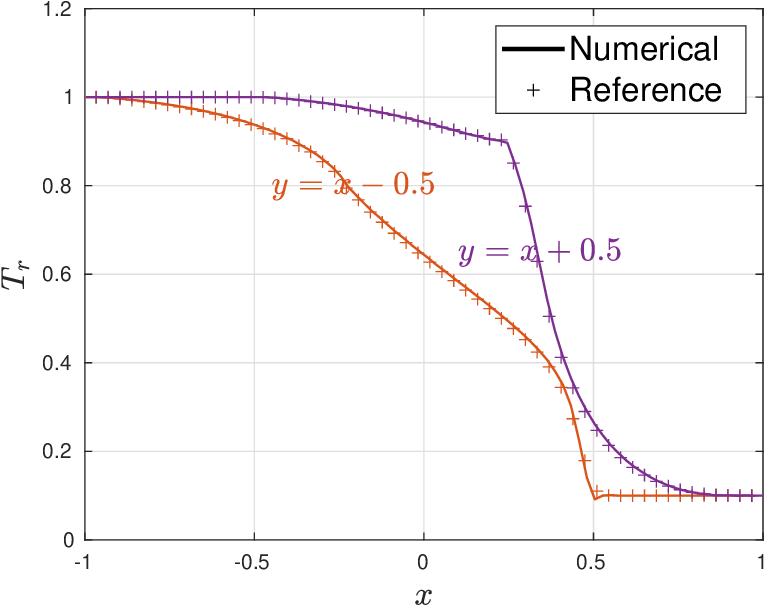}
} \hfill 
\subfloat[$T_e, \kappa = 1000$]{
\includegraphics[width=0.3\linewidth]{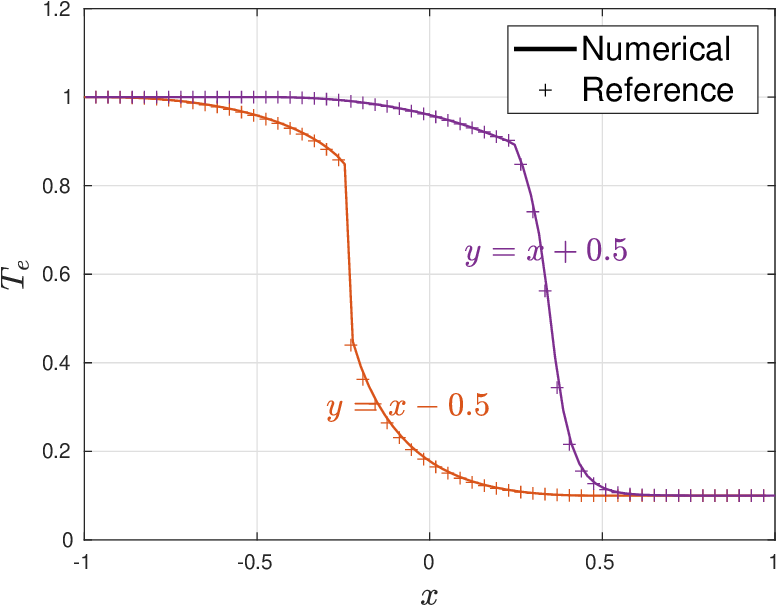}
}\hfill
\subfloat[$T_i, \kappa = 1000$]{
\includegraphics[width=0.3\linewidth]{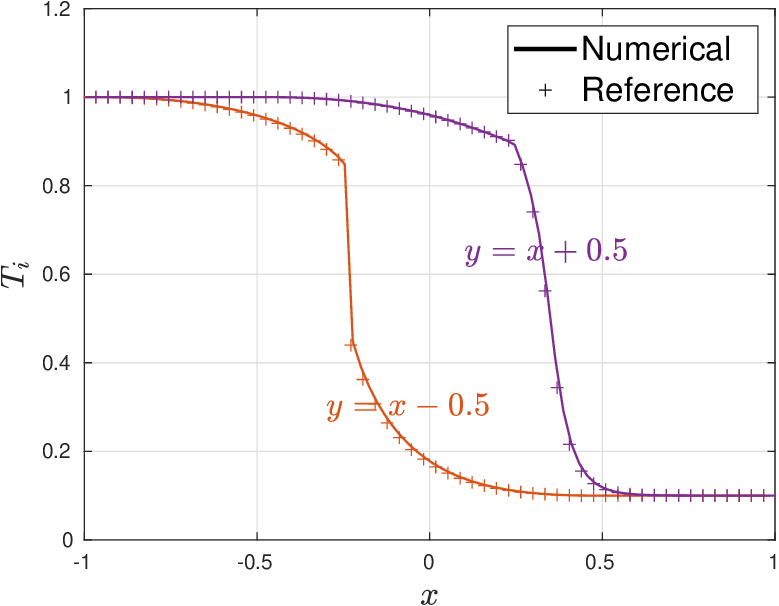}
} 
\caption{(2D Riemann problem in Sec. \ref{sec:2D}) The radiative temperature $T_r$ and the temperature of electron and ion $T_e, T_i$ at $t = 1$ along $y = x - 0.5$ and $y = x+0.5$ for $\kappa = 1000$. The solid line is the numerical solution obtained by the AP $P_N$ method, and the symbol is the reference solution obtained by the $S_N$ method.}
\label{fig:2D_sol_1000_local}
\end{figure}

To verify the AP property of this numerical scheme, the numerical solution is examined for an extremely large $\kappa$ value, specifically $\kappa = 1000$. The same expansion order for the $\Pn$ method, mesh size, and time step size are utilized. The numerical solutions of $T_r, T_e$, and $T_i$ at $t = 1$ are presented in Fig. \ref{fig:2D_sol_1000}, along with the reference solution obtained using the $S_N$ method. As depicted in Fig. \ref{fig:2D_sol_1000}, the numerical solution closely aligns with the reference solution, even with $\kappa = 1000$. Furthermore, Fig. \ref{fig:2D_sol_1000_local} displays the numerical solution along $y = x - 0.5$ and $y = x + 0.5$. The comparison with the reference solution further illustrates the AP property of the numerical scheme introduced in Sec. \ref{sec:alg}.

\section{Conclusion}
\label{sec:conclusion}

In this work, we have introduced an AP numerical method for solving the three-temperature radiative transfer model. This method utilizes a splitting scheme to ensure the accuracy and stability of the solution in both the diffusion limit as $\epsilon$ approaches zero and the two-temperature limit as $\kappa$ tends to infinity. Additionally, the numerical scheme has been demonstrated to possess the AP property and conserve energy. We have presented various numerical examples, including homogeneous model problems, the 1D Marshak wave problem, and the 2D Riemann problems, to validate the AP property and the efficiency of the proposed numerical scheme. Future work may involve exploring more complex numerical tests and applications of this method in different physical contexts such as the frequency-dependent problems \cite{zhang2023fully}. Furthermore, this splitting method can be generalized to provide an asymptotic-preserving coupling method of radiation and fluid, as well as radiation and plasma.

\section*{Acknowledges}
We thank Prof. Peng Song from the Institute of Applied Physics and Computational Mathematics and Prof. Tao Xiong from Xiamen University for their valuable suggestions. The work of Ruo Li is partially supported by the National Natural Science Foundation of China (Grant No. 12288101). W. Li is partially supported by the National Natural Science Foundation of China (Grant No. 12001051). M. Tang is partially supported by the National Natural Science Foundation of China (Grant No. 12031013), Shanghai pilot innovation project 21JC1403500, and the Strategic Priority Research Program of Chinese Academy of Sciences Grant No. XDA25010401. This work of Y. Wang is partially supported by the Foundation of the President of China Academy of Engineering Physics (YZJJZQ2022017) and 
the National Natural Science Foundation of China (Grant No. 12171026, U2230402, and 12031013).

\appendix
\section{Order analysis}
\label{app:order_analy}
In this section, we present an order analysis of $\psi_l$, where $l = 0, 1,\cdots $. A similar analysis is available in \cite{Fu2022}, and we included it here for the completeness of the algorithm. As stated in \cite{Fu2022}, it is confirmed that for $\psi_l$, the following relationship holds:
\begin{equation}
    \label{eq:oa}
    \psi_l = \mO(\epsilon^l), \qquad l = 0, 1, \cdots .
\end{equation}
To prove this, the Chapman-Enskog expansion is utilized, and the specific intensity $\psi$ is expanded in power series of $\epsilon$ as
\begin{equation}
	\psi=\psi^{(0)}+\epsilon\psi^{(1)}+\epsilon^2\psi^{(2)}+\cdots. 
\end{equation}
Define
\begin{equation}
	\psi_l^{(k)}= \int_{-1}^{1} P_l(\mu)\psi^{(k)}(t,x,\mu)\dd\mu, \qquad l = 0, 1, \cdots .
\end{equation}
Then, it holds that  
\begin{equation}
    \label{eq:oa_1}
    \psi_l = \psi_l^{(0)} +\epsilon \psi_l^{(1)} +\epsilon^2 \psi_l^{(2)} + \cdots .
\end{equation}
We will show that for any fixed $k$, 
\begin{equation}
	\psi_l^{(k)} = 0,\qquad \forall l > k.
	\label{eq:claim}
\end{equation}
First, if $k=0$, matching the order of $\mathcal{O}(1)$ in \eqref{eq:pn_sys1} we can obtain
\begin{equation}
	\psi^{(0)}=acT_e^4.
\end{equation}
Then, by the orthogonality of the spherical harmonic functions, it holds that
\begin{equation}
	\psi_l^{(0)}=0,\qquad \forall l \geqslant 1.
\end{equation}
Hence, for the case of $k=0$, \eqref{eq:claim} remains valid.
Subsequently, assuming that \eqref{eq:claim} holds for all $k$ up to and including $n$, as stated
\begin{equation}
	\psi_l^{(k)}=0,\qquad \forall k\leqslant n,\qquad  l>k.
	\label{eq:assum}
\end{equation}
When considering the case where $k=n+1$, it is established that for $\psi_l$ with $l>0$ that 
\begin{equation}
    \sigma \psi_l = \frac{\epsilon^2}{c}\fpp{\psi_l}{t}+  \frac{\epsilon l}{2l+1} \fpp{\psi_{l-1}}{x} + \frac{\epsilon (l+1)}{2l+1} \fpp{\psi_{l+1}}{x}. 
\end{equation}
Matching the order of $\mathcal{O}(\epsilon^{n+1})$ yields
\begin{equation}
	\sigma \psi_l^{(n+1)} = \frac{\epsilon^2}{c}\fpp{\psi_l^{(n-1)}}{t}+  \frac{\epsilon l}{2l+1} \fpp{\psi_{l-1}^{(n)}}{x} + \frac{\epsilon (l+1)}{2l+1} \fpp{\psi_{l+1}^{(n)}}{x}. 
\end{equation}
According to \eqref{eq:assum}, it can be observed that $\psi_l^{(n+1)}=0$ for values of $l$ greater than $n+1$. With this established, we have successfully demonstrated \eqref{eq:claim}, and consequently, \eqref{eq:oa} can be derived naturally.

\section{WENO reconstruction}
\label{app:weno}

The 3-order WENO reconstruction \cite{WENO}:
\begin{equation}
    \begin{aligned}
    &f^{L,1}=\frac{3}{2}f_{j}-\frac{1}{2}f_{j-1},\quad 
    f^{L,2}=\frac{1}{2}f_{j}+\frac{1}{2}f_{j+1},\\ 
    &f^{R,1}=\frac{3}{2}f_{j}-\frac{1}{2}f_{j+1},\quad 
    f^{R,2}=\frac{1}{2}f_{j}+\frac{1}{2}f_{j-1},\\
    &\omega_{L,1}=\frac{1/3}{\left[10^{-6}+(f_{j}-f_{j-1})^2\right]^2},\quad 
    \omega_{L,2}=\frac{2/3}{\left[10^{-6}+(f_{j+1}-f_{j})^2\right]^2},\\
    &\omega_{R,1}=\frac{1/3}{\left[10^{-6}+(f_{j+1}-f_{j})^2\right]^2},\quad 
    \omega_{R,2}=\frac{2/3}{\left[10^{-6}+(f_{j}-f_{j-1})^2\right]^2},\\
    &f_{j+1/2}^L=\frac{\omega_{L,1}f^{L,1}+\omega_{L,2}f^{L,2}}{\omega_{L,1}+\omega_{L,2}},\quad 
    f_{j-1/2}^R=\frac{\omega_{R,1}f^{R,1}+\omega_{R,2}f^{R,2}}{\omega_{R,1}+\omega_{R,2}}.
    \end{aligned}
\end{equation}
Then, the numerical flux can be obtained from $f_{j+1/2}^L$ and $f_{j-1/2}^R$. 



\bibliographystyle{plain}
\bibliography{article}

\end{document}